\documentclass{bhamthesis}
\title{Complementarity and Related Problems}
\author{Lianghai Xiao}
\date{March~2020}

\usepackage{imakeidx} 
\usepackage[pdftex]{graphicx}
\usepackage{ctable}
\usepackage{amssymb, amsmath,amsfonts,amscd, amsxtra,color}
\usepackage{ntheorem}
\usepackage[active]{srcltx}
\usepackage{verbatim}
\usepackage{url}
\usepackage{cite}
\usepackage{enumerate}
\usepackage[margin=2cm,nohead]{geometry}
\usepackage{multicol}
\usepackage{graphics}
\usepackage{threeparttable}
\usepackage{indentfirst}
\usepackage{rotating}
\usepackage{longtable}
\usepackage[acronym]{glossaries}
\usepackage[stdsubgroups,nocfg]{nomencl}
\usepackage{setspace}
\usepackage{etoolbox}
\usepackage{xcolor}
\theoremstyle{break}
\theoremheaderfont{\normalfont\bfseries\hspace{-\theoremindent}}

\newtheorem{definition}{Definition}[section]

\newtheorem{lemma}{Lemma}[section]
\newtheorem{corollary}{Corollary}[section]
\newtheorem{proposition}{Proposition}[section]
\newtheorem{remark}{Remark}[section]
\newtheorem{example}{Example}[section]
\newtheorem{theorem}{Theorem}[section]
\newtheorem*{proof}{ \it Proof:}

\makeatletter

\newcommand\listdefinitionname{List of Definitions}
\newcommand\listofdefinitions{%
  \section*{\listdefinitionname}\@starttoc{def}}
\makeatother

\DeclareMathOperator{\C}{\mathcal C}

\DeclareMathOperator{\CP}{CP}
\DeclareMathOperator{\diag}{diag}
\DeclareMathOperator{\grad}{grad}

\DeclareMathOperator{\ICP}{ICP}
\DeclareMathOperator{\inte}{int}
\DeclareMathOperator{\LCP}{LCP}
\DeclareMathOperator{\MixCP}{MixCP}
\DeclareMathOperator{\MixICP}{MixICP}
\DeclareMathOperator{\NCP}{NCP}

\DeclareMathOperator{\SCP}{SOL-CP}
\DeclareMathOperator{\SICP}{SOL-ICP}
\DeclareMathOperator{\SLCP}{SOL-LCP}
\DeclareMathOperator{\SMixCP}{SOL-MixCP}
\DeclareMathOperator{\SMixICP}{SOL-MixICP}
\DeclareMathOperator{\SNCP}{SOL-NCP}

\DeclareMathOperator{\SStLCP}{SOL-SLCP}
\DeclareMathOperator{\SStMCP}{SOL-SMixCP}
\DeclareMathOperator{\StLCP}{SLCP}
\DeclareMathOperator{\StMCP}{SMixCP}
\DeclareMathOperator{\sgn}{sgn}

\newcommand{\beq}{\begin{equation}}

\newcommand{\bs}{\left(\begin{smallmatrix}}

\newcommand{\E}{\mathbb E}
\newcommand{\eeq}{\end{equation}}

\newcommand{\es}{\end{smallmatrix}\right)}
\newcommand{\f}{\frac}

\newcommand{\lf}{\left}
\newcommand{\lng}{\langle}

\newcommand{\mc}{\mathcal}

\newcommand{\p}{\partial}

\newcommand{\pr}{^\perp}
\newcommand{\R}{\mathbb R}
\newcommand{\rg}{\right}
\newcommand{\rng}{\rangle}
\newcommand{\SP}{\mathbb S}
\newcommand{\tp}{^\top}

\def\argmin{\operatorname{argmin}}

\def\min{\operatorname{min}}
\def\max{\operatorname{max}}

\renewcommand\nomgroup[1]{%
  \item[\bfseries
  \ifstrequal{#1}{S}{Spaces}{%
  \ifstrequal{#1}{M}{Matrices}{%
  \ifstrequal{#1}{Sc}{Scalars}{%
  \ifstrequal{#1}{V}{Vectors}{%
  \ifstrequal{#1}{F}{Functions}{%
  \ifstrequal{#1}{Se}{Sets}{%
  \ifstrequal{#1}{P}{Problem Classes and Fundamental Objects}{%
  \ifstrequal{#1}{Mc}{Matrix Classes}{}}}}}}}}%
]}

\newacronym{kkt}{KKT}{Karush-Kuhn-Tucker}
\newacronym{esoc}{ESOC}{extended second order cone}
\newacronym{esoclcp}{ESOCLCP}{linear complementarity problems on extended second order cones}
\newacronym{cvar}{CVaR}{Conditional Value at Risk}
\newacronym{mv}{MV}{mean variance model}
\newacronym{vi}{VI}{variational inequality}
\newacronym{cp}{CP}{complementarity problem}
\newacronym{lcp}{LCP}{linear complementarity problem}
\newacronym{icp}{ICP}{implicit complementarity problem}
\newacronym{mixcp}{MixCP}{mixed complementarity problem}
\newacronym{mixicp}{MixICP}{mixed implicit complementarity problem}
\newacronym{fb}{FB}{Fischer-Burmeister}
\newacronym{cf}{C-function}{complementarity function}
\newacronym{lm}{LM}{Levenberg-Marquardt}
\newacronym{sesoclcp}{S-ESOCLCP}{the stochastic linear complementarity problem on extended second order cones}
\newacronym{slcp}{SLCP}{stochastic linear complementarity problem}
\newacronym{stcp}{SCP}{stochastic complementarity problem}
\newacronym{ev}{EV}{Expected value}
\newacronym{erm}{ERM}{Expected residual minimisation}
\newacronym{smpec}{SMPEC}{Stochastic mathematical programs with equilibrium constraints}
\newacronym{sp}{SP}{Stochastic programming}
\newacronym{cm}{CM}{CVaR minimisation}
\newacronym{var}{VaR}{Value at risk}
\newacronym{stmcp}{SMixCP}{stochastic mixed complementarity problem}
\newacronym{chks}{CHKS}{Chen-Harker-Kanzow-Smale}
\newacronym{saa}{SAA}{Sample Average Approximation}
\newacronym{aloc}{ALoC}{average loss of complementarity}
\newacronym{men}{MEN}{mean-Euclidean  norm model}
\newacronym{mad}{MAD}{Mean-Absolute Deviation model}
\newacronym{ssd}{SSD}{second degree stochastic dominance}
\newacronym{capm}{CAPM}{Capital Asset Pricing Model}
\newacronym{apt}{APT}{Arbitrage Pricing Model}
\newacronym{cdo}{CDO}{Collateralised Debt Obligation}
\newacronym{ncp}{NCP}{nonlinear complementarity problem}
\newacronym{li}{l.i.}{linearly independent}
\newacronym{as}{a.s.}{almost surely}
\makeglossaries
\makenomenclature
\makeindex

\prefixappendix
\makeatletter\@addtoreset{chapter}{part}\makeatother%

\begin{document}

\frontmatter
\maketitle

\begin{abstract}

{
In this thesis, we present results related to complementarity problems. 

We study the linear complementarity problems on extended second order cones. We convert a linear complementarity problem on an extended second order cone into a mixed complementarity problem on the non-negative orthant. We present algorithms for this problem, and exemplify it by a numerical example. Following this result, we explore the stochastic version of this linear complementarity problem. Finally, we apply complementarity problems on extended second order cones in a portfolio optimisation problem. In this application, we exploit our theoretical results to find an analytical solution to a new portfolio optimisation model. 

We also study the spherical quasi-convexity of quadratic functions on spherically self-dual convex sets. We start this study by exploring the characterisations and conditions for the spherical positive orthant. We present several conditions characterising the spherical quasi-convexity of quadratic functions. Then we generalise the conditions to the spherical quasi-convexity on spherically self-dual convex sets. In particular, we highlight the case of spherical second order cones. 
}

\end{abstract}

\tableofcontents
\listoffigures
\listoftables

\mainmatter

\chapter{Basic Concepts of Complementarity Problems}\label{cht:basic_esoc}

\section{Introduction}
The concept of complementarity is firstly introduced by Karush \cite{karush1939minima} and considered by Dantzig and Cottle in a technical report \cite{dantzig1963positive}, for the non-negative orthant. In 1968, Cottle and Dantzig \cite{cottle1968complementary} connected the linear programming problem, the quadratic programming problem and the bimatrix game problem to the complementarity problem, which attracted many researchers' attentions to this field (see \cite{mangasarian1976linear, garcia1973some, borwein1989linear,alizadeh2003second,FacchineiPang2003}).

As a cross-cutting problem, complementarity problem provides a powerful framework for the study of optimisation and equilibrium problems, and hence has a wide range of applications in engineering and economics. Earlier works in cone complementarity problems present the theory for a general cone and the practical applications merely for the non-negative orthant only (similarly to the books \cite{FacchineiPang2003, MR2503647}). These are related to equilibrium problems in economics, engineering, physics, finance and traffic. Examples in economics are Walrasian price equilibrium models, price oligopoly models, Nash-Cournot production/distribution models, models of invariant capital stock, Markov perfect equilibria, models of decentralised economy and perfect competition equilibrium, models with individual markets of production factors. Engineering and physics applications are frictional contact problems, elastoplastic structural analysis and nonlinear obstacle problems. An example in finance is the discretisation of the differential complementarity formulation of the Black-Scholes models for the American options \cite{jaillet1990variational}. An application to congested traffic networks is the prediction of steady-state traffic flows. In the recent years several applications have emerged where the complementarity problems are defined by cones essentially different from the non-negative orthant such as positive semidefinite cones, second order cones and direct product of these cones (for mixed complementarity problems containing linear subspaces as well). Recent applications of second order cone complementarity problems are in elastoplasticity \cite{MR2925039,MR3010551}, robust game theory \cite{MR2568432,MR2522815} and robotics \cite{MR2377478}. All these applications come from the Karush-Kuhn-Tucker conditions of second order conic optimisation problems.

N\'emeth and Zhang extended the concept of second order cone in \cite{NZ20151} to the \gls{esoc}. Their extension seems the most natural extension of second order cones. Sznajder showed that the extended second order cones in \cite{NZ20151} are irreducible cones (i.e., they cannot be written as a direct product of simpler cones) and calculated the Lyapunov rank of these cones \cite{RS2016}. The applications of second order cones and the elegant way of extending them suggest that the extended second order cones will be important from both theoretical and practical point of view. Although conic optimisation problems with respect to ESOC can be reformulated as conic optimisation problems with respect to second order cones (SOC), we expect that for several such problems using the particular inner structure of the second order cones provides a more efficient way of solving them than solving the transformed conic optimisation problem with respect to second order cones. Indeed, such a particular problem is the projection onto an extended second order cone which is much easier to solve directly than solving the reformulated second order conic optimisation problem \cite{FN2016}.

Until now the extended second order cones of N{\'e}meth and Zhang were used as a working tool only for finding the solutions of mixed complementarity problems on general cones \cite{NZ20151} and variational inequalities for cylinders whose base is a general convex set \cite{NZ2016a}. The applications above for second order cones show the importance of these cones and motivates considering conic optimization and complementarity problems on extended second order cones. In this thesis we develop an application  to portfolio optimisation problems \cite{markowitz1952portfolio,roy1952safety} described in Chapter \ref{cht:psp}. 

We further extend our study to the existance of the solution to a nonlinear complementarity problem. The  existence of the solution to a nonlinear complementarity problem can be converted to a problem of minimising a quadratic function on the intersection between a cone and a sphere, according to \cite[Theorem~18]{Nemeth2006} and \cite[Corollary~8.1]{IsacNemeth2006}. We also study the spherical convexity problem. It is started by exploring the characterisations and conditions for the spherical positive orthant, then it is extended to the spherical quasi-convexity on spherically self-dual convex sets. 

The thesis is organised as follows: In the rest of this chapter, we illustrate the main terminologies and definitions used in this thesis. The terminologies, definitions and basic results of complementarity problem, extended second order cone,and convex sets on the sphere are in Section \ref{sec:cp}, Section \ref{sec:esoc}, and Section \ref{sec:int.1}, respectively. 

In Chapter \ref{cht:lesoc}, we reformulate the linear complementarity problem as a mixed (implicit, mixed implicit) complementarity problem (MixCP) on the non-negative orthant. Our main result is Theorem \ref{thm:main}, which discusses the connections between an linear complementarity problems on extended second order cones (ESOCLCP) and mixed (implicit, mixed implicit) complementarity problems on nonnegative orthant (MixCP). Based on the above, we use some algorithms to solve the MixCP. A solution to this MixCP is equivalent to a solution to the corresponding ESOCLCP. In the last section of this chapter, we provide an example of ESOCLCP corresponding to the cases in Item (iv) of Proposition \ref{prop:cs-esoc}. 

In Chapter \ref{cht:stesoc}, we study the stochastic linear complementarity problems on extended second order cones (stochastic ESOCLCP). We first convert the problem to a stochastic mixed complementarity problem on the nonnegative orthant (SMixCP). Enlightened by the idea of Chen and Lin \cite{chen2011cvar}, we introduce the \gls{cvar} method to measure the loss of complementarity in the stochastic case. A CVaR - based minimisation problem is introduced to achieve a solution which is ``good enough" for the complementarity requirement of the original SMixCP. Smoothing function and sample average approximation methods are introduced and the the problem is converted to a form which can be solved by Levenberg-Marquardt smoothing SAA algorithm. At the end of this chapter, a numerical example will be used to illustrates our results.

In Chapter \ref{cht:psp}, we present an application of extended second order cones to portfolio optimisation  problems. Based on the mean-absolute deviation (MAD) model, we  introduce the mean-Euclidean norm (MEN) model for portfolio optimisation. This new setting has advantages of low computational cost because we work out its analytical solution. 

{
In Chapter \ref{cht:sphere}, we turn to study the spherical convexity as we are motivated by the fact that such questions are related to the existance of the solution of nonlinear complementarity problem. 
In Section~\ref{sec:qcqfcs} we characterise the quadratic spherically quasi-convexity of functions on a general spherically convex set. 
In Section~\ref{sec:qcqfpo} we study the conditions and the properties of spherically quasi-convex quadratic functions defined on the spherical positive orthant. The results of this chapter is published in our paper \cite{FerreiraNemethXiao2018}.
The results in Section \ref{sec:qcqfsdcs} are based on the previous sections. It provides derivations of many useful properties of spherically quasi-convex functions on spherically subdual convex sets. In particular, the spherical positive orthant studied in Section \ref{sec:qcqfpo} is also a spherically subdual convex set. The results obtained on subdual convex sets are applied on the second order cone (Lorentz cone) in Section \ref{sec:qcqflc}. We also prove a condition partially characterising the spherical quasi-convexity of quadratic functions on spherically convex sets associated to the second order cone. 
}

We conclude the study in this thesis by making final remarks in Chapter \ref{cht:FinalRemarks}.

\section{Preliminaries}\label{sec:pre}

In this section, the notations and auxiliary results used throughout this thesis will be presented. Let us start with the definitions of inner products and cones. For $\R^n$, the Euclidian space whose elements are column vectors, the definition of the canonical inner product\index{Canonical inner product} $\langle \cdot, \cdot \rangle$  is given by
\[
	\langle x, y\rangle =\sum\limits_{i=1}^n x_i y_i, \quad x, y \in \R^n,
\]

and the definition of norm\index{Norm} $\|\cdot\|$ is given by
\[
	\|x\| = \sqrt{\langle x, x\rangle}.
\]

Denote by $\R^n_+$ the nonnegative orthant\index{Nonnegative orthant} and by $\R^n_{++}$ the positive orthant \index{Nonnegative orthant!Positive orthant}, that is,
\[
	\R^n_+=\{x=(x_1,\dots,x_n)\tp:x_1\ge0,\dots,x_n\ge0\},
\]
and
\[
	\R^n_{++}=\{x=(x_1,\dots,x_n)\tp:x_1>0,\dots,x_n>0\}.
\]

{
Denote by ${\cal L}$ the second order cone\index{Cone!Second order cone, Lorentz cone} (Lorentz cone) 
\begin{equation}\label{eq:soc}
	{\cal L} := \lf\{x= ( x_1,\dots,x_n)\tp\in\R^n :x_1\geq \sqrt{x_2^2 + \dots + x_n^2}\rg\}.
\end{equation}
}

It should be noted that the Lorentz cone ${\cal L}$, the nonnegative orthant and the positive orthant are self-dual cones. 

Let $k$, $l$ be positive integers. The inner product of pairs of {vectors $\bs x\\u\es,~ \bs y\\v\es \in \R^k\times\R^{\ell}$}, where $x\in\R^k$ and $u\in\R^{\ell}$, is defined by
\[
	\langle\bs x\\u\es ,\bs y\\v\es \rangle =\langle x,y\rangle + \langle u,v\rangle .
\]

Let $\R^n$ be a Euclidian space.  A set ${\cal K} \subseteq \R^n$ is called a \emph{convex cone}\index{Cone!Convex cone} if for any $\alpha , \beta > 0$, and $x, y\in {\cal K}$, we have
\[
	\alpha x + \beta y \in {\cal K} .
\]
In other words, a convex cone is a set which is invariant under multiplication of vectors with positive scalars and addition of vectors. {The {\it dual cone}\index{Cone!Dual cone} of  cone ${\cal K}  \subseteq \R^n$ is the convex cone ${\cal K} ^* \!\!:=\!\{ x\in \R^n : \langle x, y \rangle\!\geq\! 0, ~ \forall \, y\!\in\! {\cal K} \}.$ }  A convex cone \( {\cal K}  \subseteq \R^{n}\)  is called {\it pointed}\index{Cone!Pointed cone} if \({\cal K}  \cap \{-{\cal K} \} \subseteq \{0\}\), or equivalently, if \({\cal K} \) does not contain straight  lines through the origin. A convex cone\index{Cone!Closed cone} which is a closed set is called a \emph{closed convex cone}\index{Cone!Closed convex cone}. Any pointed closed convex cone with nonempty interior will be called {\it proper cone}\index{Cone!Proper cone}. The cone ${\cal K} $ is called {\it subdual}\index{Cone!Subdual cone} if ${\cal K} \subseteq {\cal K} ^*$, {\it superdual}\index{Cone!Superdual cone} if ${\cal K} ^*\subseteq {\cal K} $, and {\it self-dual}\index{Cone!Self-dual cone} if ${\cal K} ^*={\cal K} $.

The matrix $I_n$ denotes the $n\times n$ identity matrix\index{Identity matrix}. We denote by $ \R^{k\times\ell}$ the set of matrices with $k$ rows and $\ell$ columns with real elements. In particular $\R^k$ can be identified with $\R^{k\times 1}$.

Let $x\in\R^n$, then the {\it projection} ${\rm P}_{\cal K}(x)$ of the point $x$ onto the cone ${\cal K}$ is defined by
\[
	{\rm P}_{\cal K}(x):= \argmin_{y} \{\|x-y\|:y\in{\cal K}\}.
\]

For any $x\in {\cal K}$, we define the nonnegative part of $x$, nonpositive part of $x$ and the absolute value of $x$ with respect to ${\cal K}$ by
\begin{equation} \label{eq:npav}
	x_+^{\cal K}:={\rm P}_{\cal K}(x), \qquad x_-^{\cal K}:={\rm P}_{{\cal K}^*}(-x), \qquad  |x|^{\cal K}:=x_+^{\cal K}+x_-^{\cal K},
\end{equation}
respectively. We recall   from   Moreau's  decomposition\index{Moreau's  decomposition} theorem \cite{Moreau1962} (see also  \cite[Theorem 3.2.5]{HiriartLemarecal1}),   that
for a closed convex cone ${\cal K}$  there hold:
\begin{equation} \label{eq:cmdt}
	x=x_+^{\cal K}-x_-^{\cal K}, \qquad \left\langle x_+^{\cal K}, x_-^{\cal K} \right\rangle=0,   \qquad \qquad  ~x\in\R^n.
\end{equation}
For any $z\in \mathbb{R}\times
{\mathbb{R}}^{n-1}$, let{
\(z:=(z_1, {z^{(2)}})  \in \mathbb{R}\times {\mathbb{R}}^{n-1} \)}, where {${z^{(2)}}:= (z_2, z_3, \dots, z_n)^{\tp}$}.
An explicit formula for the projection\index{Projection, projection mapping} mapping    \({\rm P}_{\cal{L}}\)  onto the  Lorentz cone  ${\cal L}$ is given in  \cite[Proposition~3.3]{FukushimaTseng2002}, which is recalled for the case when $x\notin{\cal L}\cup -{\cal L}$ in the following lemma.
\begin{lemma} \label{l:projude-}
Let {\(x=(x_1, {x^{(2)}})  \in  \{(y_1, {y^{(2)}}) \in \R \times \R^{n-1}: ~ |y_1|<\|{y^{(2)}}\|\}  \)} and \({\cal{L}}\) be the Lorentz cone. Then,
{
\[
x_+^{\cal L}= \lf(\frac{x_1+\|{x^{(2)}}\|}{2\|x^{(2)}\|}\rg)\left(\|x^{(2)}\|,x^{(2)}\right), \qquad
x_-^{\cal L}= \lf(\frac{-x_1+\|{x^{(2)}}\|}{2\|x^{(2)}\|}\rg)\left(\|x^{(2)}\|,-x^{(2)}\right)
\]}
and, as a  consequence, the absolute value of $x$ with respect to ${\cal L}$ is given by
{
\[
	|x|^{\cal L} = \frac1{\|x^{(2)}\|}\left(\|{x^{(2)}}\|^2,x_1x^{(2)}\right).
\]}
\end{lemma}
For a general nonzero vector {\(x=\lf(x_1, {x^{(2)}}\rg)  \in  \R \times \R^{n-1}\)} the absolute value of $x$ with respect to ${\cal L}$ is given in
the next
lemma, which follows immediately from Lemma \ref{l:projude-} and equations \eqref{eq:cmdt}.
\begin{lemma} \label{l:projude}
Consider a nonzero vector {\(x=\lf(x_1, {x^{(2)}}\rg)  \in  \R \times \R^{n-1}\)} and let \({\cal{L}}\) be the Lorentz cone. Then, the absolute value of $x$  is given by
{
\[|x|^{\cal L} = \frac1{\|x^{(2)}\|}\Big{(}\max\left(|x_1|,\|{x^{(2)}}\|\right)\|x^{(2)}\|,~\min(|x_1|,\|x^{(2)}\|)\sgn(x_1)x^{(2)}\Big{)},\]}
where  $\sgn(x_1)$ is equal to $-1$, $0$ or $1$ whenever  $x_1$ is negative, zero or positive, respectively.  
\end{lemma}

{
\bigskip \begin{definition}[$S_0$ matrix]\label{def:s0_mat} 
	A matrix $A \in \R^{n\times n}$ is said to be an {\it $S_0$ matrix} if there exists a vector $x\in \R^n_+$ such that 
	\[
		A x \ge 0.
	\]
\end{definition}

\bigskip \begin{definition}[$P_0$ matrix]\label{def:p_mat}\emph{\cite[Definition 2.2]{yamashita1997modified}}
	A matrix $A \in \R^{n\times n}$ is said to be a {\it $P_0$ matrix} if every principal minor of $A$ is non-negative. In particular, if every principal minor of $A$ is positive, $A$ is said to be a {\it $P$ matrix}.
\end{definition}
}

\bigskip \begin{definition}[Schur complement]\label{def:schur_comp}\emph{\cite{zhang2006schur}}
	The {\it Schur complement} for a matrix $M =\left(\begin{smallmatrix} A & B\\C & D\end{smallmatrix}\right)$ in nonsingular matrix $D$ is
	\[
		\left( M/D \right) = A - BD^{-1}C.
	\]
\end{definition}
\index{Schur complement}

\bigskip
In this study, both smooth function and semi-smooth function will be carefully reviewed. We will introduce relevant concepts about both continuously differentiable and Lipschitz continuous.

\bigskip \begin{definition}[Fr{\'e}chet differentiable]\label{def:Frechetdif}
	Let ${\cal K}$ be an open subset with ${\cal K} \subseteq \R^{\ell}$ and $ f: \R^{\ell} \supseteq {\cal K} \rightarrow \R^k$, $k$ is not necessary different from $\ell$. We say that $f$ is a {\it differentiable function} on ${\cal K}$, if {there is a  linear map} $ J: {\cal K} \rightarrow \R^k$ such that
	\[
		\lim_{t\rightarrow 0}\frac{\|f(x+t)-f(x)-J(t)\|}{\|t\|} = 0,
	\]
	for any $x\in {\cal K}$.
\end{definition}
\index{Fr{\'e}chet differentiable}

\bigskip \begin{definition}[Continuously differentiable function]
	Let ${\cal K}$ be an open subset with $ {\cal K} \subseteq \R^{\ell}$ and $ f: \R^{\ell} \supseteq {\cal K} \rightarrow \R^k$, $k$ is not necessary different from $\ell$. We say that $f$ is a {\it continuously differentiable function} on ${\cal K}$, if there is a linear map $ J: {\cal K} \rightarrow \R^k$
	\[
		\lim_{t\rightarrow 0}\frac{\|f(x+t)-f(x)-J(t)\|}{\|t\|} = 0,
	\]
	such that, for any $x\in {\cal K}$ the map $J$ is continuous.
\end{definition}
\index{Continuously differentiable function}

\bigskip \begin{definition}[{Lipschitz continuous function}]\label{def:lips_func}\cite[Definition 4.6.2]{sohrab2003basic}
	\item[(i)]Let ${\cal K}$ be an open subset with $ {\cal K} \subseteq \R^{\ell}$ and $ f: \R^{\ell} \supseteq {\cal K} \rightarrow \R^k$, $k$ is not necessary different from $\ell$. We say that $f$ is a {\it Lipschitz continuous function} on ${\cal K}$, if there is a constant $ \lambda >0$ such that
	\begin{equation}
		\| f(x) - f(x') \| \le \lambda\|x - x'\| \quad \forall x, x' \in {\cal K}
	\end{equation}
	\item[(ii)] We say that $f$ is {locally Lipschitz continuous} if for any $ x \in {\cal K}$, there exists $ \epsilon > 0 $ such that $f$ is Lipschitz on $ {\cal K}\cap \bar{\textbf{B}}(x,\epsilon)$, where $\bar{\textbf{B}}(x,\epsilon):= \left\{y\in \R^k : \|x - y\|\leq \epsilon\right\}$ is the closed ball centered at x.
\end{definition}
\index{Lipschitz function}

\bigskip \begin{definition}[Semismooth function]\label{def:semi_smu}\cite[Definition 1]{mifflin1977semismooth}
	A function $f(x)$ is {\it semismooth} at $x\in\R^n$ if
	\item[(i)] $f(x)$ is a Lipschitz function on $\bar{\textbf{B}}(x,\epsilon)$, the closed ball centered at x, and
	\item[(ii)] for each $d\in\R^n$ and for any sequences $\{t_m\}\subseteq \R_+$, $\{\epsilon_m\}\subseteq\R^n$ and $\{g_m\}\subseteq\R^n$ such that $\{t_m\}\downarrow 0$, $\{\frac{\epsilon_m}{t_m}\}\rightarrow 0\in\R^n$ and {$g_m\in \partial f(x+t_m d+\epsilon_m)$}, the sequence $\{\langle g_m,d\rangle \}$ has exactly one accumulation point.
\end{definition}
\index{Semismooth function}

%
%

\section{Complementarity problems}\label{sec:cp}

\
This section briefly overviews the terminologies and definitions of complementarity problems before needed in following chapters. Some basic results about complementarity problems will be presented.

\bigskip \begin{definition}[Complementarity set]\label{def:cs}
	 Let ${\cal K}\subseteq\R^m$ be a nonempty closed convex cone and ${\cal K}^*$ its dual. The set \textcolor[rgb]{1.00,1.00,1.00}{$\C({\cal K})$}
	\[
		\C({\cal K}):=\left\{(x,y)\in {\cal K}\times {\cal K}^*:\langle x, y\rangle = 0\right\}
	\]
	is called the \emph{complementarity set} of cone ${\cal K}$.
\end{definition}
\index{Complementarity!Complementarity set}

\bigskip \begin{definition}[Complementarity function]\label{def:cfunc}
	A function $\phi(a,b)$ is called \emph{\gls{cf}} if it satisfies:
	\[
		\phi(a,b) = 0 \quad \Leftrightarrow\quad a\ge 0,\quad b\ge 0,\quad ab = 0.
	\]
\end{definition}
\index{Complementarity!Complementarity function}

\bigskip \begin{definition}[Variational inequalities]\label{def:vi}
	Let ${\cal K}\subseteq\R^n$ be a nonempty closed convex cone and $F:\R^n\to\R^n$ be a mapping. The \emph{\gls{vi}} defined by $F$ and ${\cal K}$ is the problem
	\begin{equation}\label{eq:vi}
		VI(F,{\cal K})\left\{
		\begin{array}{l}
			Find \; x\in \R^n, \; such \; that\\
			\langle y - x, F(x)\rangle\ge 0, \quad \forall y\in {\cal K}.
		\end{array}
		\right.
	\end{equation}
\end{definition}
\index{Complementarity!Variational inequalities}

\bigskip \begin{definition}[Complementarity problem]\label{def:cp}
	Let $F:\R^n\to\R^n$ be a mapping. Let ${\cal K}\subseteq\R^n$ be a nonempty closed convex cone and ${\cal K}^*$ its dual.
	Defined by ${\cal K}$ and $F$ the \emph{\gls{cp}} is:
	\begin{equation}\label{eq:cp}
		CP(F,{\cal K})\left\{
		\begin{array}{l}
			Find \; x\in \R^n, \; such \; that\\
			\left(x, F(x)\right)\in\C({\cal K}).
		\end{array}
		\right.
	\end{equation}
	The solution set of $\CP(F,{\cal K})$ is denoted by $\SCP(F,{\cal K})$:
	\[
		\SCP(F,{\cal K}) =  \{x\in \R^n: (x, F(x)) \in \C({\cal K})\}.
	\]
\end{definition}
\index{Complementarity problem}

\bigskip
In particular, the definition of the \emph{\gls{lcp}} is:
\index{Complementarity problem!Linear complementarity problem}
\[
	LCP(F,{\cal K})\left\{
	\begin{array}{l}
		Find \; x\in \R^n, \; such \; that\\
		\left(x, F(x)\right)\in\C({\cal K}).
	\end{array}
	\right.
\]
where $F(x)$ defined by $F(x)=Tx+r$ is a linear function, where $T\in\R^{n\times n}$ is a matrix and $r\in\R^n$. The solution set of $\LCP(T,r,{\cal K})$ is denoted by $\SLCP(T,r,{\cal K})$.


\bigskip \begin{definition}[Implicit complementarity problem]\label{def:icp}
	Let $G,F:\R^n\to\R^n$ be mappings. The {\it \gls{icp}} defined by $G$, $F$, and the cone ${\cal K}$ is:
	\[
		\ICP(G,F,{\cal K})\left\{
		\begin{array}{l}
			Find \; x\in\R^n, \; such \; that\\
			\left(G(x), F(x)\right)\in\C({\cal K}).
		\end{array}
		\right.
	\]
	The solution set of $\ICP(G,F,{\cal K})$ is denoted by $\SICP(G,F,{\cal K})$:
	\[
		\SICP(G,F,{\cal K}) = \{x\in \R^n:  (G(x), F(x)) \in \C({\cal K})\}.
	\]
\end{definition}
\index{Complementarity problem!Implicit complementarity problem}

\bigskip \begin{definition}[Mixed complementarity problem]\label{def:mixcp}
	Consider the mappings $F_1:\R^k\times\R^\ell\to \R^k$ and $F_2:\R^k\times\R^\ell\to \R^\ell$. Let ${\cal S} \in \R^k$ be a nonempty closed convex cone. The {\it \gls{mixcp}} defined by $F_1$, $F_2$, and ${\cal S}$ is:
	\[
		\MixCP(F_1,F_2,{\cal S}):\left\{
		\begin{array}{l}
			Find \; \bs x\\u\es\in\R^k\times\R^{\ell}, \; such \; that\\
			F_2(x,u)=0, \; and\; (x,F_1(x,u))\in\C({\cal S}).
		\end{array}
		\right.
	\]
	The solution set of $\MixCP(F_1,F_2,{\cal S})$ is denoted by $\SMixCP(F_1,F_2,{\cal S})$:
	\[
		\SMixCP(F_1,F_2,{\cal S}) =\{\bs x\\u\es\in\R^k\times\R^{\ell}: F_2(x,u)=0,(x, F_1(x,u)) \in \C({\cal S})\}.
	\]
\end{definition}
\index{Complementarity problem!Mixed complementarity problem}

\bigskip
The mixed complementarity problem is one of the most important problem formulations in mathematical programming. Many well-studied optimisation problems can be converted into a mixed complementarity problem.

\bigskip \begin{definition}[Mixed implicit complementarity problem]\label{def:mixicp}
	Consider the mappings $F_1$, $G_1:\R^k\times\R^\ell\to \R^k$, $F_2:\R^k\times\R^\ell\to \R^{\ell}$, and a proper cone ${\cal S} \in \R^k$ . The {\it \gls{mixicp}} defined by $F_1$, $F_2$, $G_1$, and ${\cal S}$ is
	\[
		\MixICP(F_1,F_2,G_1,{\cal S}):\left\{
		\begin{array}{l}
			Find \; \bs x\\u\es\in\R^k\times\R^{\ell}, \; such \; that\\
			F_2(x,u)= 0, \; and\; (G_1(x,u),F_1(x,u))\in\C({\cal S}).
		\end{array}
		\right.
	\]
	The solution set of $\MixICP(F_1,F_2,G_1,{\cal S})$ is denoted by $\SMixICP(F_1,F_2,G_1,{\cal S})$:
	\begin{align*}
		\SMixICP&( F_1,F_2,G_1,{\cal S})  = \\
		&\{\bs x\\u\es\in\R^k\times\R^{\ell}: F_2(x,u)= 0,(G_1(x,u), F_1(x,u)) \in \C({\cal S})\}.
	\end{align*}
\end{definition}
\index{Complementarity problem!Mixed implicit complementarity problem}

Based on the definitions above, we get the following propositions straightforwardly. 

\bigskip \begin{proposition}\label{prop:cp_mixcp}
	Let $n,k,\ell$ be nonnegative integers such that $n=k+l$, ${\cal S}\in\R^k$ be a nonempty closed convex cone and
${\cal K}={\cal S}\times\R^\ell$. Denote by ${\cal S}^*$ the dual of ${\cal S}$ in $\R^k$ and by ${\cal K}^*\subset \R^k\times\R^\ell$ the dual of ${\cal K}$ in $\R^k\times\R^\ell$. Consider the mappings \(F_1:\R^k\times\R^\ell\to \R^k,\) \(F_2:\R^k\times\R^\ell\to \R^\ell.\) Define the
	mappings
	\(F:\R^k\times\R^\ell\to\R^k\times\R^\ell\) by \(F(x,u)=\bs F_1(x,u) \\ F_2(x,u)\es\). Then,
	\[\bs x\\u\es\in\SCP(F,{\cal K})\iff \bs x\\u\es\in\SMixCP(F_1,F_2,{\cal S}).\]
\end{proposition}

\begin{proof}
	It is easy to check that ${\cal K}^*={\cal S}^*\times\{0\}$. The result follows immediately from ${\cal K}^*={\cal S}^*\times \{0\}$ and the definitions of $CP(F,{\cal K})$ and
	$\MixCP(F_1,F_2,{\cal S})$.
\hfill$\square$\par \end{proof}

{
\bigskip \begin{proposition}\label{prop:icp_mixicp}
	Let $n,k,\ell$ be nonnegative integers such that $n=k+l$, ${\cal S}\in\R^k$ be a nonempty closed convex cone and
${\cal K}={\cal S}\times\R^\ell$. Denote by ${\cal S}^*$ the dual of ${\cal S}$ in $\R^k$ and by ${\cal K}^*\subset \R^k\times\R^\ell$ the dual of ${\cal K}$. Consider the mappings \(F_1,G_1:\R^k\times\R^\ell\to \R^k,\) \(F_2,G_2:\R^k\times\R^\ell\to \R^\ell.\) Define the mappings
	\(F,G:\R^k\times\R^\ell\to\R^k\times\R^\ell\) by \(F(x,u)=(F_1(x,u),F_2(x,u)),\) \(G(x,u)=(G_1(x,u),G_2(x,u)),\) respectively.
	Then,
	\[\bs x\\u\es\in\SICP(F,G,{\cal K})\iff \bs x\\u\es\in\SMixICP(F_1,F_2,G_1,{\cal S}).\]
\end{proposition}
}

\begin{proof}
	The result follows immediately from ${\cal K}^*={\cal S}^*\times \{0\}$ and the definitions of $\ICP(F,G,{\cal K})$ and $\MixICP(F_1,F_2,G_1,{\cal S})$.
\hfill$\square$\par \end{proof}

\section{Extended second order cone}\label{sec:esoc}

The extended second order cones (ESOC), introduced by N\'emeth and Zhang\cite{NZ20151}, are natural extensions of the second order cones (or Lorentz cones). The researches about ESOCs can be found in \cite{NZ2016a,kong2017isotonicity,FN2016,RS2016}. The particular structure of ESOCs {provides} a more direct method for solving these problems without reformulations, such important examples are mixed complementarity problems on general cones \cite{NZ20151} and variational inequalities on cylinders\cite{NZ2016a, kong2017isotonicity}. Moreover, from purely conic analysis point of view, the ESOCs cannot be trivially reduced to second order cones because the ESOCs are irreducible\cite{RS2016} (i.e., they cannot be written as a direct product of simpler cones). Employing the forerunners' results associated to ESOCs, this section offers some fundamental definitions and properties before needed in the later discussion of linear complementarity problems on ESOCs.

Let $n,k,\ell$ be nonnegative integers such that $n=k+l$. The definitions of the mutually dual extended second order cone $L(k,\ell)$ and $M(k,\ell)$ in $\mathbb{R}^n\equiv\R^k\times\R^\ell$ are:
\index{Extended second order cone}
{
\begin{equation}\label{esoc}
	L(k,\ell) = \{\bs x\\u\es \in \mathbb{R}^k\times \mathbb{R}^\ell : x \geq \|u\|e\},
\end{equation}
}
\begin{equation}\label{desoc}
	M(k,\ell) = \{\bs x\\u\es \in \mathbb{R}^k\times \mathbb{R}^\ell : e\tp x\geq \|u\|, \;x \ge 0 \},
\end{equation}
where  $e=(1, \dots, 1)\tp \in \mathbb{R}^k $. If there is no ambiguity about the dimensions, then we simply denote $L(k,\ell)$ and $M(k,\ell)$ by $L$ and $M$, respectively.

\begin{figure}[bt!]
	\centering
	\begin{minipage}{1\textwidth}
	\includegraphics[width=1\textwidth, height=0.5\textwidth]{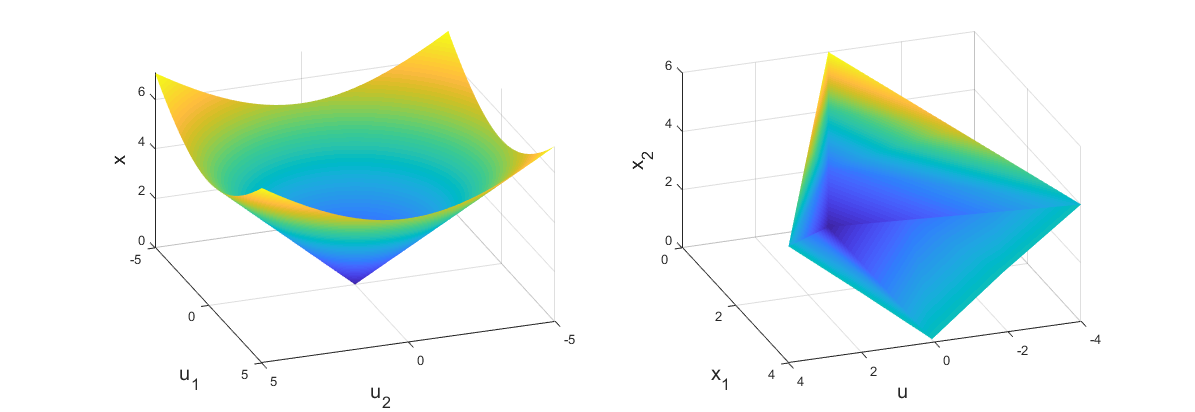}\\
	{\footnotesize Note: This figure provides the comparison of second order cone (SOC) : $C = \{ (x, u)\tp \in\R^{1+ 2}: x\ge \|u\|\}$ and extended second order cone (ESOC): $C = \{(x, u)\tp \in\R^{2+ 1}: x\ge \|u\|e\}$. It shows that the SOC is symmetric, whereas the ESOC is not symmetric. Both cones are in $\R^3$, but the shapes of them are very different. \par}
	\end{minipage}
	\caption{Second order cone and extended second order cone}
	\label{fig:soc}
\end{figure}

We remark that when $k = 1$ the ESOC is equivalent to a second order cone (defined in \eqref{eq:soc} in $\R\times \R^{\ell}$. Hence in the following study we assume that the integer $k\ge 2$. Figure \ref{fig:soc} provides the comparison of three-dimensional second order cone and three-dimensional extended second order cone. Both cones can be expressed as a pair of two vectors. Though both cones are in $\R^3$, the shapes of them are very different. This is because that the second order cone is in $ \R^{1}\times\R^{ 2}$, whilst the extended second order cone is in $ \R^{2}\times\R^{1}$.

\bigskip \begin{proposition}
	The extended second order cone is a pointed, closed convex cone with nonempty interior (henceforth it is a proper cone). 
\end{proposition}

\bigskip \begin{proposition}\label{prop:selfdual}
	Given any nonnegative integers $k$ and $\ell$ with $k\ge 2$, we have $L \subset M$, and $M \not\subseteq L$.
\end{proposition}

Proposition \ref{prop:selfdual} shows that an ESOC is subdual but not superdual. Hence, the ESOC is not self-dual.

\bigskip \begin{proposition}\label{prop:cs-esoc}
	Let $x,y\in\R^k$ and $u,v\in\R^\ell\setminus\{0\}$.
	\begin{enumerate}
		\item[(i)] {$(x,0,y,0):=(\bs x\\0\es,\bs y\\0\es)\in\C(L)$ if and only if $(x,y)\in\C(\R^k_+)$.}
		\item[(ii)] $(x,0,y,v)\in\C(L)$ if and only if $e\tp y\ge\|v\|$ and $(x,y)\in\C(\R^k_+)$.
		\item[(iii)] $(x,u,y,0)\in\C(L)$ if and only if $x\ge\|u\|e$ and $(x,y)\in\C(\R^k_+)$.
		\item[(iv)] $(x,u,y,v)\in\C(L)$ if and only if there exists $\lambda>0$ such that $v=-\lambda u$,
			$e\tp y=\|v\|$ and $(x-\|u\|e,y)\in C(\R^k_+)$.
	\end{enumerate}
\end{proposition}

\begin{proof}
	Item (i) follows definition \eqref{esoc} and \eqref{desoc}. Let $(x,0,y,0)\in\C(L)$, it is trivial to have that $x \ge 0$, $y\ge 0$ and $\langle x,y \rangle =0 $, i.e., $(x,y)\in\C(\R^k_+)$.

	Item (ii) follows definition \eqref{desoc}. Let $(x,0,y,v)\in\C(L)$, naturally we have $\langle \bs x\\0\es , \bs y\\v\es\rangle = \langle x,y \rangle + \langle 0, v \rangle = \langle x, y\rangle = 0$, as well as $ e\tp y\ge\|v\|$ because $(y,v)\in M$, so that we conclude $(x,y)\in\C(\R^k_+)$.

	Item (iii) follows definition \eqref{esoc}. Let $(x,u,y,0)\in\C(L)$, an easy consequence is $x\ge\|u\|e$; and $\langle \bs x\\u\es, \bs y\\0\es\rangle = \langle x, y\rangle = 0$, then $(x,y)\in\C(\R^k_+)$.

	Item (iv)follows from \cite[Proposition 1]{FN2016}. For the completeness of the results we will
	repeat its proof here. First assume that there exists $\lambda>0$ such that $v=-\lambda u$,
	$e\tp y=\|v\|$ and $(x-\|u\|e,y)\in C(\R^k_+)$. Thus, $\bs x\\u\es\in L$ and $\bs y\\v\es\in M$. On the other
	hand, \[\lng \bs x\\u\es,\bs y\\v\es\rng=x\tp y+u\tp v=\|u\|e\tp y-\lambda\|u\|^2=\|u\|\|v\|-\lambda\|u\|^2=0.\]
	Thus, $(x,u,y,v)\in C(L)$.

	Conversely, if $(x,u,y,v)\in C(L)$, then $\bs x\\u\es\in L$, $\bs y\\v\es\in M$ and
	\[0=\lng \bs x\\u\es,\bs y\\v\es\rng=x\tp y+u\tp v\ge\|u\|e\tp y+
	u\tp v\ge\|u\|\|v\|+u\tp v\ge0.\]
	This implies the existence of a $\lambda>0$ such that $v=-\lambda u$, $e\tp y=\|v\|$ and
	$(x-\|u\|e)\tp y=0$. It follows that $(x-\|u\|e,y)\in C(\R^k_+)$.
\hfill$\square$\par \end{proof}

The following corollary generalises all the cases shown in Proposition \ref{prop:cs-esoc}:

\bigskip \begin{corollary}
	Let $x$, $y\in \R^k$, and $u$, $v\in\R^{\ell}$. Then, $(x, u, y, v)\in\C (L)$ if and only if  there exists $\lambda >0$ such that
	\begin{enumerate}
	  \item[(1)] $\|u\|\|v\|(v+\lambda u) = 0$,
	  \item[(2)] $\|u\|(x - \|u\|e)\ge 0$,
	  \item[(3)] $(\|u\|\|v\|, e\tp y- \|v\|)\in\C(\R^2_+)$,
	  \item[(4)] $\lf[1 - sgn(\|u\|\|v\|)\rg](x,y)\in C(\R^k_+)$, and
	  \item[(5)]$\|v\|(x-\|u\|e,y)\in \C(\R^k_+)$.
	\end{enumerate}
\end{corollary}

\begin{proof}
	Given the four cases in Proposition \ref{prop:cs-esoc}, we will examine each item in the corollary accordingly.

	If $u=v=0$, then the group of items (1)-(5) is equivalent to the group of items (3)-(4), because items (1)-(3) and (5) trivially hold. In turn
	the group of items (3)-(4) is equivalent to $(x,y)\in\C(\R^k_+)$. Hence, the result in this case follows from Proposition \ref{prop:cs-esoc}
	item (i).

	If $u=0$ and $v\ne 0$, then the group of items (1)-(5) is equivalent to the group of items (3)-(5), because items (1) and (2) trivially hold.
	In turn the group of items (3)-(5) is equivalent to $e\tp y\ge\|v\|$ and $(x,y)\in\C(\R^k_+)$. Hence, the result in this case follows from
	Proposition \ref{prop:cs-esoc} item (ii).

	If $u\ne0$ and $v=0$, then the group of items (1)-(5) is equivalent to the group of items (2)-(4), because items (1) and (5) trivially hold.
	In turn the group of items (2)-(4) is equivalent to $x\ge\|u\|e$ and $(x,y)\in\C(\R^k_+)$. Hence, the result in this case follows from
	Proposition \ref{prop:cs-esoc} item (iii).

	If $u\ne0$ and $v\ne 0$, then item (1) is equivalent to $v=-\lambda u$. In turn item (3) is equivalent to $e\tp y = \|v\|$, the group of items
	$\{$(2), (5)$\}$ is equivalent to $(x-\|u\|e,y)\in \C(\R^k_+)$ and item (4) trivially holds. In conclusion the group of items (1)-(5) is
	equivalent to $v=-\lambda u$ for some $\lambda$, $e\tp y = \|v\|$ and $(x-\|u\|e,y)\in \C(\R^k_+)$. Hence, the result in this case follows from
	Proposition \ref{prop:cs-esoc} item (iv).

\hfill$\square$\par \end{proof}

\section{Convex sets on the sphere} \label{sec:int.1}

\
This section gives the definitions about the convex sets on the sphere. Some results in this chapter are based on the results in \cite{Nemeth1998}, but we provide more explicit statements and proofs herein.
 We start this section with the definition of copositive matrix and Z-matrix. 

\bigskip \begin{definition}[${\cal K}$-Copositive matrix]\label{def:k_copo}
	A matrix $A \in \R^{n\times n}$ is {\it ${\cal K}$-copositive} if
	\[
		\langle Ax, x \rangle\!\geq\! 0
	\]
	for any $x\in {\cal K}$.
\end{definition}\index{Copositive!${\cal K}$-Copositive}

Particularly, we give the following definition: 

\bigskip \begin{definition}[Copositive matrix]\label{def:copo}
	A matrix $A \in \R^{n\times n}$ is {\it copositive} if
	\[
		\langle Ax, x \rangle\!\geq\! 0
	\]
	for any $x\in \R^n_+$.
\end{definition}\index{Copositive}

According to the two definitions above, the Definition \ref{def:copo} is equivalent to Definition \ref{def:k_copo} when ${\cal K} = \R^n_+$. 

Let \( {\cal K} \subseteq \R^{n}\)  be a (not necessarily convex) cone. The Lorentz cone \({\cal{L}}\) can be written as\index{Second order cone, Lorentz cone}
\begin{equation*} \label{eq:LorentzConeJ}
	{\cal L}:=\left\{x=(x_1,\dots,x_n)^\top\in\R^n:~x_1\geq 0,\textrm{ }\lng Jx,x\rng\ge0\right\},
\end{equation*}
where $J=\diag(1,-1,\dots,-1)\in\R^{n\times n}$.
It is easy to see that 
\[
	{\cal L}\cup-{\cal L}=\left\{x=(x_1,\dots,x_n)^\top\in\R^n:\lng Jx,x\rng\ge0\right\}.
\]

	This straightforwardly implies that  $A\in \R^{n \times n}$ is ${\cal L}$-copositive if and only if it is ${\cal L}\cup-{\cal L}$-copositive.
Hence, the S-Lemma (see \cite{yakubovich1997s,polik2007survey}) implies:
\begin{lemma}\label{lorcop}
	$A\in \R^{n \times n}$ is ${\cal L}$-copositive if and only if there exist a $\rho\ge0$ such that $A-\rho J$ is positive
	semidefinite.
\end{lemma}

The matrix $I_n$ denotes the $n\times n$ identity matrix\index{Identity matrix}. We denote by $ \R^{k\times\ell}$ the set of matrices with $k$ rows and $\ell$ columns with real elements. In particular $\R^k$ can be identified with $\R^{k\times 1}$.

Recall that  $A=(a_{ij})\in \R^{n \times n}$ is {\it positive} if  $a_{ij}>0$ and {\it nonnegative}  if $a_{ij}\geq 0$ for any $i,j=1, \ldots, n$.  A matrix   $A \in \R^{n \times n}$ is {\it reducible} if there is permutation matrix  $P \in \R^{n \times n}$ such that
\[
	P^{T}AP=
	\begin{bmatrix}
	B_{11} & B_{12} \\
	0 & B_{22}
	\end{bmatrix},
\]
\[
	B_{11}\in \R^{m \times m}, ~B_{22}\in \R^{(n-m) \times (n-m)}, ~ B_{12}\in  \R^{m \times (n-m)},  \quad m<n.
\]
A matrix   $A \in \R^{n \times n}$ is {\it irreducible} if  it not reducible.  In the following  we state a version of   {\it Perron-Frobenius
theorem} for both positive matrices and nonnegative irreducible matrices,  its proof  can be found in \cite[Theorem~8.2.11]{HornJohnson85} and  \cite[Theorem~8.4.4]{HornJohnson85}, respectively.

\bigskip \begin{theorem}[Perron-Frobenius Theorem] \label{Perron-Frobenius theorem1}
Let $A\in \R^{n \times n}$  be  either nonnegative and  irreducible or positive. Then  $A$ has a dominant eigenvalue  $\lambda_{max}(A)\in \R$ with  associated eigenvector $v\in \R^n$ which satisfies the following properties:
\begin{itemize}
\item[i)] The  eigenvalue $\lambda_{max}(A)>0$ and its  associated eigenvector $v\in \R^n_{++}$;
\item[ii)] The eigenvalue $\lambda_{max}(A)$ has multiplicity one;
\item[iii)] Every other  eigenvalue  $\lambda$ of $A$ is less that $\lambda_{max}(A)$ in absolute value, i.e, $|\lambda|<\lambda_{max}(A)$;
\item[iii)] There are no other positive  or  non-negative eigenvectors  of $A$ except positive multiples of $v$.
\end{itemize}
\end{theorem}

\bigskip \begin{definition}[Z-matrix]\label{def:z_m}
	A matrix $A \in \R^{n\times n}$ is a Z-matrix if its off-diagonal elements are all nonpositive.
\end{definition}

\bigskip \begin{definition}[${\cal K}$-Z-property]\label{def:k_z_p}
	Let ${\cal K} \subseteq \R^n $ be a proper cone, the {\it ${\cal K}$-Z-property} of a matrix $A\in\R^{n\times n}$ means that
	\[
		\langle Ax,y\rangle \le0, \qquad \forall (x,y)\in {\cal C}({\cal K})
	\]
	where ${\cal C}({\cal K})$ is the complementarity set defined as ${\cal C}({\cal K}):=\{(x,y)\in\R^n\times\R^n:~x\in {\cal K},\textrm{ }y\in {\cal K}^*,  \langle x, y \rangle=0\}$ (see Definition \ref{def:cs} in \textbf{Part I}).
\end{definition}

The following theorem proves that when ${\cal K} = \R^n_+$, Definition \ref{def:z_m} and Definition \ref{def:k_z_p} are equivalent. 

\bigskip \begin{theorem}\label{thm:k_z_z_m}
	The matrix $A\in \R^{n\times n}$ is a Z-matrix if and only if $A$ has the $\R^n_+$-Z-property.
\end{theorem}

\begin{proof}
	Suppose that A has the $\R^n_+$-Z-property, take $x = e^i$, $y = e^j$ for any $i,j \in \{1, 2, \dots, n\}$ with $i\neq j$, $e^i$ and $e^j$ are canonical vectors of $\R^n$. We have
	\[
		\langle x, y\rangle = 0, \quad\langle Ae^i, e^j\rangle \le 0, 
	\]
	which implies that $a_{ij}\le 0$ for any $i\neq j$. Hence, $A$ is a Z-matrix.
	
	Conversely, suppose $A$ is a Z-matrix, let $a = \max_i a_{ii}$, and denote $P = aI_n - A$. Then $P$ is an entrywise nonnegative matrix. Arbitrarily take  $(x,y)\in {\cal C}(\R^n_+)$, then
	\begin{align*}
		\langle Ax,y\rangle & = \langle (aI_n - P)x,y\rangle \\
		& = a\langle x, y \rangle - \langle Px, y\rangle \\
		& = - \langle Px, y\rangle \leq 0.
	\end{align*}
	Hence, A has the $\R_+^n$-Z-property. 
\end{proof}

\bigskip \begin{definition}[Euclidean sphere and its tangent hyperplane]
	Denote by $\SP^{n-1}:= \left\{ p = (p_1, \dots, p_n)\in \R^n: ~ \|p\| =1 \right\}$ the {\it n-dimensional Euclidean sphere}, the {\it tangent hyperplane} at point $x\in\SP^{n-1}$ is
	\[
		T_{x}\SP^{n-1}:=\left\{v\in \R^n :~ \langle x, v \rangle=0, ~ x\in\SP^{n-1} \right\},
	\]
\end{definition}\index{Euclidean sphere}\index{Tangent hyperplane}

\bigskip \begin{definition}[Intrinsic distance on the sphere]\label{def:idots}
	The  {\it intrinsic distance on the sphere} between two arbitrary points \(x, y \in \SP^{n-1}\)  is  defined by
	\begin{equation} \label{eq:Intdist}
		d(x, y):=\arccos \langle x , y\rangle.
	\end{equation}
\end{definition}\index{Intrinsic distance on the sphere}

By definition \ref{def:idots}, it can  be shown that the space \((\SP^{n-1}, d)\) is a complete metric space, so that $d(x, y)\ge 0$ for any $x, y \in \SP^{n-1}$, and
$d(x, y)=0$ if and only if $x=y$. It is also easy to check that $d(x, y)\leq \pi$ for any $x, y\in \SP^{n-1}$,
and $d(x, y)=\pi$ if and only if $x=-y$.

A mapping $\gamma :[x,y] \rightarrow \SP^{n-1}$ is called a { \it geodesic}, if it is the intersection curve of a plane through the origin of \(\R^{n}\) with the sphere \( \SP^{n-1}\). For any $x, y\in \SP^{n-1}$ such that $y\neq x$ and $y\neq -x$, there exists a unique segment of minimal geodesic from $x$ to $y$. The definition of the {\it minimal geodesic} is as follows: \index{Geodesic} \index{Minimal geodesic}

\bigskip \begin{definition}[Minimal geodesic]
	For any $x, y\in \SP^{n-1}$ such that $y\neq  x$ and $y\neq -x$, the unique segment of {\it minimal geodesic from to $x$ to $y$} is
\begin{equation} \label{eq:gctp}
\gamma_{xy}(t)= \left( \cos (td(x, y)) - \frac{\langle x, y\rangle \sin (td(x, y))}{\sqrt{1-\langle x, y\rangle^2}}\right) x
+ \frac{\sin (td(x, y))}{\sqrt{1-\langle x, y\rangle^2}}\;y, \qquad t\in [0, \;1].
\end{equation}
In particular, let \(x\in \SP^{n-1}\) and  \(v\in T_{x}\SP^{n-1}\) such that \(\|v\|=1\). The minimal geodesic
connecting  \(x\) to \(-x\),  starting at \(x\) with velocity $v$ at $x$ is given by
\begin{equation} \label{eq:gctpp}
 \gamma_{x\{-x\}}(t):=\cos(\pi t) \,x+ \sin(\pi t)\, v, \qquad t\in  [0, \;1].
\end{equation}
\end{definition}

\begin{figure}[h!]
	\centering
	\includegraphics[width=.6\textwidth, height=.5\textwidth]{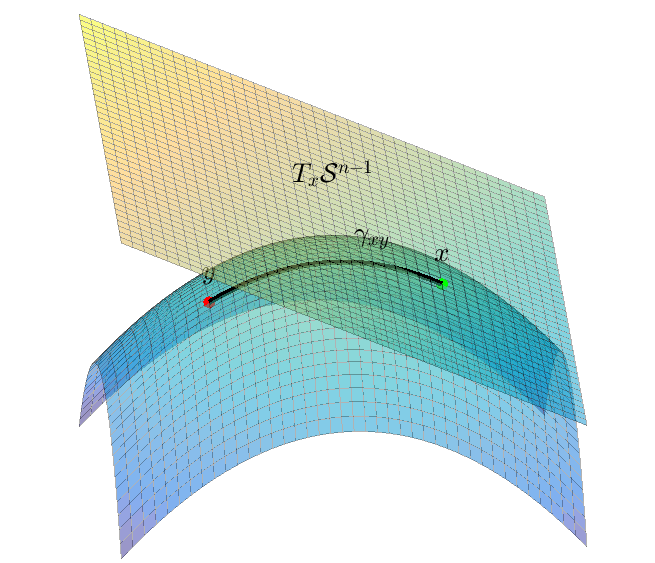} \\
	 \caption{ The geodesic and the tangent hyperplane} 
	\label{fig:geod}
\end{figure}

\bigskip \begin{definition}[Gradient on the sphere]
Let \({\cal S} \subseteq \SP^{n-1}\) be a spherically open set (i.e., a set open with respect to the induced topology in $\SP^{n-1}$). The {\it
gradient on the sphere} of a differentiable function \(f: {\cal S} \to  \R\)
at a point \(x\in {\cal S}\) is the vector defined by
\begin{equation} \label{eq:grad}
\grad f(x):= \left[I_n-xx^T \right] Df(x)= Df(x)- \langle Df(x) , x \rangle \,x,
\end{equation}
where \(Df(x) \in \R^{n}\) is the usual gradient of  \( f\) at \(x\in {\cal S}\).
\end{definition}
\index{Gradient on the sphere}

Let ${\cal D}\subseteq\R^n$ be an open set, \(I\subseteq \R\) an open interval,  \({\cal S} \subseteq
\SP^{n-1}\)  a spherically open set
and   \(\gamma:I\to {\cal S}\) a  geodesic  segment. If \(f:{\cal D} \to \R\)  is a differentiable
function, then,  since  \( \gamma'(t)\in T_{\gamma(t)} \SP^{n-1}\) for any \(t\in I\), we have $\langle \gamma'(t), \gamma(t)\rangle = 0$. The equality \eqref{eq:grad} implies
\begin{equation} \label{eq:cr1}
\frac{d}{dt}f(\gamma(t)) =\left\langle \grad f(\gamma(t)), \gamma'(t) \right\rangle=
\left\langle D f(\gamma(t)), \gamma'(t) \right\rangle, \qquad \forall ~ t\in I.
\end{equation}

\bigskip \begin{definition}[Spherically convex set] \label{def:cf}
The set \({\cal S} \subseteq \SP^{n-1}\) is said to be \emph{spherically convex} if for any $x$, $y\in {\cal S}$,
 the minimal geodesic segments from $x$ to $y$ are contained in ${\cal S}$.
\end{definition}
\index{Cone!Spherically convex set}

\bigskip \begin{example}
The set  \( S_{+}=\{(x_1,\dots, x_{n})\in \SP^{n-1}\, : \, x_{1}\geq 0,\dots, x_{n}\geq 0 \} \)	is  a closed spherically convex set.
\end{example}

We assume  for convenience that from now on all spherically convex sets are {\it nonempty proper
subsets of the sphere}.  For any set \( {\cal S} \subseteq \SP^{n-1}\), we define \({\cal K}_{\cal S}\subseteq\R^{n}\) the {\it cone spanned by} ${\cal S}$,
namely,
\begin{equation} \label{eq:pccone}
{\cal K}_{\cal S}:=\left\{ tx \, :\, x\in {\cal S}, \; t\in [0, +\infty) \right\}.
\end{equation}

\begin{figure}[h!]
	\centering
	\includegraphics[width=.48\textwidth, height=.4\textwidth]{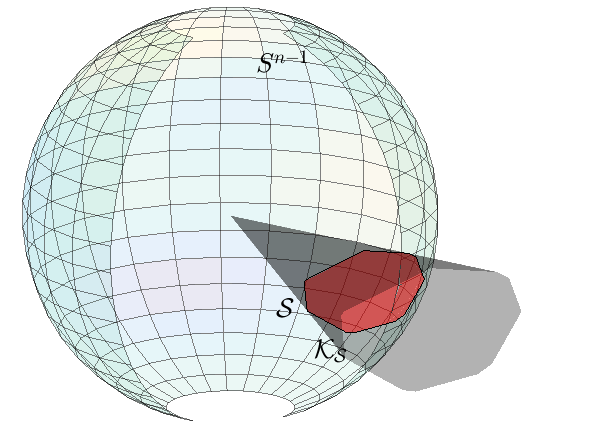} \\
	 \caption{ Closed set \({\cal S}\) and the cone \({\cal K}_{\cal S} \) spanned by \({\cal S}\).~~~~~~~~~} 
	\label{fig:conespan}
\end{figure}
Clearly, ${\cal K}_{\cal S}$ is
the smallest closed cone which contains ${\cal S}$.  The following proposition exhibits a relationship of spherically convex sets with the cones spanned by
them.
\bigskip \begin{proposition} \label{pr:ccs}\cite[Proposition~2]{FerreiraIusemNemeth2013}
The set ${\cal S}$ is spherically convex if and only if  the cone ${\cal K}_{\cal S}$ is convex and pointed.
\end{proposition}

\bigskip\begin{definition}[Spherically convex function]
	Let ${\cal S}\subseteq\SP^{n-1}$ be a spherically convex set. A function $f:{\cal S} \to \R$  is said to be \emph{(strictly) spherically convex} if for any minimal geodesic  segment $\gamma:[0, 1]\to {\cal S}$, the composition $f\circ \gamma :[0, 1]\to \R$ is (strictly) convex in the usual sense.
\end{definition}
\index{Spherically convex!Spherically convex function}

For a spherically convex set ${\cal S}\subseteq\SP^{n-1}$, the {\it sub-level sets} of a function \(f:~{\cal S} \to \R\) are denoted by
\begin{equation} \label{eq:sls}
[f\leq c]:=\{x\in {\cal S} :\; f(x)\leq c\}, \qquad c\in \R.
\end{equation}

\bigskip \begin{definition}[spherically quasi-convex function]\label{def:qcf-b}
Let \( {\cal S}\subseteq\SP^{n-1}\) be a spherically convex set.
A function \(f:{\cal S} \to \R\)  is said to be spherically quasi-convex (respectively, strictly spherically quasi-convex)
if for any minimal geodesic  segment \(\gamma: [0, 1]\to {\cal S}\), the composition
\( f\circ \gamma :[0, 1]\to \R\) is quasi-convex (respectively, strictly quasi-convex) in the usual sense, i.e.,
\(f(\gamma(t))\leq \max \{ f(\gamma(0)), f(\gamma(1))\}\) for any \(t\in [0, 1]\), (respectively,   \(f(\gamma(t))< \max \{ f(\gamma(0)), f(\gamma(1))\}\)  for any \(t\in [0, 1]\)).
\end{definition}
\index{Spherically convex!Spherically quasi-convex function}

{
From the above definition, it follows  that (strictly) spherically convex  functions are (strictly) spherically quasi-convex, but the converse is not true. It is worth to remark that the quasi-convexity concept generalises the convexity one, which was extensively studied in \cite{FerreiraIusemNemeth2014}. 

\begin{proposition}\label{pr:c2qc}
	Let ${\cal S}\subseteq\SP^{n-1}$ be a spherically convex set. If function $f:{\cal S} \to \R$  is (strictly) spherically convex, then it is (strictly) spherical quasi-convex.
\end{proposition}
\begin{proof}
	We just prove spherically convex $\Rightarrow$ spherical quasi-convex. The strict case is omitted. Suppose function $f:{\cal S} \to \R$ is spherically convex, for any $t_1$, $t_2 \in [0,~1]$, we have 
	\[
		f(\gamma(\lambda t_1+ (1-\lambda)t_2)) \le \lambda f(\gamma(t_1)) + (1- \lambda)f(\gamma(t_2))
	\]
	for any $ \lambda \in [0,~1]$. Let $t_1 =0$, $t_2 =1$, we have 
	\[
		f(\gamma( 1-\lambda)) \le \lambda f(\gamma(0)) + (1- \lambda)f(\gamma(1))\le \max \{ f(\gamma(0)), f(\gamma(1))\}
	\]
	for any $ \lambda \in [0,~1]$. 
\end{proof}
}

For subsequent use only, we denote the {\it spherically open ball} with radius $\delta >0$ and center in $x\in \SP^{n-1}$ by
\[
	\textbf{B}(x,\delta):=\{y\in \SP^{n-1} : d(x, y)<\delta \},
\]
 and the {\it spherically closed ball} with radius $\delta >0$ and center in $x\in \SP^{n-1}$ by
\[
	\bar{\textbf{B}}(x,\delta):=\{y\in \SP^{n-1} : d(x, y)\leq \delta \}.
\]

\bigskip \begin{proposition} \label{pr:charb1}
Let   \( {\cal S}\subseteq\SP^{n-1}\) be a spherically convex set.  A function \(f:{\cal S} \to \R\) is spherically quasi-convex if and only if
the sub-level sets \([f\leq c]\) (defined in \eqref{eq:sls}) are spherically convex  for any \(c\in \R\).
\end{proposition}
\begin{proof}
Suppose that $f$ is spherically quasi-convex and $c\in \R$. Arbitrarily take $x, y\in  [f\leq c]$, and let $\gamma_{xy}: [0, 1]\to  \SP^{n-1}$ be the minimal geodesic from $x$ to $y$. By \eqref{eq:gctp}, we have \(\gamma_{xy}(0)=x\) and \(\gamma_{xy}(1)=y\). Since $f$ is a spherically quasi-convex function and $x, y\in  [f\leq c]$  we have
\[
	f(\gamma_{xy}(t))\leq \max \{ f(\gamma_{xy}(0)), f(\gamma_{xy}(1))\}\leq \max \{f(x), f(y)\}\leq c,
\]
for any $t\in [0, 1]$, which implies that $\gamma(t)\in [f\leq c]$ for any \(t\in [0, 1]\). Hence we conclude that $[f\leq c]$ is a spherically convex set,  for any \(c\in \R\).

Conversely, suppose that for any \(c\in \R\), the set \([f\leq c]\) is spherically convex. For any $c$ with $f(x)\leq c$ and $f(y) \le c$, we have $x, y\in  [f\leq c]$. Without loss of generality, suppose $c = \max \{ f(x), f(y)\}$. Let $\gamma_{xy}: [0, 1]\to [f\leq c]$ be the minimal geodesic from $x$ to $y$.  By Definition \ref{def:cf}, we have $\gamma_{xy}(t)\in [f\leq c]$ for any \(t\in [0,1]\), which implies
\[
	f(\gamma_{xy}(t))\leq c =\max \{ f(x), f(y)\}= \max \{ f(\gamma(0)), f(\gamma(1))\},
\]
for any $t\in [0, 1]$. Therefore,  \(f\) is a spherically quasi-convex function.
\hfill$\square$\par \end{proof}

\bigskip \begin{proposition} \label{pr:MinPro}
Let   \( {\cal S}\subseteq\SP^{n-1}\) be a spherically convex set and \(f:{\cal S} \to \R\) be a spherically quasi-convex function. If $x^{*}\in {\cal S} $
is a strict local minimiser of $f$, then  $x^{*}$ is also a strict global minimiser of $f$ in ${\cal S}$.
\end{proposition}
\begin{proof}
Suppose that there exists $x^{*}$ is a strict local minimiser of $f$, then there exists a scalar $\delta >0$ such that
\begin{equation} \label{eq:lmg}
	f(x) > f(x^{*}), \qquad \forall ~ x\in \textbf{B}(x^*,\delta)\setminus \{x^*\}=\{y\in {\cal S} ~:~ 0 < d(y, x^*)<\delta\}.
\end{equation}
Assume by contradiction that $x^{*}$ is not a strict global minimiser of $f$ in ${\cal S}$. Thus, there exists ${\bar x}\in {\cal S} $ with
${\bar x}\neq x^*$ such that  $f({\bar x})\leq f(x^{*})$. Since $C$ is  spherically convex, we can   take  a  minimal
geodesic  segment \(\gamma_{x^*{\bar x}}:[0, 1]\to {\cal S}\) joining   $x^*$ and ${\bar x}$, then we have
\[
	\gamma_{x^*{\bar x}}(0)=x^*,~\quad \gamma_{x^*{\bar x}}(1)={\bar x}.
\]
Considering that \(f\) is spherically quasi-convex, by Definition \ref{def:qcf-b} we have
\begin{equation}\label{eq:rmg}
	f(\gamma_{x^*{\bar x}}(t))\leq  \max \{ f(x^*), f({\bar x})\}=f(x^{*})
\end{equation}
for any \(t\in [0, 1]\). On the other hand, for $t$ sufficiently small  we have $\gamma_{x^*{\bar x}}(t)\in \textbf{B}(x^*,\delta)$, which means
\[
	f(\gamma_{x^*{\bar x}}(t)) > f(x^*)
\]
Therefore, the inequality \eqref{eq:rmg} contradicts \eqref{eq:lmg}.
\hfill$\square$\par \end{proof}

\bigskip \begin{proposition} \label{pr:MinProS}
Let   \( {\cal S}\subseteq\SP^{n-1}\) be a spherically convex set and   \(f:{\cal S} \to \R\) be a strictly spherically  quasi-convex function. Then   $f$ has
at most one local  minimiser point which   is also a  global  minimiser point of $f$.
\end{proposition}
\begin{proof}
Without loss of generality assume by contradiction that the strictly  spherically quasi-convex \(f(\gamma(t)) <  \max \{ f(x^*), f({\bar x})\}\) for any \(t\in [0, 1]\).  Since we can take  $t$ sufficiently close to $0$ or $1$, the last inequality  function $f$ has two local minimiser  $x^*, {\bar x}\in {\cal S} $ with ${\bar x}\neq x^*$. Thus,
\[
	f(x) > f(x^{*}), \qquad \forall ~ x\in \textbf{B}(x^*,\delta)\setminus \{x^*\}=\{y\in {\cal S} ~:~ 0 < d(y, x^*)<\delta\},
\]
and
\[
	f(x) > f({\bar x}), \qquad \forall ~ x\in \textbf{B}({\bar x},\delta)\setminus \{{\bar x}\}=\{y\in {\cal S} ~:~ 0 < d(y, {\bar x})<\delta\}.
\]
we can take a minimal geodesic segment \(\gamma_{x^*{\bar x}}:[0, 1]\to {\cal S}\) joining $x^*$ and ${\bar x}$, then we have
\[
	\gamma(0)=x^*, ~\quad \gamma(1)={\bar x}.
\]
Due to \(f\) being  strictly  spherically quasi-convex, we have
\[
	f(\gamma(t)) <  \max \{ f(x^*), f({\bar x})\}
\]
for any \(t\in [0, 1]\).  If we can take $t$ sufficiently close to $0$ or $1$, it gives $f(\gamma(t)) > f(x^*)$ or $f(\gamma(t)) > f({\bar x})$, respectively. The last two inequalities contradicts the assumption that  $x^*, {\bar x}$  are two distinct  local minimisers. Thus, $f$ has at most one local minimiser point.  Since  $f$ is strictly quasi-convex, the local minimiser is strict. Therefore, Proposition~\ref{pr:MinPro} implies that  the  local minimiser point is global.
\hfill$\square$\par \end{proof}
\bigskip \begin{proposition} \label{pr:CharDiff}
Let   \( {\cal S}\subseteq\SP^{n-1}\) be an open spherically convex set and   \(f:{\cal S} \to \R\) be a differentiable function. Then $f$ is spherically quasi-convex if and only if
\begin{equation} \label{eq:cqqudf2}
f(x)\leq f(y) \Longrightarrow \langle Df(y) , x \rangle - \langle x , y \rangle\langle Df(y) , y \rangle \leq 0 , \qquad \forall ~ x,  y\in {\cal S} .
\end{equation}
\end{proposition}
\begin{proof}
Let  \(\gamma:I\to {\cal S}\)  be a geodesic  segment and consider  the composition \( f\circ \gamma :I\to \R\).	The usual characterisation of scalar quasi-convex functions implies that \( f\circ \gamma\) is quasi-convex  if and  only if
\begin{equation} \label{eq:cqqud}
 f(\gamma(t_1))\leq f(\gamma(t_2)) \Longrightarrow  \frac{d}{dt}f(\gamma(t_2))(t_1-t_2)\leq 0, \qquad  \forall ~ t_1, t_2 \in I.
\end{equation}
On the other hand,  for each  $ x,  y\in {\cal S}$ with $y\neq x$  we have from  \eqref{eq:gctp} that  $\gamma_{xy}$ is the  minimal  geodesic segment from \(x=\gamma_{xy}(0)\) to \(y=\gamma_{xy}(1)\) and
$$
\gamma_{xy}'(1)=\displaystyle \frac{\arccos \langle x, y\rangle}{\sqrt{1-\langle x, y\rangle ^2}}
\left(yy^T-I_n\right)x\in T_{y}\SP^{n-1} , \qquad y\neq -x.
$$
Note that letting $x=\gamma(t_1)$ and  $y=\gamma(t_2)$ we have that $\gamma_{xy}(t)=\gamma(t_1+t(t_2-t_1))$.
Therefore, by using \eqref{eq:cr1}  we can rewrite the right hand side of \eqref{eq:cqqud} as
\begin{align*}
  \frac{d}{dt}\left(f(\gamma(t_2))\right)(t_1-t_2) &  = \langle \grad f(\gamma(t_2)), \gamma'(t_2)\rangle (t_1-t_2)\\
   & = \Big\langle \grad f(\gamma(t_2)), \frac{\gamma(t_2) - \gamma(t_1)}{t_2 - t_1}\Big\rangle (t_1-t_2)\\
   & = \langle \grad f(\gamma(t_2)), \gamma(t_1)\rangle - 0\\
   & = \Big\langle Df(t_2) - \langle Df(t_2),\gamma(t_2)\rangle \gamma(t_2), ~\gamma(t_1)\Big\rangle\\
   & = \langle Df(t_2),\gamma(t_1)\rangle - \langle Df(t_2),\gamma(t_2)\rangle \langle \gamma(t_2), \gamma(t_1)\rangle \le 0
\end{align*}
which is equivalent to \eqref{eq:cqqudf2}.
\hfill$\square$\par \end{proof}

\chapter{Linear Complementarity Problems on Extended Second Order Cones}\label{cht:lesoc}

In this chapter we elaborate the formulation and the solution to \gls{esoclcp}. We present the major characterisation of ESOCLCP in Theorem \ref{thm:main}. Based on Theorem \ref{thm:main}, an ESOCLCP can be converted into a mixed complementarity problem on the nonnegative orthant. We state necessary and sufficient conditions for a point to be a solution to the converted problem. We also present solution strategies for this problem, as well as some numerical examples. The results in this chapter are published in the paper \cite{nemeth2018linear}, co-worked with my supervisor. In order to improve the readability of our results, in this chapter, we give more explicit proofs and more detailed explanations about these results. 

\section{Problem formulation}
 Let $T=\bs A & B\\C & D \es$, with $A\in\R^{k\times k}$,
	$B\in\R^{k\times\ell}$, $C\in\R^{\ell\times k}$ and $D\in\R^{\ell\times\ell}$. Let $r = \bs p \\ q\es$ with $p\in\R^k$ , $q\in\R^\ell$. The linear complementarity problem defined by{ the} extended second order cone $L$ and a linear function $F(x,u) = T\bs x\\u\es +r$ is: 
\begin{equation}\label{eq:lcp}
		\LCP(F,L)\left\{
		\begin{array}{l}
			Find \; \bs x\\u\es\in L,\; such \; that\\
			F(x,u) \in M \; and \; \langle \bs x\\u\es, F(x,u) \rangle = 0.
		\end{array}
		\right.
\end{equation}

Based on the idea of Proposition \ref{prop:cp_mixcp} and Proposition \ref{prop:icp_mixicp}, using Proposition \ref{prop:cs-esoc}, we developed the following theorem, which shows the equivalence {of} various complementarity problems associated with ESOC. For convenience, let $n = k +\ell$.

\bigskip \begin{theorem}\label{thm:main}
	Denote $z=\bs x\\u\es$, {$\hat z=\bs \hat x\\u\es:=\bs x-\|u\|e\\u\es$}, $\widetilde{z}=\bs \tilde{x} \\u\\t\es := \bs x-t\\u\\t\es$ and $r=\bs p\\ q\es$ with $x,p\in\R^k$ , $u,q\in\R^\ell$, and $t\in \R$. Let $T=\bs A & B\\C & D \es$ with $A\in\R^{k\times k}$,
	$B\in\R^{k\times\ell}$, $C\in\R^{\ell\times k}$ and $D\in\R^{\ell\times\ell}$. The square matrices $T$, $A$ and $D$ are nonsingular. Let $L$ be the extended second order come.
	\begin{enumerate}
		\item[(i)] Suppose $u= 0$. We have
			\begin{align*}
				z\in\SLCP(T,r,L)& \\  \iff & x\in\SLCP(A,p,\R^k_+)\mbox{ and }e\tp(Ax+p)\ge\|Cx+q\|.
			\end{align*}
		\item[(ii)] Suppose $Cx+Du+q= 0$. Then, \[z\in\SLCP(T,r,L)\iff z\in\SMixCP(F_1,F_2,\R^k_+)
			\mbox{ and }x\ge\|u\|,\] where $F_1(x,u)=Ax+Bu+p$ and $F_2(x,u)=0$.
		\item[(iii)] Suppose $u\ne 0$ and $Cx+Du+q\ne 0$.
			We have \[z\in\SLCP(T,r,L)\iff z\in\SMixICP(G_1,F_1,F_2,\R^k_+),\] where
			\[F_2(x,u)=\lf(\|u\|C+ue\tp A\rg)x+ue\tp(Bu+p)+\|u\|(Du+q),\] $G_1(x,u)=x-\|u\|e$ and $F_1(x,u)=Ax+Bu+p$.
		\item[(iv)] Suppose $u\ne 0 $ and $Cx+Du+q\ne 0$.
			We have \[z\in\SLCP(T,r,L)\iff \hat z\in\SMixCP(F_1,F_2,\R^k_+),\] where
			\[F_2(\hat{x},u)=\lf(\|u\|C+ue\tp A\rg)(\hat{x}+\|u\|e)+ue\tp(Bu+p)+\|u\|(Du+q)\] and $F_1(\hat{x},u)=A(\hat{x}+\|u\|e)+Bu+p$.
		\item[(v)] Suppose $u\ne 0$, $Cx+Du+q\ne 0$ and $\|u\|C+u\tp e A$ is a nonsingular matrix. We have
			\[z\in\SLCP(T,r,L)\iff \hat z\in\SICP(F_1,F_2,\R^k_+),\] where
			\[F_1(u)=A\lf(\lf(\|u\|C+ue\tp A\rg)^{-1}\lf(ue\tp(Bu+p)+\|u\|(Du+q)\rg)\rg)+Bu+p\] and
			\[F_2(u)=\lf(\|u\|C+ue\tp A\rg)^{-1}\lf(ue\tp(Bu+p)+\|u\|(Du+q)\rg).\]
		\item[(vi)] Suppose $u\ne 0$, $Cx+Du+q\ne 0$. We have
			\[z\in\SLCP(T,r,L)\iff \exists t>0,\] such that \[\tilde{z}\in\MixCP(\widetilde{F}_1,\widetilde{F}_2,\R^k_+),\] where
			\begin{equation}\label{eq:mathcal_F1}
				\widetilde{F}_1(\tilde{x},u,t)=A(\tilde{x}+te)+Bu+p
			\end{equation}
			and
			\begin{equation}\label{eq:mathcal_F2}
				\widetilde{F}_2(\tilde{x},u,t) =
				\begin{pmatrix}
					& \lf(tC+ue\tp A\rg)(\tilde{x}+te)+ue\tp(Bu+p)+t(Du+q) \\
					& t^2 - \|u\|^2
				\end{pmatrix}.
			\end{equation}

	\end{enumerate}
\end{theorem}

\begin{proof}
	\begin{enumerate}
		\item[(i)] We have that $z\in\SLCP(T,r,L)$ is equivalent to
			$(x,0,Ax+p,Cx+q)\in\C(L)$ or, by item (i) and (ii) of Proposition \ref{prop:cs-esoc}, equivalent to
			$(x,Ax+p)\in\C(\R^k_+)$ and $e\tp(Ax+p)\ge\|Cx+q\|$.
		\item[(ii)] We have that $z\in\SLCP(T,r,L)$ is equivalent to
			$(x,u,Ax+Bu+p,0)\in\C(L)$ or, by item (i) and (iii) of Proposition \ref{prop:cs-esoc}, equivalent to
			$(x,Ax+Bu+p)\in\C(\R^k_+)$ and $x\ge\|u\|$, or equivalent to
			\[z\in\SMixCP(F_1,F_2,\R^k_+)\mbox{ and }x\ge\|u\|,\] where $F_1(x,u)=Ax+Bu+p$ and
			$F_2(x,u)= 0$.
		\item[(iii)] Suppose that $z\in\SLCP(T,r,L)$. Then, $(x,u,y,v)\in\C(L)$, where $y=Ax+Bu+p$ and $v=Cx+Du+q$. Then, by item (iv) of Proposition \ref{prop:cs-esoc} we obtain that $\exists\lambda>0$ such that
			\begin{equation}\label{parall-eq}
				Cx+Du+q=v=-\lambda u,
			\end{equation}
			\begin{equation}\label{prod-eq}
				e\tp (Ax+Bu+p)=e\tp y=\|v\|=\|Cx+Du+q\|=\lambda\|u\|,
			\end{equation}
			\begin{equation}\label{cpset-eq}
				\lf(G_1(x,u),F_1(x,u)\rg)=(x-\|u\|e,Ax+Bu+p)=(x-\|u\|e,y)\in\C(\R^k_+).
			\end{equation}
			From equation \eqref{parall-eq} we obtain $\|u\|(Cx+Du+q)=-\lambda\|u\| u$, which by equation \eqref{prod-eq} implies
			$\|u\|(Cx+Du+q)=-ue\tp (Ax+Bu+p)$, which after some algebra gives
			\begin{equation}\label{zero-eq}
				F_2(x,u)= 0.
			\end{equation}
			From equations \eqref{cpset-eq} and \eqref{zero-eq} we conclude that $z\in\SMixICP(F_1,F_2,G_1)$.
			\medskip
			
			Conversely suppose that $z\in\SMixICP(F_1,F_2,G_1)$. Then,
			\begin{equation}\label{zero-eq2}
				\|u\|v+ue\tp y=\|u\|(Cx+Du+q)+ue\tp (Ax+Bu+p)=F_2(x,u)= 0
			\end{equation}
			and
			\begin{equation}\label{cpset-eq2}
				(x-\|u\|e,y)=(x-\|u\|e,Ax+Bu+p)=(G_1(x,u),F_1(x,u))\in\C(\R^k_+),
			\end{equation}
			where $v=Cx+Du+q$ and $y=Ax+Bu+p$. Equations \eqref{cpset-eq2} and \eqref{zero-eq2} imply
			\begin{equation}\label{parall-eq2}
				v=-\lambda u,
			\end{equation}
			where
			\begin{equation}\label{lambda-eq}
				\lambda=(e\tp y)/\|u\|>0.
			\end{equation}
			Equations \eqref{parall-eq2} and \eqref{lambda-eq} imply
			\begin{equation}\label{norm-v-eq}
				e\tp y=\|v\|
			\end{equation}
			By item (iv) of Proposition \ref{prop:cs-esoc}, equations \eqref{parall-eq2}, \eqref{norm-v-eq} and \eqref{cpset-eq2}
			imply $(x,y,u,v)\in C(L)$ and therefore $z\in\SLCP(T,r,L)$.
		\item[(iv)] It is a simple reformulation of item (iii) by using the change of variables
			\Large
			\[\bs x\\u\es\to \bs \hat{x} \\u\es :=  \bs x-\|u\|e\\u\es.\]
			\normalsize
		\item[(v)] It is a simple reformulation of item (iv) by using that $\|u\|C+u\tp e A$ is a nonsingular matrix.
		\item[(vi)]  Suppose that $z\in\SLCP(T,r,L)$. Then, $(x,u,y,v)\in\C(L)$, where $y=Ax+Bu+p$ and $v=Cx+Du+q$. Let $t = \|u\|$, Then, by item (iv)
			of Proposition \ref{prop:cs-esoc} we have that $\exists\lambda>0$ such that
			\begin{equation}\label{parall-eq3}
				Cx+Du+q=v=-\lambda u,
			\end{equation}
			\begin{equation}\label{prod-eq3}
				e\tp (Ax+Bu+p)=e\tp y=\|v\|=\|Cx+Du+q\|=\lambda t,
			\end{equation}
			\begin{equation}\label{cpset-eq3}
				(\tilde{x},\widetilde{F}_1(\tilde{x},u,t))=(x-te,Ax+Bu+p)=(x-te,y)\in\C(\R^k_+)
			\end{equation}
			where $\widetilde{z} = \bs \tilde{x}\\u\\t\es := \bs x-t\\u\\t\es\in\R^k\times\R^{\ell}\times\R$. From equation \eqref{parall-eq3} we obtain $t(Cx+Du+q)=-t \lambda u$, which by equation \eqref{prod-eq3} implies
			$t(Cx+Du+q)=-ue\tp (Ax+Bu+p)$, which after some algebra gives
			\begin{equation}\label{zero-eq3}
				\widetilde{F}_2(\tilde{x},u,t)= 0.
			\end{equation}
			Equations \eqref{cpset-eq3} and \eqref{zero-eq3} yield $z\in\SMixCP(\widetilde{F}_1,\widetilde{F}_2,\R^k_+)$.
			\medskip

	\end{enumerate}
\hfill$\square$\par \end{proof}

{
\textbf{Comment:} Many well-developed methods to a complementarity problem \eqref{eq:cp} are based on a smooth function $F$ \cite{FacchineiPang2003, chen1996class, chen1997smooth}. Hence, by modifying the semi-smooth function $F_1(\hat{x},u)$ in item (iv), we introduce item (vi) with a smooth function $\widetilde{F}_1(\tilde{x},u,t)$.

\textbf{Comment:} Converting an ESOCLCP to a MixCP will very likely increase its complexity, because it converts a linear problem to a nonlinear one. However, due to lacking methods for solving an ESOCLCP, we have to use other available methods to solve it. Given the fact that there are many methods, especially the complementarity function (C-function) method, proposed for solving the complementarity problem on nonnegative orthant \cite{kanzow1996nonlinear, fischer1995newton, fischer1992special, mangasarian1976equivalence}, we are therefore motivated to implement such conversion.}
 As it is converted to a MixCP, we will be able to solve the ESOCLCP by means of a C-function.
 
The scalar form of \emph{\gls{fb} C-function} \cite{ fischer1995newton} is defined as:

\[
	\psi_{FB}(a,b) = \sqrt{a^2+b^2} - (a+b) \quad \forall (a,b) \in \mathbb{R}^2.
\]

The equivalent FB-based formulation of MixCP is:

\begin{equation}\label{FBform}
	\mathbb{F}^{\MixCP}_{FB}(x,u,t) :=
	\begin{pmatrix}
		 \psi_{FB}\lf(x_1,(\widetilde{F}_1)_1(x,u,t)\rg) \\
		 \vdots \\
		 \psi_{FB}\lf(x_k,(\widetilde{F}_1)_k(x,u,t)\rg) \\
		 \widetilde{F}_2(x,u,t)
	\end{pmatrix}
\end{equation}

The FB-based formulation of $\MixCP$ is semi-smooth. Based on the property of FB C-function, if there is a point $ \bs x^*\\u^*\\t^*\es$ such that
\begin{equation}\label{eq:FBCFque}
	\mathbb{F}^{\MixCP}_{FB}(x^*,u^*,t^*) = 0,
\end{equation}
then $\bs x^*\\u^*\\t^*\es$ is a solution to  $\MixCP$. The equation \eqref{eq:FBCFque} is semi-smooth, but it still can be solved by using semi-smooth Newton's method. Denote by $\partial\mathbb{F}^{\MixCP}_{FB}(x,u,t)$ the generalised Jacobian set of $\mathbb{F}^{\MixCP}_{FB}(x,u,t)$ . 
Since $\mathbb{F}^{\MixCP}_{FB}(x,u,t)$ is semi-smooth, we have that $\partial\mathbb{F}^{\MixCP}_{FB}(x,u,t)$ satisfies

\[
	\partial\mathbb{F}^{\MixCP}_{FB}(x,u,t)\subseteq 
	\begin{pmatrix} 
		\mathcal{D}_a(x,u,t) + \mathcal{D}_b(x,u,t)J_x\widetilde{F}_1(x,u,t) & \mathcal{D}_b(x,u,t)J_{\bs u\\t\es}\widetilde{F}_1(x,u,t)\\ 
		J_x\widetilde{F}_2(x,u,t) & J_{\bs u\\t\es}\widetilde{F}_2(x,u,t)
	\end{pmatrix}.
\]
where $ \mathcal{D}_a$ and $ \mathcal{D}_b$ are $k\times k$ diagonal matrices respectively denoted by $diag (a_1(x,u,t)$, $\dots$, $a_k(x,u,t) )$ and $diag(b_1(x,u,t)$, $\dots$, $b_k(x,u,t))$ , with $\bar{\textbf{B}}(x,1)$ denoting a closed unit ball centered at the point $x$:
\begin{equation}\label{eq:clball}
	\lf(a_i(x,u,t),b_i(x,u,t)\rg)= \left\{
			\begin{array}{ll}
				=\frac{(x_i,(\widetilde{F}_1)_i(x,u,t))}{\sqrt{x_i^2 + (\widetilde{F}_1)_i^2(x,u,t)}} - (1,1) & if\; \lf(x_i, (\widetilde{F}_1)_i(x,u,t)\rg) \neq (0,0)\\
				\in \bar{\textbf{B}}\lf((0,0),1\rg) - (1,1) & if\; \lf( x_i, (\widetilde{F}_1)_i(x,u,t)\rg) = (0,0)
			\end{array}
			\right.
\end{equation}
Specifically, for $i \in (1,\dots,k)$ such that $(x_i,\widetilde{F}_1^i(x,u,t)) \neq (0,0)$ the i-th FB-based formulation $\lf( \mathbb{F}^{\MixCP}_{FB}\rg)_i $ is differentiable at $\bs x \\u \\t \es$. Take an element
\[
	 \mathcal{A} \in \partial\mathbb{F}^{\MixCP}_{FB}(x,u,t).
\]

Denoting by $e^i = (0,0,\dots,1,\dots,0)\tp$ the $i$-th coordinate vector, and the i-th row of the Jacobian $\mathcal{A}_i $,  which is the derivative of $\psi_{FB}\lf(x_i,(\widetilde{F}_1)_i(x,u,t)\rg)$ with respect to $x$, is shown as follows:
\begin{align*}
	\lf(\mathcal{A}_x\rg)_i (x,u,t) = & \frac{\partial\psi_{FB}\lf(x_i,(\widetilde{F}_1)_i(x,u,t)\rg)}{\partial x} = \; a_i(x,u,t)e^i + b_i(x,u,t)J_x(\widetilde{F}_1)_i(x,u,t) \\
	=&\lf(\frac{x_i}{\sqrt{x_i^2+(\widetilde{F}_1)_i^2(x,u,t)}}-1\rg)e^i  +\lf(\frac{(\widetilde{F}_1)_i(x,u,t)} {\sqrt{x_i^2+(\widetilde{F}_1)_i^2(x,u,t)}}-1\rg) J_x(\widetilde{F}_1)_i(x,u,t)			
\end{align*}

Similarly, for $i \in (1,\dots,k)$ such that $(x_i,\widetilde{F}_1^i(x,u,t))\neq (0,0)$, the i-th row of Jacobian $\mathcal{A}_i $ with respect to $\bs u \\t\es$ is:
\begin{align*}
	\lf(\mathcal{A}_{\bs u \\t\es}\rg)_i (x,u,t) = & \frac{\partial\psi_{FB}\lf(x_i,(\widetilde{F}_1)_i(x,u,t)\rg)}{\partial \bs u \\t\es} = \;b_i(x,u,t)J_{\bs u \\t\es}(\widetilde{F}_1)_i(x,u,t) \\
	= & \lf(\frac{(\widetilde{F}_1)_i(x,u,t)}{\sqrt{x_i^2+(\widetilde{F}_1)_i^2(x,u,t)}} - 1\rg) J_{\bs u \\t\es}(\widetilde{F}_1)_i(x,u,t).
\end{align*}

By \eqref{eq:clball}, for $i \in (1,\dots,k)$, if the pair  $(x_i,\widetilde{F}_1^i(x,u,t)) = (0,0)$, since $\mathbb{F}^{\MixCP}_{FB}$ is semi-smooth at origin, the Jacobian $\mathcal{A}_i$ at the origin will be a generalised Jacobian of a composite function provided that $\partial \|(0,0)\| = \bar{\textbf{B}}\lf((0,0),1\rg)$. We have
\[
	\lf(\mathcal{A}_x\rg)_i (x,u,t) = \lf\{\lf(\tilde{a}_ie^i +\tilde{b}_iJ(\widetilde{F}_1)_i(x,u,t)\rg) : (\tilde{a},\tilde{b}) \in \bar{\textbf{B}}\lf((0,0),1\rg)- (1,1)\rg\},
\]
for $i \in (1,\dots,k)$ and
\[
	\lf(\mathcal{A}_{\bs u \\t\es}\rg)_i (x,u,t) = \lf\{\lf(\tilde{a}_i\cdot 0 +\tilde{b}_iJ(\widetilde{F}_1)_i(x,u,t)\rg) : (\tilde{a},\tilde{b}) \in \bar{\textbf{B}}\lf((0,0),1\rg) - (1,1)\rg\},
\]
for $i \in (k+1,\dots,n+1)$. 

For convenience, for $i \in (1,\dots,k)$ satisfying $(x_i,\widetilde{F}_1^i(x,u,t)) = (0,0)$, we choose
\[
	\tilde{a}_i(x,u,t) = 0-1, \qquad \tilde{b}_i(x,u,t) = 0-1.
\]
It is easy to prove that $(\tilde{a},\tilde{b}) \in \bar{\textbf{B}}\lf((0,0),1\rg)- (1,1)$. Then we conclude
\begin{align*}
	\left(\mathcal{A}_x\right)_i (x,u,t) = &\;  \tilde{a}_i(x,u,t)e^i + \tilde{b}_i(x,u,t)J_x(\widetilde{F}_1)_i(x,u,t) \\
	=&-e^i - J_x(\widetilde{F}_1)_i(x,u,t),
\end{align*}

and
\begin{align*}
	\left(\mathcal{A}_{\bs u\\t\es}\right)_i (x,u,t) = & \;\tilde{b}_i(x,u,t)J_{\bs u\\t\es}(\widetilde{F}_1)_i(x,u,t) \\
	= & - J_{\bs u\\t\es}(\widetilde{F}_1)_i(x,u,t).
\end{align*}

Moreover, by the continuous differentiability of $\widetilde{F}_2(x,u,t)$, for $i \in ( k+1,\dots,n+1)$, the Jacobian $\mathcal{A}_i$ is:
\[
	 {\mathcal{A}_i} = \begin{pmatrix}J_x(\widetilde{F}_2)_i(x,u,t) & J_{\bs u \\t\es}(\widetilde{F}_2)_i(x,u,t)\end{pmatrix}.
\]


Hence, the Jacobian matrix for $\mathbb{F}^{\MixCP}_{FB}(x,u,t) $ can be written as:
\begin{equation}\label{eq:jmerit}
	\mathcal{A} =
		\begin{pmatrix}
			D_a+D_bJ_x\widetilde{F}_1(x,u,t) & D_bJ_{\bs u \\t\es}\widetilde{F}_1(x,u,t)\\
			J_x\widetilde{F}_2(x,u,t) & J_{\bs u \\t\es}\widetilde{F}_2(x,u,t)
		\end{pmatrix},
\end{equation}
where $D_a$ and $D_b$ are nonpositive definite diagonal matrices:
\[
	(D_a)_{ii}:= \left\{
		\begin{array}{ll}
			\frac{x_i}{\sqrt{x_i^2+( \widetilde{F}_1)_i^2(x,u,t)}}-1, &
			 if\; \lf(x_i, (\widetilde{F}_1)_i(x,u,t)\rg) \neq (0,0),\\
			-1, &
			if\; \lf( x_i, (\widetilde{F}_1)_i(x,u,t)\rg) = (0,0),
		\end{array}\right.
	\qquad i\in\{1, \dots , k\},
\]

\[
	(D_b)_{ii}:=\left\{
	\begin{array}{ll}
		 \frac{( \widetilde{F}_1)_i(x,u,t)} {\sqrt{x_i^2+( \widetilde{F}_1)_i^2(x,u,t)}}-1,   &
			 if\; \lf(x_i, (\widetilde{F}_1)_i(x,u,t)\rg) \neq (0,0),\\
		-1, &
			if\; \lf( x_i, (\widetilde{F}_1)_i(x,u,t)\rg) = (0,0),
	\end{array}\right.
	 \qquad  i\in\{1, \dots , k\}.
\]
Detailed methods of solving the semi-smooth equation \eqref{eq:FBCFque} will be introduced in next section. 

{
\section{Find the solution}

Many methods are proposed to solve the complementarity problem. Two of the most popular methods are: to reformulate the complementarity problem to a system of nonlinear equations; and, to reformulate it to an unconstrained minimisation problem. In this section, we will present the details for both methods. 

\subsection{Reformulate to a system of nonlinear equations}

In the previous section, we provided the equivalent FB-based formulation of MixCP \eqref{FBform}. Equation \eqref{eq:FBCFque} enables us to solve the MixCP as an unconstrained system of nonlinear equation: For convenience, we restate the equation \eqref{eq:FBCFque} here: 
\[
	\mathbb{F}^{\MixCP}_{FB}(x^*,u^*,t^*) = 0.
\]

The study about solving the complementarity problem in such nonlinear equation system is abundant. Since the FB C-function is not continuously differentiable, many semi-smooth methods are developed for solving the complementarity problem. These includes Newton-type methods \cite{andreas1995local, harker1990finite, qi1993convergence, facchinei1996inexact}, Levenberg-Marquardt methods\cite{facchinei1997nonsmooth, ma2007globally}, etc. Theses methods are proved to have at least a linear rate of convergence under certain assumptions, among which the nonsingularity is the most important assumption to guarantee the convergence of these algorithms. We use the following proposition to obtain conditions for the nonsingularity of the generalised Jacobian of $\mathbb{F}^{\MixCP}_{FB}(x^*,u^*,t^*)$.

Before stating the proposition, we define the following index sets:
\[
	\begin{array}{l}
	\alpha := \lf\{i: x_i = 0 < (\widetilde{F}_1)_i(x,u,t)  \rg\}, \\
	\beta  := \lf\{i: x_i = 0 = (\widetilde{F}_1)_i(x,u,t)  \rg\},  \\
	\gamma := \lf\{i: x_i > 0 = (\widetilde{F}_1)_i(x,u,t)  \rg\},  \\
	\delta := \lf\{1, \dots, k \rg\}\setminus\lf(\alpha\cup\beta\cup\gamma\rg).  \\
	\end{array}
\]

\begin{proposition}
	\cite[Proposition 9.4.2]{FacchineiPang2003} If $\widetilde{F}_1(x,u,t)$  and $\widetilde{F}_1(x,u,t)$ are continuously differentiable, given $\bs x \\ u \\ t \es\in \R^{k+\ell+1}$. Let $\overline{\alpha}:= \gamma\cup\beta\cup\delta$ be the complement of $\alpha$ in $\{ 1, \dots, k\}$. Assume that
	\begin{enumerate}
		\item[(i)] the submatrices
		\[
			\begin{pmatrix}
				J_{\bs u \\ t\es} \widetilde{F}_2(x,u,t) & 
				J_{x_{\tilde{\gamma}}} \widetilde{F}_2(x,u,t) \\
				J_{\bs u \\ t\es} (\widetilde{F}_1)_{\tilde{\gamma}}(x,u,t) & 
				J_{x_{\tilde{\gamma}}} (\widetilde{F}_1)_{\tilde{\gamma}}(x,u,t)
			\end{pmatrix}
		\]
		are nonsingular for all $\tilde{\gamma}$ satisfying 
		\[
			\gamma \subseteq \tilde{\gamma}\cup\gamma \cup \beta,
		\]
		\item[(ii)] the Schur complement of 
		\[
			\begin{pmatrix}
				J_{\bs u \\ t\es} \widetilde{F}_2(x,u,t) & 
				J_{x_{\gamma}} \widetilde{F}_2(x,u,t) \\
				J_{\bs u \\ t\es} (\widetilde{F}_1)_{\gamma}(x,u,t) & 
				J_{x_{\gamma}} (\widetilde{F}_1)_{\gamma}(x,u,t)
			\end{pmatrix}
		\]
		in
		\[
			\begin{pmatrix}
				J_{\bs u \\ t\es} \widetilde{F}_2(x,u,t) & 
				J_{x_{\overline{\alpha}}} \widetilde{F}_2(x,u,t) \\
				J_{\bs u \\ t\es} (\widetilde{F}_1)_{\overline{\alpha}}(x,u,t) & 
				J_{x_{\overline{\alpha}}} (\widetilde{F}_1)_{\overline{\alpha}}(x,u,t)
			\end{pmatrix}
		\]
		is a $P_0$ matrix,
	\end{enumerate}
	then the Jacobian of $\mathbb{F}^{\MixCP}_{FB}(x,u,t) $ \eqref{eq:jmerit} is nonsingular.
\end{proposition}

First , we  illustrate the semi-smooth inexact Newton's Method.

\begin{flushleft}
\textbf{Algorithm 1} (Semi-smooth Inexact Newton's method)\cite{facchinei1996inexact}:
\end{flushleft}

\textbf{Input}: the initial point $z_0 :=\bs x_0 \\ u_0 \\ t_0 \es \in \R^{k+\ell+1}$, and the tolerance $\eta_0\in\R_{+}$.

\textbf{Step 1}: Set $j = 0$.

\textbf{Step 2}: If $ \mathbb{F}^{\MixCP}_{FB}(z_j) = 0$, then STOP.

\textbf{Step 3}: Select an element $\mathcal{A}$ in the generalised Jacobian set $ \partial\mathbb{F}^{\MixCP}_{FB}(x,u,t)$, and find a direction $d_j \in \R^{k+\ell+1} $ such that
\[
	\mathbb{F}^{\MixCP}_{FB}(z_j) + \mathcal{A}\tp (z_j)d_j = r_j,
\]
where the residual vector $r_j\in \R^{k+\ell+1}$ satisfying 
\[
	\|r_j\|\le \eta_j\|\mathbb{F}^{\MixCP}_{FB}(z_j)\|.
\]

\textbf{Step 4}: Choose $\eta_{j+1}\ge 0$; set $z_{j+1} := z_j + d_j $ and $j := j + 1$; go to \textbf{Step 2}.

The above algorithm is a modification of the semi-smooth algorithm introduced by Qi and Sun \cite{qi1993nonsmooth}. It is worth noting that there are many other Newton-type algorithms for solving a complementarity problem as a system of nonlinear equations. We refer interested readers to \cite{sun1999regularization, chen2008regularization, zhang2004nonmonotone}.

The following theorem is from \cite[Theorem~3.2]{facchinei1996inexact}. It proves that the semi-smooth inexact Newton's Method at least Q-linearly converges to a solution to \eqref{eq:FBCFque}.

\bigskip \begin{theorem}\label{thm:newton}
	Let $\mathbb{F}^{\MixCP}_{FB}(x,u,t)$ be semi-smooth in $\textbf{B}\lf(z^*,\delta\rg)$, where $\delta >0$, and $z^*:= \bs x^* \\ u^* \\ t^* \es$ satisfies $\mathbb{F}^{\MixCP}_{FB}(x^*,u^*,t^*)=0$. If $\partial\mathbb{F}^{\MixCP}_{FB}(x^*,u^*,t^*)$ is nonsingular. Then the following statements hold:
	\begin{enumerate}
		\item[(i)] There exists $\bar{\eta}>0$ such that, if $ z_0 \in \textbf{B}\lf(z^*,\delta\rg)$ and $\eta_j\le \bar{\eta}$ for all $j$, then \textbf{Algorithm~1} is well defined and the sequence $\{z_j\}$ converges Q-linearly to $z^*$.
	
		\item[(ii)] If $\lim\limits_{ j\rightarrow\infty}\eta_j = 0$, then the sequence $\{z_j\}$ converges Q-superlinearly to $z^*$.
		
		\item[(iii)] If $\mathbb{F}^{\MixCP}_{FB}(x,u,t)$ is strongly semi-smooth at $z^*$, and if there exists $\tilde{\eta}>0$ such that $\eta_j\le \tilde{\eta}\|\mathbb{F}^{\MixCP}_{FB}(x_j,u_j,t_j)$ for all $j$, then the sequence $\{z_j\}$ converges Q-quadratic to $z^*$.
	\end{enumerate}
\end{theorem}
}

As mentioned earlier, A disadvantage of employing Newton's method for finding solutions to \eqref{eq:FBCFque} is that, it requires nonsingularity assumption of the Jacobian $\mathcal{A} $. Hence, it is worthy to present the widely-used \gls{lm} algorithm\cite{marquardt1963algorithm}. LM algorithm have least a linear rate of convergence without requiring all the Jacobian matrices in the iteration to be nonsingular. LM algorithm approximates the Hessian matrix by:

\[
	\mathcal{H}(z)=\mathcal{A}\tp (z)\mathcal{A}(z),
\]
and it approximates the gradient by:
\[
	\mathcal{G}(z)=\mathcal{A}\tp (z)\mathbb{F}^{\MixCP}_{FB}(z).
\]
Hence, its upgraded step will be
\[
	z_{j+1} = z_j - \left[\mathcal{A}\tp (z_j)\mathcal{A}(z_j) +\mu \mathbb{I}\right]^{-1}\mathcal{A}\tp (z_j)\mathbb{F}^{\MixCP}_{FB}(z_j).
\]

The parameter $\mu$ is used to prevent $d_j$ from being too large when $\mathcal{A}\tp (z_j)\mathcal{A}(z_j)$ is nearly singular. When $\mu$  equals to zero, the upgrading step is just the same as a Newton's method which uses the Hessian matrix for approximation. 

However, noting that the LM algorithm sacrifices calculation speed for the compatibility of singular Jacobian matrix, 
its number of iteration is probably greater than that of a Newton's method.  It should be worth noting that the settings of parameters also influence the calculation speed. A greater value of parameter $\mu$ will lead to a longer calculation time and larger number of iteration. LM algorithm is demonstrated as follows:

\begin{flushleft}
\textbf{Algorithm 2} (Semismooth Inexact Levenberg-Marquardt Method):
\end{flushleft}

\textbf{Input}: the initial point $z_0 = \bs x_0 \\u \\t\es \in \R^{k+\ell+1}$, the LM parameter $\mu_0\in\R_{+}$, and the tolerance $\eta_0\in\R_{+} $.

\textbf{Step 1}: Set $k = 0$.

\textbf{Step 2}: If $ \mathbb{F}^{\MixCP}_{FB}(z_j) =0 $, stop.

\textbf{Step 3}: Select an element $\mathcal{A}$ in the generalised Jacobian set $ \partial\mathbb{F}^{\MixCP}_{FB}(x,u,t)$, and find a direction $d_j \in \R^{k+\ell+1} $ such that
\[
	 \mathcal{A}(z_j)\tp \mathbb{F}^{\MixCP}_{FB}(z_j) + \left[\mathcal{A}\tp (z_j)\mathcal{A}(z_j) +\mu \mathbb{I}\right]d_j = r_j,
\]
where the residual vector $r_j\in \R^{k+\ell+1}$ satisfying
\[
	\|r_j\| \le \eta_j\|\mathcal{A}\tp (z_j)\mathbb{F}^{\MixCP}_{FB}(z_j)\|.
\]

\textbf{Step 4}: Choose $\eta_{j+1}\ge 0$ and $\mu_{j+1}\ge 0$; set $z_{j+1} := z_j + d_j $ and $j := j + 1$; go to \textbf{Step 2}.

Admittedly, Levenberg-Marquardt algorithm is also an efficient algorithm for finding solutions, because it also converges at least quadratically to a numerical solution.

{
\bigskip \begin{theorem}\label{LevenbergM} \emph{\cite{facchinei1997nonsmooth}}
	Let $\mathbb{F}^{\MixCP}_{FB}(x,u,t)$ be semi-smooth in $\textbf{B}\lf(z^*,\delta\rg)$, where $\delta >0$, and $z^*:= \bs x^* \\ u^* \\ t^* \es$ satisfies $\mathbb{F}^{\MixCP}_{FB}(x^*,u^*,t^*)=0$. If $\partial\mathbb{F}^{\MixCP}_{FB}(x^*,u^*,t^*)$ is nonsingular. Then the following statements hold:
	\begin{enumerate}
		\item[(i)] If $ z_0 \in \textbf{B}\lf(z^*,\delta\rg)$, then \textbf{Algorithm~2} is well defined and the sequence $\{z_j\}$ converges Q-linearly to $z^*$.
	
		\item[(ii)] If $\lim\limits_{ j\rightarrow\infty}\eta_j = 0$ and $\lim\limits_{ j\rightarrow\infty}\mu_j = 0$, then the sequence $\{z_j\}$ converges Q-superlinearly to $z^*$.
		
		\item[(iii)] If $\mathbb{F}^{\MixCP}_{FB}(x,u,t)$ is strongly semi-smooth at $z^*$, and if there exists $\tilde{\eta}>0$ and $\tilde{\mu}>0$ such that $\eta_j\le \tilde{\eta}\|\mathbb{F}^{\MixCP}_{FB}(x_j,u_j,t_j)$ and $\mu_j\le \tilde{\mu}\|\mathbb{F}^{\MixCP}_{FB}(x_j,u_j,t_j)$ for all $j$, then the sequence $\{z_j\}$ converges Q-quadratic to $z^*$.
	\end{enumerate}
\end{theorem}
}


\subsection{Reformulate to an unconstrained minimisation problem}

{
Another approach of solving the MixCP is to reformulate it to an unconstrained minimisation problem.} We will investigate the associated merit function of MixCP:
\index{merit function}
\begin{equation}\label{eq:amf}
	\theta^{\MixCP}_{FB}(x,u,t):=\frac{1}{2}\mathbb{F}^{\MixCP}_{FB}(x,u,t)\tp \mathbb{F}^{\MixCP}_{FB}(x,u,t).
\end{equation}

Obviously, based on \eqref{eq:FBCFque}, if there is a point $ \bs x^*\\u^*\\t^*\es$ such that
\begin{equation}\label{eq:Meritque}
	\theta^{\MixCP}_{FB}(x^*,u^*,t^*) = 0,
\end{equation}
then $\bs x^*\\u^*\\t^*\es$ is a solution to  $\MixCP$. FB C-function has its squared function $\psi_{FB}^2(a,b)$ to be continuously differentiable on $\mathbb{R}^2$ \cite{kanzow1994unconstrained}. Hence, it is easy to verify that the merit function $\theta^{\MixCP}_{FB}(x,u,t)$ is continuously differentiable if both $\widetilde{F}_1(x,u,t)$ and $\widetilde{F}_2(x,u,t)$ are. Since the merit function $\theta^{\MixCP}_{FB}(x,u,t)$ is nonnegative, if there exists a solution to \eqref{eq:Meritque}, then the point $\bs x^*\\u^*\\t^*\es$ will be a global minimiser of \eqref{eq:amf}. {Hence, the MixCP can be reformulated as the following unconstrained minimisation problem: 

\begin{equation}\label{eq:merit_min}
	\min\limits_{\bs x \\ u\\ t \es\in \R^{k+\ell+1}} \theta^{\MixCP}_{FB}(x,u,t).
\end{equation}

 In order to find a global minimiser of \eqref{eq:merit_min},} the general gradient of the merit function $\theta^{\MixCP}_{FB}(x,u,t)$ will be introduced:
\[
	\nabla \theta^{\MixCP}_{FB}(x,u,t) = \mathcal{A}\tp \mathbb{F}^{\MixCP}_{FB}(x,u,t),
\]
for any $ \mathcal{A} \in \partial\mathbb{F}^{\MixCP}_{FB}(x,u,t)$. It is not complicate to find the gradient of the merit function, but the difficulty lies in the lack of efficient tools to work out the non-convexity of $\theta^{\MixCP}_{FB}(x,u,t)$. For this purpose, we give the notion of stationary point. {The point $\bs x^*\\u^*\\t^*\es$ is said to be a stationary point of \eqref{eq:amf} if it satisfies the following inequality:
\begin{equation}\label{eq:VI}
	\lf\langle \bs x ~-~x^*\\ u~-~u^*\\ t~-~t^* \es, \nabla\theta^{\MixCP}_{FB}(x^*,u^*,t^*)\rg\rangle = 0, \quad \forall \bs x \\ u\\ t \es \in \R^{k +\ell +1}_+.
\end{equation}

The problem \eqref{eq:VI} is a variational inequality problem\cite{fukushima1992equivalent}. However, though $\bs x^*\\u^*\\t^*\es$ can be a stationary point of \eqref{eq:amf}, it does not guarantee that $\bs x^*\\u^*\\t^*\es$ is a global minimiser of \eqref{eq:amf}.} The discussion below is associated with the Jacobian $\partial\mathbb{F}^{\MixCP}_{FB}(x,u,t)$ at a stationary point of $ \theta^{\MixCP}_{FB}(x,u,t)$. Before introducing the notion of {\it FB regular\cite{FacchineiPang2003} } point, we define the following index sets:
\[
	\begin{array}{lcl}
	\C := \lf\{i: x_i \ge 0, (\widetilde{F}_1)_i\ge 0, x_i (\widetilde{F}_1)_i(x,u,t) = 0\rg\}, &  & complementarity~index, \\
	\mathcal{R} := \lf\{1, \dots, k\rg\}\setminus \C, &  & residual~index, \\
	\mathcal{P} := \lf\{i\in \mathbb{R} : x_i > 0 , (\widetilde{F}_1)_i(x,u,t) > 0\rg\}, &  & positive~index, \\
	\mathcal{N} := \mathcal{R}\setminus\mathcal{P}, &  & negative~index. \\
	\end{array}
\]

{
\textbf{Comment:} The FB regularity of the point $\bs x^*\\u^*\\t^*\es$ is defined by the Jacobian of $\widetilde{F}_1$ and $\widetilde{F}_2$ at $\bs x^*\\u^*\\t^*\es$. The motivation of introducing the notion of FB regularity is to avoid the case when $\bs x^*\\u^*\\t^*\es$ is a stationary point but the Jacobian $\partial\mathbb{F}^{\MixCP}_{FB}(x,u,t)$ is singular.  The property of FB regularity is tailored to the FB C-function as suggested by its name. In the following we will introduce the definition of FB regularity, and use Theorem \ref{thm:Sp} to show the connection between FB regularity and the solution to \eqref{eq:merit_min}.}

\bigskip \begin{definition}[FB regular]\label{FB_regular}

	A point $\bs x\\u\\t\es \in \R^k\times\R^{\ell}\times\R$ is called {\it FB-Regular} for the merit function $ \theta^{\MixCP}_{FB}$ if the Jacobian $ J_x\widetilde{F}_2(x,u,t) $  is nonsingular, and if for any $w \in \R^k$, $w\neq 0 $ with
	\[
	   w_i
		 \lf\{
		\begin{array}{ll}
			>0, \quad & if\;\; i\in \mathcal{P}, \\
			<0, \quad & if\;\; i\in \mathcal{N}, \qquad i\in\{ 1, \dots, k\},\\
			=0, \quad & if\;\; i\in \C. \\
		\end{array}
		\rg.
	\]
	there exists a nonzero vector $ v \in \R^k $ such that
	\begin{equation}\label{z_ineq}
		v_i
		\lf\{
		\begin{array}{ll}
			\ge0, \quad & if\;\; i\in \mathcal{P}, \\
			\le0, \quad & if\;\; i\in \mathcal{N}, \qquad i \in \{ 1, \dots, k\},\\
			=0  , \quad & if\;\; i\in \C. \\
		\end{array}
		\rg.
	\end{equation}
	and
	\begin{equation}\label{sch_c}
		w\tp\lf( M(x,u,t)/J_{\bs u\\ t\es}\widetilde{F}_2(x,u,t)\rg)v \ge 0,
	\end{equation}
	where
	\begin{equation}\label{m_m}
		M(x,u,t) :=
		\begin{pmatrix}
			 J_x\widetilde{F}_1(x,u,t)  &  J_{\bs u\\ t\es}\widetilde{F}_1(x,u,t)\\
			 J_x\widetilde{F}_2(x,u,t)  &  J_{\bs u\\ t\es}\widetilde{F}_2(x,u,t)
		\end{pmatrix}
		\in \R^{(k+\ell+1)\times(k+\ell+1)}
	\end{equation}
	and $ M(x,u,t)/J_{\bs u\\ t\es}\widetilde{F}_2(x,u,t) $ is the Schur complement of the block $ J_{\bs u\\t\es}\widetilde{F}_2(x,u,t)$ of the matrix $ M(x,u,t) $.
\end{definition}
\index{FB regular}

 Based on the original problem ESOCLCP, we can find the Jacobian of $ \widetilde{F}_1$ and $\widetilde{F}_2$ for MixCP:

\[
	J\widetilde{F}_1(x,u,t) :=
	\begin{pmatrix}
		J_x\widetilde{F}_1(x,u,t)  &  J_{\bs u\\ t\es}\widetilde{F}_1(x,u,t)
	\end{pmatrix}
	=
	\begin{pmatrix}
		A & \widetilde{B}
	\end{pmatrix},
\]
\[
	J\widetilde{F}_2 (x,u,t) :=
	\begin{pmatrix}
		J_x\widetilde{F}_2(x,u,t)  &  J_{\bs u\\ t\es}\widetilde{F}_2(x,u,t)
	\end{pmatrix}
	=
	\begin{pmatrix}
		\widetilde{C} & \widetilde{D}
	\end{pmatrix},
\]

where
\[
	\widetilde{B} :=
	\begin{pmatrix}
	B & Ae
	\end{pmatrix}, \qquad
		\widetilde{C} :=
	\begin{pmatrix}
		tC+ue^\top A \\
		0
	\end{pmatrix},
\]
and
\[
	\widetilde{D} :=
	\begin{pmatrix}
		\lf[A(x+te)+Bu+p\rg]\tp eI_{\ell} + ue\tp B + tD & Cx+2tCe+ue^\top Ae+Du \\
		-2u^\top & 2t
	\end{pmatrix}.
\]

In our case, if the Jacobian
\[
	 J_{\bs u\\ t\es}\widetilde{F}_2(x,u,t) = \widetilde{D}
\]
is nonsingular, then the Schur complement exists:

\begin{equation}\label{sch_eq}
	\left( M(x,u,t)/J_{\bs u\\ t\es}\widetilde{F}_2(x,u,t)\right) =
	A-\widetilde{B}\widetilde{D}^{-1}\widetilde{C}.
\end{equation}

The following theorem is based on \cite[Theorem 9.4.4]{FacchineiPang2003}. We made a slight modification in order to fit the context of the problem $ \MixCP(\widetilde{F}_1, \widetilde{F}_2,\R^k)$.  As there is no proof given in the source for this theorem, for the sake of completeness, a rewritten proof based on Definition \ref{FB_regular} is provided as follows:

\bigskip \begin{theorem}\label{thm:Sp}
	Let $\widetilde{F}_1:\R^k\times\R^{\ell}\times\R\rightarrow \R^{k}$ and $\widetilde{F}_2:\R^k\times\R^{\ell}\times\R\rightarrow \R^{\ell}\times\R$ be continuously differentiable. If $\bs x^*\\u^*\\t^*\es\in \R^k\times\R^{\ell}\times\R$ is a stationary point of $ \theta_{FB}^{\MixCP}$, then $\bs x^*\\u^*\\t^*\es$ is a solution to \eqref{eq:merit_min} if and only if $\bs x^*\\u^*\\t^*\es $ is an FB regular point of $ \theta_{FB}^{\MixCP}$ .
\end{theorem}

\begin{proof}
	Suppose that $\bs x^*\\u^*\\t^*\es\in \SMixCP(\widetilde{F}_1,\widetilde{F}_2,\R^k)$. It then follows that $\bs x^*\\u^*\\t^*\es$ is a global minimum and hence a stationary point of $ \theta_{FB}^{\MixCP}$. Thus, $(x^{*},\widetilde{F}_1(z^{*}))\in\C(\R^k_+)$, then we have $\mathcal{P}=\mathcal{N}=\emptyset$. Therefore the FB regularity of $x^{*}$ holds  since $x^{*}=x_{\C}$, and we cannot find a nonzero vector $x$ satisfying conditions \eqref{z_ineq}.
	Conversely, suppose that $x^{*}$ is FB regular and $ \bs x^*\\u^*\\t^*\es$ is a stationary point of $ \theta_{FB}^{\MixCP}$. It follows that $\nabla \theta_{FB}^{\MixCP} = 0$, i.e.:
	\[
		\mathcal{A}^\top \mathbb{F}_{FB}^{\MixCP} =
		\begin{pmatrix}
			D_a+\lf(J_x\widetilde{F}_1(x^*, u^*, t^*)\rg)\tp D_b & \lf(J_x\widetilde{F}_2(x^*, u^*, t^*)\rg)\tp\\
			\lf(J_{\bs u\\ t\es}\widetilde{F}_1(x^*, u^*, t^*)\rg)\tp D_b & \lf(J_{\bs u\\ t\es}\widetilde{F}_2(x^*, u^*, t^*)\rg)\tp
		\end{pmatrix}
		\mathbb{F}_{FB}^{\MixCP}=0, 
	\]
	where $D_a$ and $D_b$ are diagonal matrices with
	\[
		(D_a)_{ii}= \left\{
			\begin{array}{ll}
				 \frac{x_i}{\sqrt{x_i^2+( \widetilde{F}_1)_i^2(x,u,t)}}-1,  &
				 if\; \lf(x_i, (\widetilde{F}_1)_i(x,u,t)\rg) \neq (0,0),\\
				-1, &
				if\; \lf( x_i, (\widetilde{F}_1)_i(x,u,t)\rg) = (0,0),
			\end{array}\right.
		\qquad i\in\{1, \dots , k\},
	\]

	\[
		(D_b)_{ii}=\left\{
		\begin{array}{ll}
			 \frac{( \widetilde{F}_1)_i(x,u,t)} {\sqrt{x_i^2+( \widetilde{F}_1)_i^2(x,u,t)}}-1,   &
				 if\; \lf(x_i, (\widetilde{F}_1)_i(x,u,t)\rg) \neq (0,0),\\
			 -1,  &
				if\; \lf( x_i, (\widetilde{F}_1)_i(x,u,t)\rg) = (0,0),
		\end{array}\right.
		 \qquad  i\in\{1, \dots , k\}.
	\]
	Hence, for any $w \in \R^k\times\R^{\ell}\times\R$, we have
	\begin{equation}\label{p_thm3_0}
		w\tp
		\begin{pmatrix}
			D_a+\lf(J_x\widetilde{F}_1(x^*, u^*, t^*)\rg)\tp D_b & \lf(J_x\widetilde{F}_2(x^*, u^*, t^*)\rg)\tp\\
			\lf(J_{\bs u\\ t\es}\widetilde{F}_1(x^*, u^*, t^*)\rg)\tp D_b & \lf(J_{\bs u\\ t\es}\widetilde{F}_2(x^*, u^*, t^*)\rg)\tp
		\end{pmatrix}
		\mathbb{F}_{FB}^{\MixCP}=0.
	\end{equation}
	
	Suppose that $\bs x^*\\u^*\\t^*\es$ is not a solution to $\MixCP $, we have that the index set $\mathcal{R}$ is not empty.
	Define $v:= D_b\mathbb{F}_{FB}^{\MixCP}$, we have
	\[
		v_{\mathcal{C}} =0,
		\qquad
		v_{\mathcal{P}} >0,
		\qquad
		v_{\mathcal{N}} <0.
	\]

	Take $w$ with
	\[
		w_{\mathcal{C}} =0,
		\qquad
		w_{\mathcal{P}} >0,
		\qquad
		w_{\mathcal{N}} <0.
	\]
	From the definition of $D_a$ and $D_b$, we know that $D_a\mathbb{F}_{FB}^{\MixCP}$ and $D_b\mathbb{F}_{FB}^{\MixCP}$ have the same sign. Therefore
	\begin{equation}\label{p_thm3_1}
		w\tp (D_a\mathbb{F}_{FB}^{\MixCP}) = w\tp_{\mathcal{C}}(D_a\mathbb{F}_{FB}^{\MixCP})_{\mathcal{C}} + w\tp_{\mathcal{P}}(D_a\mathbb{F}_{FB}^{\MixCP})_{\mathcal{P}} +w\tp_{\mathcal{N}}(D_a\mathbb{F}_{FB}^{\MixCP})_{\mathcal{N}} > 0.
	\end{equation}
	
	By regularity of $J\widetilde{F}_1(z)\tp $, we conclude
	\begin{equation}\label{p_thm3_2}
		w\tp J\widetilde{F}_1(z)\tp (D_a\mathbb{F}_{FB}^{\MixCP})= w\tp J\widetilde{F}_1(z)\tp w \ge 0.
	\end{equation}
	
	The inequalities \eqref{p_thm3_1} and \eqref{p_thm3_2} together contradict condition \eqref{p_thm3_0}. Hence $ \mathcal{R} = \emptyset$. It means that $z^{*}$ is the solution to $\MixCP(\widetilde{F}_1,\widetilde{F}_2,\R^k)$.

\hfill$\square$\par \end{proof}

If the Schur complement $ M(x,u,t)/J_{\bs u\\ t\es}\widetilde{F}_2(x,u,t)\in\R^{k\times k}$ is a signed $S_0$ matrix, then $\bs x\\u\\t\es$ is an FB regular point of $ \theta_{FB}^{\MixCP}$. Before proving this assertion, the definition of signed $S_0$ matrix is given as follows:

\bigskip \begin{definition}[Signed $S_0$ matrix]\label{signed_s0}
	Let $\widetilde{F}_1:\R^k\times\R^{\ell}\times\R\rightarrow \R^{k}$ and $\widetilde{F}_2:\R^k\times\R^{\ell}\times\R\rightarrow \R^{\ell}\times\R$ be continuously differentiable. We say that the Schur Complement $ M(x,u,t)/J_{\bs u\\ t\es}\widetilde{F}_2(x,u,t) $ is a {\it signed $S_0$ matrix} if
	\begin{equation}\label{ss0_mat}
		\Xi := \Lambda\lf(M(x,u,t)/J_{\bs u\\ t\es}\widetilde{F}_2(x,u,t)\rg) \Lambda
	\end{equation}
	is a $S_0$ matrix (see Definition \ref{def:s0_mat}), where $\Lambda\in\R^{k\times k}$ is the diagonal matrix whose diagonal entries $\lambda_i$, $i \in\{1, \dots, k\}$ satisfy
	\begin{equation}\label{deflm}
		\lambda_i := \lf\{
		\begin{array}{rl}
			1, \quad & if\;\; i\in \mathcal{P}, \\
			-1, \quad & if\;\; i\in \mathcal{N}, \\
			0, \quad & if\;\; i\in \C. \\
		\end{array}
		\rg.
	\end{equation}
\end{definition}
\index{Signed $S_0$ matrix}

\bigskip \begin{theorem}\label{SP2}
	Let $\widetilde{F}_1:\R^k\times\R^{\ell}\times\R\rightarrow \R^{k}$ and $\widetilde{F}_2:\R^k\times\R^{\ell}\times\R\rightarrow \R^{\ell}\times\R$ be continuously differentiable functions, and suppose the Jacobian matrix $J_{\bs u\\ t\es}\widetilde{F}_2$ is nonsingular. If the Schur complement $ M(x^*,u^*,t^*)/J_{\bs u\\ t\es}\widetilde{F}_2(x^*,u^*,t^*) $ is a signed $S_0$ matrix, then $\bs x^*\\u^*\\t^*\es$ is an FB regular point of $ \theta_{FB}^{\MixCP}$.
\end{theorem}

\begin{proof}
	Suppose that  $ M(x^*,u^*,t^*)/J_{\bs u\\ t\es}\widetilde{F}_2(x^*,u^*,t^*) $ is a signed $S_0$ matrix, denoted by
	\[
		M_2 = M(x^*,u^*,t^*)/J_{\bs u\\ t\es}\widetilde{F}_2(x^*,u^*,t^*),
	\]
then Definition \ref{signed_s0} implies:
	\begin{equation}\label{intermed1}
		 \lf(\Lambda  M_2  \Lambda\rg) u \ge 0.
	\end{equation}
	where $\Lambda$ is defined in \eqref{ss0_mat}. It is clear that there exists $0\neq u \ge 0$ such that $u$ is a solution to \eqref{intermed1}. Denote by $\lf(M_2\rg)_i$ the i-th row of matrix $M_2$, we can rewrite \eqref{intermed1} as
	\begin{equation}\label{intermed2}
	\begin{array}{rl}
		1\cdot \lf(M_2\rg)_i \Lambda u \ge 0 ,& \qquad \forall i \in \mathcal{P},\\
		-1\cdot \lf(M_2\rg)_i \Lambda u \ge 0,& \qquad \forall i \in \mathcal{N},\\
		0\cdot \lf(M_2\rg)_i \Lambda u = 0   ,& \qquad \forall i \in \C.
	\end{array}
	\end{equation}
	Take any $w\in \R^k\setminus \{0\}$ with
	\[
	   w_i
		 \lf\{
		\begin{array}{ll}
			>0, \quad & if\;\; i\in \mathcal{P}, \\
			<0, \quad & if\;\; i\in \mathcal{N}, \qquad i \in\{1, \dots, k\},\\
			=0, \quad & if\;\; i\in \C, \\
		\end{array}
		\rg.
	\]
	multiplies with inequalities/equality \eqref{intermed2}, then there exists a nonzero vector $v = \Lambda u$ such that
	\begin{equation}\label{intermed3}
	\begin{array}{rl}
		w_i\cdot \lf(M_2\rg)_i v \ge 0 ,& \qquad \forall i \in \mathcal{P},\\
		w_i\cdot \lf(M_2\rg)_i v \ge 0,& \qquad \forall i \in \mathcal{N},\\
		w_i\cdot \lf(M_2\rg)_i v = 0   ,& \qquad \forall i \in \C.
	\end{array}
	\end{equation}
	 hold. Noting that
	\[
	v := \lf\{
		\begin{array}{rl}
			u_i, \quad & if\;\; i\in \mathcal{P}, \\
			-u_i, \quad & if\;\; i\in \mathcal{N},\\
			0, \quad & if\;\; i\in \C. \\
		\end{array}
		\rg.
	\]

	 Summing up the three iequalities/equality above in \eqref{intermed3}, we get
	\begin{equation}\label{eq:s02fb}
		\sum_{i\in\mathcal{P}}w_i\lf(M_2 v\rg)_i +\sum_{i\in\mathcal{N}}w_i\lf(M_2 v\rg)_i + \sum_{i\in\mathcal{C}}w_i\lf(M_2 v\rg)_i\ge 0.
	\end{equation}

	That is equivalent to
	\[
		w\tp\lf(M_2\rg)v \ge 0.
	\]

	Hence, $\bs x^*\\u^*\\t^*\es$ is an FB regular point of $ \theta_{FB}^{\MixCP}$.
\hfill$\square$\par \end{proof}

\bigskip \begin{example}
	Taking the notation in Theorem \ref{SP2}. Suppose that $M_2$ is a signed $S_0$ matrix. Let the Schur complement be
	\[
		M_2 =
		\begin{pmatrix}
			m_{11} & m_{12} & m_{13} \\
			m_{21} & m_{22} & m_{23} \\
			m_{31} & m_{32} & m_{33}
		\end{pmatrix}.
	\]
	with the index set $\mathcal{P} = \{1\}$,  $\mathcal{N} = \{2\}$,  and $\C = \{3\}$. There exists a nonzero vector $u = \lf(u_1, u_2, u_3\rg)\tp$ with $0\neq u\ge 0$ such that

	\begin{align*}
		&
		 \begin{pmatrix}
			1 &  0 & 0\\
			0 & -1 & 0\\
			0 &  0 & 0
		\end{pmatrix}
		\begin{pmatrix}
			m_{11} & m_{12} & m_{13} \\
			m_{21} & m_{22} & m_{23} \\
			m_{31} & m_{32} & m_{33}
		\end{pmatrix}
		\begin{pmatrix}
			1 &  0 & 0\\
			0 & -1 & 0\\
			0 &  0 & 0
		\end{pmatrix}
			\begin{pmatrix}
		u_1 \\
		u_2 \\
		u_3
		\end{pmatrix}\\
		=& \begin{pmatrix}
			 1\cdot m_{11} &  1\cdot m_{12} & 1\cdot m_{13} \\
			-1\cdot m_{21} & -1\cdot m_{22} & -1\cdot m_{23} \\
			0\cdot m_{31} & 0\cdot m_{32} & 0\cdot m_{3}
		\end{pmatrix}
			\begin{pmatrix}
		 u_1 \\
		-u_2 \\
		 0
		\end{pmatrix}
		\ge \mathbf{0}.
	\end{align*}

	Denote $v = \lf(u_1, -u_2, 0\rg)\tp$, the above linear system can be written as:
	\[
		\begin{array}{rl}
			 1\cdot \lf( m_{11}, m_{12},m_{13}\rg)v &\ge 0,\\
			-1\cdot \lf( m_{21}, m_{22},m_{23}\rg)v &\ge 0,\\
			 0\cdot \lf( m_{31}, m_{32},m_{33}\rg)v & =  0.
		\end{array}
	\]

	Take  $w=\lf(w_1, -w_2,0\rg)$ with $w_1$, $w_2\in \R_+$, the following linear system holds:
	\[
		\begin{array}{rl}
			 w_1\cdot \lf( m_{11}, m_{12},m_{13}\rg)v &\ge 0,\\
			-w_2\cdot \lf( m_{21}, m_{22},m_{23}\rg)v &\ge 0,\\
			 0  \cdot \lf( m_{31}, m_{32},m_{33}\rg)v & =  0.
		\end{array}
	\]

	Hence, we can say for any vector $w = \bs w_1\\-w_2\\0\es $, there exists a vector  $v = \bs u_1\\-u_2\\0\es $ such that
	\[
		 w\tp M_2 v \ge 0
	\]
	holds. Hence,$ (x^*)$ is a FB regular point of the merit function $\theta(x)$.
\end{example}

The Theorem \ref{SP2} has proved that the signed $S_0$ property of the Schur complement $ M(x,u,t)/J_{\bs u\\ t\es}\widetilde{F}_2(x,u,t) $ is a sufficient condition for $x$ to be an FB regular point.
This condition outperforms the FB regularity condition, because that the verification of the $S_0$ property of the
matrix \eqref{ss0_mat} can be accomplished by simple linear programming. On the other hand, verifying FB regularity will be in general more complex and computationally expensive. Together with Theorem \ref{thm:Sp} and \ref{SP2}, the conclusion of the signed $S_0$ property enables us to find a solution to $\MixCP$ with algorithms.

\begin{flushleft}
\textbf{Algorithm 3} (FB line search method):
\end{flushleft}

\textbf{Input}: the initial point $z_0 =(x,u,t) \in \R^{k}\times\R^{\ell}\times\R$, $\rho >0$, $\gamma\in (0,1)$, and the tolerance $r\in\R_{++}$.

\textbf{Step 1}: Set $k = 0$.

\textbf{Step 2}: If $\|\nabla\theta^{\MixCP}_{FB}(z_j)\| \le r$, then STOP.

\textbf{Step 3}: Select an element $\mathcal{A}\in \partial\mathbb{F}^{\MixCP}_{FB}(x,u,t)$, and find a direction $d_j \in \R^{k}\times\R^{\ell}\times\R$ such that
\begin{equation}\label{algo3}
	\mathbb{F}^{\MixCP}_{FB}(z_j) + \mathcal{A}\tp (z_j)d_j = 0.
\end{equation}

If the system \eqref{algo3} is not solvable or if the condition
\begin{equation}\label{algo3_1}
	\nabla \theta^{\MixCP}_{FB}(z_j)d_j \le -\rho \|d_j\|
\end{equation}
is not satisfied, reset $d_j := -\nabla \theta^{\MixCP}_{FB}(z_j)$.

\textbf{Step 4}: Find the smallest nonnegative integer $i_j$ such that, with $i = i_j$, we have
\[
	\theta^{\MixCP}_{FB}(z_j + 2^{-i}d_j) \le \theta^{\MixCP}_{FB}(z_j) + \gamma 2^{-i}\theta^{\MixCP}_{FB}(z_j)\tp d_j;
\]
set $\tau_j := 2^{-i_j}$.

\textbf{Step 5}: Set $z_{j+1} := z_j + d_j $ and $j := j + 1$, go to \textbf{Step 2}.

In general, matrix $\mathcal{A}$ in \textbf{Step 3} of this algorithm is not necessary to be a generalised Jacobian of $\mathbb{F}^{\MixCP}_{FB}$ at the iterate $z_j$, but it still is required to be a nonsingular matrix. In this algorithm, \textbf{Step 2} is to examine whether an iterate $z_j$ is a stationary point of $ \theta^{\MixCP}_{FB}$ or not; whereas the condition \eqref{algo3_1} and \textbf{Step 4} are used to ensure that the iterate $z_j$ is a FB regular point of $ \theta^{\MixCP}_{FB}$.

\section{A numerical example}

In this section, we will illustrate a numerical example corresponding to item (iv) of Proposition \ref{prop:cs-esoc}. Let $L(3,2)$ and $M(3,2)$ be an ESOC and its dual cone defined by \eqref{esoc} and \eqref{desoc}, respectively. Denote
\[
	z=\bs x\\u\es\in\R^3\times\R^2,\quad \hat z=\bs\hat{x}\\u\es :=\bs x-\|u\|e\\u\es\in\R^3\times\R^2,
\]
and 
\[
	\tilde{z}=\bs \tilde{x}\\u\\t\es:=\bs x-t\\u\\t\es\in\R^3\times\R^2\times\R.
\]

Repeat for convenience an ESOCLCP defined by an extended second order cone $L\in\R^3\times\R^2$ and a linear function $F:\R^3\times\R^2\rightarrow \R^3\times\R^2$, $F(x,u) = T\bs x\\u\es +r$, is:
\[
	\LCP(F,L)\left\{
	\begin{array}{l}
		Find \; \bs x\\u\es\in L,\; such \; that\\
		F(x,u) \in M \; and \; \langle \bs x\\u\es, F(x,u) \rangle = 0.
	\end{array}
	\right.
\]
where $T=\bs A & B\\C & D \es$, $r = \bs p \\ q\es$, with $A\in\R^{3\times 3}$, $B\in\R^{3\times 2}$, $C\in\R^{2\times 3}$, $D\in\R^{2\times 2}$, $p\in\R^3$, and $q\in\R^2$.
The solution to ESOCLCP is equivalent to the solution to a corresponding MixCP converted by employing item (vi) of Theorem \ref{thm:main}. For convenience the $\MixCP(F_1,F_2,\R_+^3)$, defined by $F_1$, $F_2$, and $\R_+^3$, is provided:
\[
	\MixCP(F_1,F_2,\R_+^3):\left\{
	\begin{array}{l}
		Find \; \bs x\\u\\t\es\in\R^3\times\R^2\times\R^1, \; such \; that\\
		\widetilde{F}_2(x,u,t)=0, \; and\; (x,\widetilde{F}_1(x,u))\in\C(\R_+^3).
	\end{array}
	\right.
\]
where
\[
	\widetilde{F}_1(x,u,t)=A(x+te)+Bu+p
\]
 and
\[
	\widetilde{F}_2(x,u,t) =
	\begin{pmatrix}
		& \lf(tC+ue\tp A\rg)(x+te)+ue\tp(Bu+p)+t(Du+q) \\
		& t^2 - \|u\|^2
	\end{pmatrix}.
\]

As the propose of this section is not comparing the efficiency of the algorithms, we will only employ the Levenberg-Marquardt algorithm (\textbf{Algorithm 2}). Our objective is to solve the following FB-based equation formulation \eqref{FBform}:

\[
	\mathbb{F}^{\MixCP}_{FB}(x,u,t) =
	\begin{pmatrix}
		 \psi_{FB}\lf(x_1,(\widetilde{F}_1)_1(x,u,t)\rg) \\
		 \vdots \\
		 \psi_{FB}\lf(x_k,(\widetilde{F}_1)_k(x,u,t)\rg) \\
		 \widetilde{F}_2(x,u,t)
	\end{pmatrix} = 0.
\]

We set the tolerance parameter $r = 10^{-7}$, the initial LM parameter $\mu = 10^{-2}$ and multiply by $10^{-1}$ after each iteration.

Consider
\[
	T=\begin{pmatrix} A & B\\C & D \end{pmatrix} =
	\lf(
	\begin{array}{rrrrr}
		41  &  -3 &  -31 &   18  &  19  \\
		28  &  22 &  -33 &   25  & -29  \\
	   -23  & -29 &   11 &  -21  & -43  \\
		-9  & -31 &  -20 &  -12  &  47  \\
		-8  &  46 &   50 &  -22  &  21
	\end{array}
	\rg), \quad 
	r=\lf(\begin{array}{r}p\\ q\end{array}\rg) = \lf(
	\begin{array}{r}
	   -26  \\
		 4  \\
		23  \\
		44  \\
	   -19
	\end{array}\rg),
\]

By using \textbf{Algorithm 2}, the sequence $\{z_j\}$ converges to a numerical solution in 11 iterations (Table \ref{tab:esoc_ex1}).
\begin{table}[!ht]
\centering
\begin{threeparttable}
\begin{tabular}{c c c r r}
\toprule
  Iteration & The value of $\|\nabla\theta^{\MixCP}_{FB}(z_j)\|$ & Optimality value of $\theta^{\MixCP}_{FB}(z_j)$ & $\mu$ & $d$ \\[0.4ex]
\midrule
	 0 &  1.1760e+09 & 2.43e+08  &	   0.01  &			   \\
	 1 &  3.9793e+07 & 7.18e+06  &	  0.001  &	  86.1516  \\
	 2 &  1.0348e+07 & 1.50e+06  &	  1e-04  &	  15.3125  \\
	 3 &  640568.382 & 2.06e+05  &	  1e-05  &	  9.24632  \\
	 4 &  37188.3244 & 2.96e+04  &	  1e-06  &	  4.70891  \\
	 5 &  4932.01875 & 5.83e+03  &	  1e-07  &	  2.09599  \\
	 6 &  394.966198 & 975.2398  &	  1e-08  &	  1.15061  \\
	 7 &  19.3492349 & 115.2587  &	  1e-09  &	  0.51351  \\
	 8 &  0.20392831 & 12.92491  &	  1e-10  &	  0.23590  \\
	 9 &  1.5122e-06 & 0.034555  &	  1e-11  &	  0.03290  \\
	10 &  7.9130e-15 & 2.00e-06  &	  1e-12  &	  5.3e-05  \\
	11 &  8.2224e-15 & 3.95e-15  &	  1e-13  &	  1.1e-06  \\
\bottomrule
\end{tabular}
\begin{tablenotes}
\footnotesize
\item Note: This table shows the iteration of $\{z_j\}$ converges to a numerical solution by using \textbf{Algorithm 2}. It takes 11 iterations to reach a numerical solution with tolerance $r = 10^{-7}$.
\end{tablenotes}
\end{threeparttable}
\caption{Numerical Example: the iteration of $\{z_j\}$} \label{tab:esoc_ex1}

\end{table}

The solution to the $\MixCP$ is $\tilde{z}^* = \lf(\tilde{x}^*, u^*, t^*\rg)\tp = \lf(\frac{781}{641}, 0, \frac{999}{1328}, \frac{333}{2693}, -\frac{619}{2428}\rg)\tp$. Verifying the complementarity: 

\[
	\tilde{x}^* = \lf(\frac{781}{641}, 0, \frac{999}{1328}\rg)\tp \ge 0,  \qquad \widetilde{F}_1(\tilde{z}^*) = \lf(0, \frac{8349}{292}, 0\rg)\tp \ge 0,
\]
\[
	\langle \tilde{x}^*, \widetilde{F}_1(\tilde{z}^*)\rangle = 0.
\]
Therefore we have $(\tilde{x}^*, \widetilde{F}_1(\tilde{z}^*))\in \C(\R^3_+)$. By the item (vi) of Theorem \ref{thm:main}, we get the solution to $\LCP(T,r,L(3,2))$ from $\tilde{z}^*$, that is $z^* = \lf(\tilde{x}^*+t^*, u^*\rg)\tp = \lf(\frac{428}{285}, \frac{325}{1147}, \frac{1716}{1657}, \frac{333}{2693}, -\frac{619}{2428}\rg)\tp$. 

For verifying the complementarity, we show that $z^* \in L(3,2)$ because $\tilde{x}^*+t^*\ge \frac{325}{1147} = \sqrt{\frac{333}{2693}^2 + \frac{619}{2428}^2} = \|u^*\|$; and
\[
	F(x,u) =
	\lf(
	\begin{array}{rrrrr}
		41  &  -3 &  -31 &   18  &  19  \\
		28  &  22 &  -33 &   25  & -29  \\
	   -23  & -29 &   11 &  -21  & -43  \\
		-9  & -31 &  -20 &  -12  &  47  \\
		-8  &  46 &   50 &  -22  &  21
	\end{array}
	\rg)\lf(
	\begin{array}{c}
		 \frac{428}{285}  \\
		 \frac{325}{1147}  \\
		 \frac{1716}{1657}  \\
		 \frac{333}{2693}  \\
		-\frac{619}{2428}
	\end{array}\rg)
	+ \lf(
	\begin{array}{r}
	   -26  \\
		 4  \\
		23  \\
		44  \\
	   -19
	\end{array}\rg)
	= \lf(
	\begin{array}{r}
		0  \\
		\frac{8349}{292}  \\
		0  \\
		-\frac{3943}{316}  \\
		\frac{4039}{157}
	\end{array}\rg).
\]
We obtain that $F(x,u) \in M(3,2)$ as $ 0+\frac{8349}{292}+0 \ge \sqrt{(-\frac{3943}{316})^2 + \frac{4039}{157}^2} = \frac{8349}{292}$, and hence $\lf(z^*, F(x,u)\rg) \in {\cal C}\lf(L(3,2)\rg)$. Therefore, we can confirm that this is a solution to the problem $\LCP(T,r,L(3,2))$.

\section{Conclusions and comments}

In this chapter, we study the  linear complementarity problem on extended second order cones (ESOCLCP). Our main result is Theorem \ref{thm:main}, which discusses the connections between an ESOCLCP and mixed (implicit, mixed implicit) complementarity problems. Under some mild conditions, we can rewrite an ESOCLCP to a mixed complementarity problems (MixCP) on the nonnegative orthant. In the new formulation, both $\widetilde{F}_1(\tilde{x},u,t)$ and $\widetilde{F}_2(\tilde{x},u,t)$ are smooth functions, which simplifies the process of finding solutions to MixCP.
The conversion from an ESOCLCP to a MixCP on the nonnegative orthant reduces the complexity of finding solutions to the original problem. 

The process of solving MixCP is straightforward. We introduced the FB C-function to reformulate the complementarity problem. By the FB C-function, the complementarity problem can be reformulated to either a system of nonlinear equations or an unconstrained minimisation problem. For the nonlinear equation approach, we introduced proposition for the nonsingularity of the Jacobian. The semi-smooth inexact Newton method and the Levenberg-Marquardt method are illustrated. Further, we provided theorems to verify the rate of convergence of both algorithms.
For the minimisation approach, we introduced and proved Theorem \ref{thm:Sp} and Theorem \ref{SP2} for the difficulty of non-convexity. A point is sufficiently to be a solution to a MixCP if it satisfies specific conditions related to stationarity, FB regularity (Theorem \ref{thm:Sp}), and Signed $S_0$ property (Theorem \ref{SP2}). These theorems can be used to determine whether a point is a solution to the MixCP or not. Based on the above, a solution to MixCP will be equivalent to a solution to the corresponding ESOCLCP. We use FB Line Search Algorithm for the minimisation approach. In the final section, we illustrate a numerical example corresponding to item (iv) of Proposition \ref{prop:cs-esoc} and item (vi) of Theorem \ref{thm:main}.

\chapter{Stochastic Linear Complementarity Problems on Extended Second Order Cones} \label{cht:stesoc}

In many practical situations, uncertainty is a common and realistic problem that results from inaccurate measurement or stochastic variation of data such as price, capacities, loads, etc. In fact, the inaccuracy or uncertainty of these real-world data are inevitable. When these data are applied as parameters in mathematical models, the constraints of models may be violated because of their stochastic characters. These violations may finally cause some difficulties that the optimal solutions obtained from the stochastic data are no longer optimal, even infeasible. Amongst approaches proposed for modeling uncertain quantities, the stochastic models outstand because of their solid mathematical foundations, theoretical richness, and sound techniques of using real data. Complementarity problems imbedded with stochastic models occur in many areas such as finance, telecommunication and engineering. Hence, considering $\LCP$ with uncertainty will be meaningful for practical treatments. If partial or all of the coefficients in the $\LCP$ are uncertain, the $\LCP$ will be turned into a \gls{slcp}, which is firstly introduced by Chen and Fukushima \cite{chen2005expected}. Articles about SLCP can be found in \cite{fang2007stochastic,gurkan1999sample,lin2009combined,chen2011cvar}.

Even though the fact that only limited number of results have been obtained on the stochastic complementarity problems, there are still some meaningful results. One of them is the CVaR (conditional value-at-risk, which is also called expected shortfall) minimisation reformulation of stochastic complementarity problem \cite{xu2014cvar}. In this chapter, \gls{sesoclcp} will be studied. Based on the results in previous chapter, a method of finding solutions to S-ESOCLCP will be elaborated, and a numerical example will be presented.

\section{Problem formulation}

Let $(\Omega, \mathcal{F},\mathcal{P})$ be a probability space\index{Probability space} defined by:
\begin{enumerate}
	\item $\Omega\subseteq \R^{n}$, the sample set of possible outcomes\index{Sample set of possible outcomes};
	\item {$\mathcal{F}\subseteq 2^{\Omega}$, a $\sigma$-algebra generated by $\Omega$ (all subsets of $\Omega$)\index{$\sigma$-algebra}}; and
	\item $\mathcal{P}: \mathcal{F}\rightarrow [0,1]$, a function {map} from events to probabilities.
\end{enumerate}

The following is the definition of a \gls{stcp}:
\index{Complementarity problem!Stochastic complementarity problem}

\bigskip \begin{definition}[Stochastic complementarity problem]\label{defn:slcp}
Given a random vector valued function $F(x,\omega): \R^n \times \Omega \rightarrow \R^n$, where $\omega\in\Omega$ is an $n$-dimensional random vector. A {\it stochastic complementarity problem} is defined by
\begin{equation}\label{eq:scp}
	SCP(F,\R^n_+,\omega)\left\{
		\begin{array}{l}
			Find \; x\in \R^m_+,\; such \; that\\ F(x,\omega)\ge 0, \quad x\tp F(x,\omega)=0,\quad\omega\in\Omega, \quad a.s.
		\end{array} \right.
	\end{equation}
\end{definition}

The abbreviation \gls{as} means  $F(x,\omega)\ge 0$ and $\quad x\tp F(x,\omega)=0$ hold almost surely for any $\omega\in \Omega$.

If $F(x,\omega)$ is a linear function of the form $F(x,\omega) = T(\omega)x+r(\omega)$, then we call problem \eqref{eq:scp} a \emph{ stochastic linear complementarity problem (SLCP) }\index{Complementarity problem!Stochastic linear complementarity problem}, specifically:
\begin{equation}\label{eq:slcp}
	SLCP(T(\omega),r(\omega),\R^n_+,\omega)\left\{
	\begin{array}{l}
		Find \; x\in \R^n_+,\; such \; that\\ T(\omega)x+r(\omega) \ge 0, x\tp (T(\omega)x+r(\omega)) = 0, \omega\in \Omega, \quad a.s.
	\end{array}
	\right.
\end{equation}

In this chapter, we assume that the coefficients $T(\omega)$ and $r(\omega)$ are measurable functions of $\omega$ with the following property:
\[
	\E[\|T(\omega)\tp T(\omega)\|] < \infty\quad and \quad \E[\|r(\omega)\|] < \infty
\]
where $\E[\cdot]$ represents the expected value of the random vector in the square bracket.

It should be mentioned that if the possible outcome set $\Omega$ contains only one single realisation (and this unique outcome definitely happens), problem \eqref{eq:slcp} will degenerate to problem \eqref{eq:lcp}.

The stochastic linear complementarity problems are very useful in solving practical problems. However, because of the existence of the random vector $\omega$ in the function $F (x, \omega)$, it is very difficult and sometimes impossible to find a solution $x$ satisfying all possible outcomes of $\omega\in \Omega$. One plausible idea to improve the viability of finding a solution to SLCP is to associate the problems with probability models, and then persuasive solutions to SLCP are obtainable by finding the solutions to the associated probability models.

Xu and Yu \cite{xu2014cvar} summarised 6 different probability models for finding solutions to SLCP:


\begin{enumerate}
	\item[(i)] \textbf{\gls{ev} method, introduced by G{\"u}rkan et. al in \cite{gurkan1999sample}}. By using the expectation value $\E[F(x,\omega)]$ to replace the stochastic term $F(x,\omega)$, this method ultimately reformulates \eqref{eq:slcp} to \eqref{eq:lcp}.

	\item[(ii)] \textbf{\gls{erm} method, introduced by Chen and Fukushima \cite{chen2005expected}}. This method minimises the expectation of the square norm of the residual $\Phi(x,\omega)$ defined by the following C-function:
	\begin{equation}\label{CF_max}
		\min_{x\in\R^n_+}\E \left[ \|\Phi(x,\omega)\|^2\right]
	\end{equation}
	where $\Phi:\R^n \times\Omega\rightarrow\R^n$ is a multi dimensional C-function defined as
	\[
		\Phi(x,\omega):= \lf(\phi\lf(x_1,F_1(x,\omega)\rg),\dots,\phi\lf(x_m,F_m(x,\omega)\rg)\rg)\tp .
	\]
	where $\phi : \R\times\R \rightarrow\R$ can be any scalar C-function satisfying:
	\[
		\phi(a,b) = 0 \quad \Leftrightarrow \quad a\ge 0,\quad b\ge 0, \quad ab = 0.
	\]

	\item[(iii)] \textbf{\gls{smpec} reformulation, introduced by Lin and Fukushima\cite{lin2009combined, lin2009solving, mataramvura2008risk}}. This method highlights a recourse variate $z(\omega)$ to compensate the violation of complementarity in \eqref{eq:slcp} for some outcomes of $\omega\in\Omega$, then it reformulates \eqref{eq:slcp} to the following model:
	\begin{equation}\label{eq:SMPEC}
		\begin{array}{lll}
		& \min\limits_x & \E\left[\eta\lf(z(\omega)\rg)\right] \\& s.t. & 0\le x\perp \lf(F(x,\omega) + z(\omega)\rg) \ge 0, \\ &		  & z(\omega) \ge 0, \omega\in \Omega\quad a.s.,
			\end{array}
	\end{equation}
	where $\eta(z) = e^{tp}z$. Ambiguous solutions to SCP can be obtained by minimising the objective function in \eqref{eq:SMPEC}, i.e. the expected value of the compensation to the violation of complementarity in \eqref{eq:slcp}.

	\item[(iv)] \textbf{\gls{sp} {reformulation} \cite{wang2010stochastic}}. Problem \eqref{eq:slcp} is reformulated to the following:
		\[
			\begin{array}{lll}
				& \min\limits_{x} & \E\left[\|\lf(x \circ F(x,\omega)\rg)_+\|^2 \right] \\
				& s.t. & F(x,\omega) \ge 0, \quad \omega \in \Omega \quad a.s. \\
				&				 & x\ge 0.
			\end{array}
		\]
		where $x_+:=\max\{x,0\}$, and $x\circ F(x,\omega)$ is the Hadamard product of $x$ and $F(x,\omega)$.

	\item[(v)] \textbf{Robust Optimisation \cite{ben2002robust, ben2009robust}, which is a deterministic {reformulation} of \eqref{eq:slcp}}. And, 

	\item[(vi)] \textbf{\gls{cm} reformulation \cite{chen2011cvar}}. By using this method, \eqref{eq:slcp} is reformulated to a problem that minimises the CVaR of the norm of the loss function $\theta(x,\omega)$, namely:
		\[
			\min_{x\in\R^n} CVaR_{\alpha}\lf(\| \theta(x,\omega)\|\rg).
		\]

\end{enumerate}

The reformulation in item (vi) uses the CVaR\index{Value at risk, VaR}, a measure of risk widely applied in financial industry. CVaR was built based on \gls{var} \cite{rockafellar2000optimization,meucci2009risk}. Let $ \omega\in\Omega$ be a vector with random outcomes and let $\theta(x,\omega): \R^n\times\Omega \rightarrow \R$ be a mapping, the VaR of $\omega$ for the loss function is defined as:
\begin{equation}\label{VaR}
	VaR_\alpha(\theta(x,\omega)) = min\{\Theta \in\R| \mathcal{P}[\theta(x,\omega)\ge \Theta] \le \alpha\}.
\end{equation}
where $\mathcal{P}[\cdot]\in[0,1]$ is the probability of the event in the square bracket. We call $\theta(x,\omega)$ the loss function. The probability (also called confidence level) $\alpha \in (0,1)$ quantifies the proportion of ``worst cases" (that is, $\theta(x,\omega)\ge VaR_\alpha(\theta(x,\omega))=\Theta$) in the group of all outcomes, and the other outcomes ($\theta(x,\omega)< \Theta$) would happens with probability $1-\alpha$.
Based on the definition of VaR, CVaR\index{Conditional Value at risk, CVaR} is defined as:
\begin{align}
	CVaR_\alpha(\theta(x,\omega)) & = \frac{1}{\alpha}\E\lf[\theta(x,\omega)\mathbf{1}_{[ VaR_{\alpha}\lf(\theta(x,\omega)\rg),+\infty)} \lf(\theta(x,\omega)\rg)\rg] \label{eq:var2cvarexp}\\
		& = \frac{1}{\alpha}\int_{\theta(x,\omega)\ge VaR_{\alpha}(\theta(x,\omega))} \theta(x,\omega) d\mathcal{P}(\omega) \nonumber\\
		& = \frac{1}{\alpha}\int_{0}^{\alpha}VaR_{\gamma}(\theta(x,\omega))d\gamma, \label{eq:var2cvar}
\end{align}
where $\mathbf{1}_{[ VaR_{\alpha}\lf(\theta(x,\omega)\rg),+\infty)}\lf(\theta(x,\omega)\rg)$ is an indicator function with
\[
	\mathbf{1}_{[ VaR_{\alpha}\lf(\theta(x,\omega)\rg),+\infty)}\lf(\theta(x,\omega)\rg) = \left\{
	\begin{array}{l}
		1 \qquad if\; \theta(x,\omega)\in [ VaR_{\alpha}\lf(\theta(x,\omega)\rg),+\infty) \\
		0 \qquad otherwise.
	\end{array}
	\right.
\]

\begin{figure}[!ht]
	\centering
	\begin{minipage}{1\textwidth}
	\includegraphics[width=1\textwidth, height=0.5\textwidth]{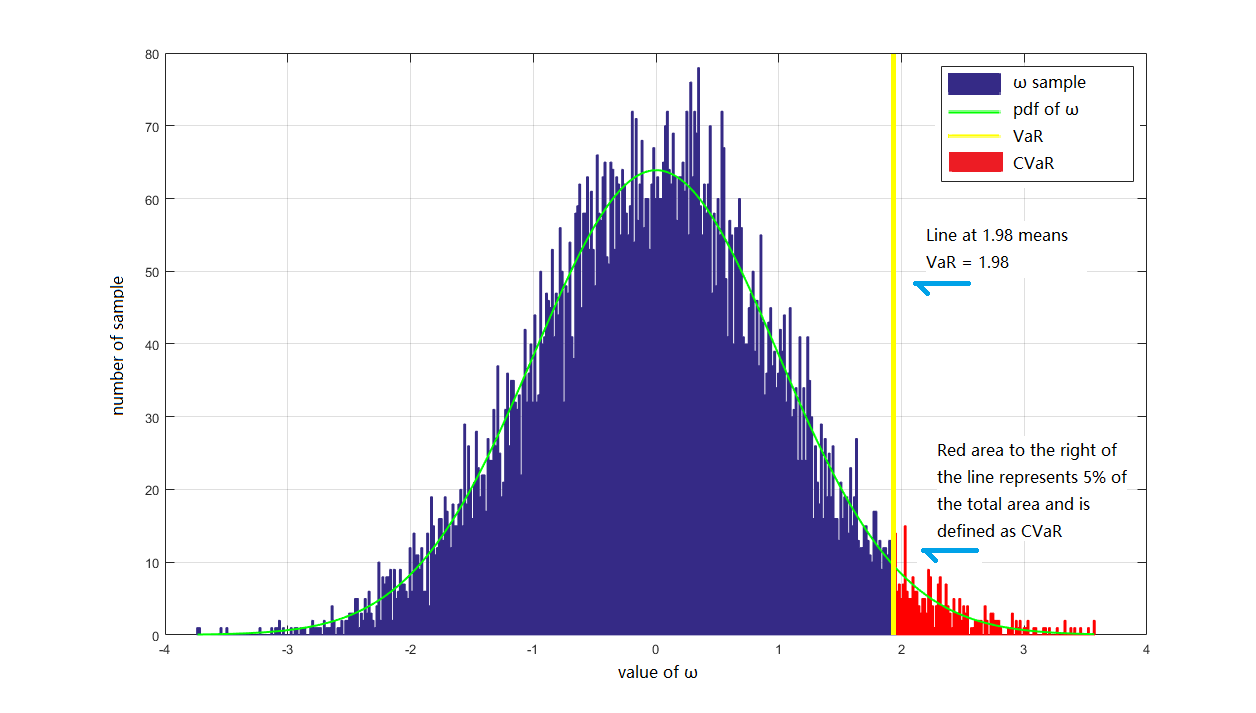}\\
	{\footnotesize Note: For a normal distributed (stochastic) event $\omega\sim N(0,1)$, the distribution of this event shows that only $5\%$ of the outcomes are above 1.98. If the confidence level is set at $95\%$, then the value of VaR equals 1.98(horizontal axis marked by yellow line), and the value of CVaR equals the integral of the area marked in red color. \par}
	\end{minipage}
	\caption{VaR and CVaR for $\theta(x,\omega) = \omega$, where $\omega \sim N(0,1)$}
	\label{fig:cvar}
\end{figure}

$CVaR_\alpha(\theta(x,\omega))$ is the conditional expectation of all outcomes with $\theta(x,\omega)\ge VaR_\alpha (\theta(x,\omega))$. For better understanding the concept of VaR and CVaR, figure \ref{fig:cvar} gives a sample of a loss function $\theta(x,\omega) = \omega$ with one-dimensional normally distributed random value $\omega \sim N(0,1)$. This figure shows that when the confidence level $(1-\alpha)$ is set at $95\%$, the value of VaR equals to the horizontal coordinate of the red vertical line, and the value of CVaR with $95\%$ confidence level equals the red area to the right of the line.

\bigskip \begin{proposition}
	A risk measure $S\lf(\theta(x,\omega)\rg)$ can have more than one of the following properties:
	\begin{enumerate}
		\item Positive homogeneity\index{Homogeneity}: $S\lf(\lambda\theta(x,\omega)\rg) = \lambda S\lf(\theta(x,\omega)\rg)$ for any $\lambda >0$ and $\omega\in \Omega$,
		\item Monotonicity\index{Monotonicity}: if $\theta(x_1,\omega) \ge \theta(x_2,\omega)$ for any $\omega \in \Omega$, we have $S\lf(\theta(x_1,\omega)\rg) \ge S\lf(\theta(x_2,\omega)\rg)$, and
		\item Sub-additivity\index{Sub-additivity}: $S\lf(\theta(x_1,\omega)+ \theta(x_2,\omega)\rg) \le S\lf(\theta(x_1,\omega)\rg) + S\lf(\theta(x_2,\omega)\rg)$  for any $\omega \in \Omega$.
	\end{enumerate}
\end{proposition}

\bigskip \begin{proposition}\label{Var_prop}{\cite{frey2002var}}
	The risk measure VaR is
	\begin{enumerate}
		\item Positive homogeneous, and
		\item Monotonic.
	\end{enumerate}
\end{proposition}

We remark that VaR is not sub-additive. A counter example shown in \cite{danielsson2005subadditivity} will be repeated here.

\bigskip \begin{example}
	Consider the function $\theta(x,\omega) = x + \omega$, where
	\[
		\omega = \epsilon + \eta, \;\; \epsilon \sim N(0,1),\;\; \eta = \lf\{
		\begin{array}{ll}
			0 & with\; probability\; 0.991\\
			10 & with\;probability\; 0.009
		\end{array}
		\rg.
	\]
	In the case when $\alpha = 0.01$, we obtain
	\[
		VaR_{\alpha}\lf(\theta(x,\omega) +\theta(y,\omega)\rg) = x+y +9.8 > VaR_{\alpha}\lf(\theta(x,\omega)\rg) + VaR_{\alpha}\lf(\theta(y,\omega)\rg) = x+3.1 + y + 3.1.
	\]
\end{example}

\bigskip \begin{proposition}{\cite{acerbi2002coherence, frey2002var}}
	The risk measure CVaR is
	\begin{enumerate}
		\item Positive homogeneous,
		\item Monotonic, and
		\item sub-additive.
	\end{enumerate}
\end{proposition}


Consider $\StLCP(F, L, \omega)$ defined by the function $F(x,u,\omega) = T(\omega)\bs x\\u\es + r(\omega)$ and the extended second order cone $L$, problem \eqref{eq:slcp} becomes:
\[
	\StLCP(T(\omega),r(\omega), L)\left\{
	\begin{array}{l}
		Find \; (x,u)\in L,\; such \; that\\ F(x,u,\omega) \in M \; and \; \langle \bs x\\u\es, F(x,u,\omega) \rangle = 0,\;\omega\in \Omega, \quad a.s.
	\end{array}
	\right.
\]
where $T(\omega)=\bs A(\omega) & B(\omega)\\C(\omega) & D(\omega) \es$, with $A(\omega)\in\R^{k\times k}$, $B(\omega)\in\R^{k\times\ell}$, $C(\omega)\in\R^{\ell\times k}$ and $D(\omega)\in\R^{\ell\times\ell}$; $r(\omega) = \bs p(\omega)\\ q(\omega)\es$, with $p(\omega)\in \R^{k}$, $q(\omega)\in \R^{\ell}$, for $\omega\in\Omega$. 

By using item (vi) of Theorem \ref{thm:main}, we reformulate $\StLCP(T(\omega),r(\omega), L\omega)$ to a \gls{stmcp}. The Stochastic mixed complementarity problem $\StMCP(\widetilde{F}_1,\widetilde{F}_2,\R^k_+,\omega)$ defined by $\widetilde{F}_1$, $\widetilde{F}_2$, and $\R^k_+$, is
\[
	\StMCP(\widetilde{F}_1,\widetilde{F}_2, \R^k_+,\omega):\left\{		
	\begin{array}{l}
		Find \; \bs x\\u\\t\es\in\R^k\times\R^{\ell}\times\R, \; such \; that\\ \widetilde{F}_2(x,u,t,\omega)=0, \; and\; (x,\widetilde{F}_1(x,u,t,\omega))\in\C(\R^k_+),\;\omega\in \Omega, \quad a.s.
	\end{array}
	\right.
\]

\bigskip \begin{theorem}\label{cor:Main_thm}

	Suppose $u\ne 0$, $Cx+Du+q\ne 0$. We have
	\[
		z\in \SStLCP(T(\omega), r(\omega), L)\iff \exists t>0,
	\]
	such that
	\[
		\tilde{z}\in \SStMCP(\widetilde{F}_1,\widetilde{F}_2,\R^k_+,\omega),
	\]
	where
	\[
		\widetilde{F}_1(x,u,t,\omega)=A(\omega)(x+te)+B(\omega)u+p(\omega)
	\]
	and
	\small
	\begin{equation}\label{mathcal_SF}
		\widetilde{F}_2(x,u,t,\omega) =
		\begin{pmatrix}
			 \lf(tC(\omega)+ue\tp A(\omega)\rg)(x+te)+ue\tp(B(\omega)u+p(\omega))+t(D(\omega)u+q(\omega)) \\
			 t^2 - \|u\|^2
		\end{pmatrix}.
	\end{equation}
	\normalsize
\end{theorem}

The proof is omitted here, as it inherits the idea of Theorem \ref{thm:main}. Theorem \ref{cor:Main_thm} provides an alternative way to find the solutions to the $\StLCP(T(\omega),r(\omega),L,\omega)$, by converting it to the $\StMCP(\widetilde{F}_1,\widetilde{F}_2, \R^k_+,\omega)$. Such conversion enables us to study $\StLCP(T(\omega),r(\omega),L,\omega)$ through a C-function.

Similar to the process in Chapter 2, \emph{Fischer-Burmeister C-function} will be associated with the problem $\StMCP(\widetilde{F}_1,\widetilde{F}_2,\R^k_+,\omega)$. The FB-based equation formulation of $\StMCP(\widetilde{F}_1,\widetilde{F}_2,\R^k_+,\omega)$ is:
\begin{equation}\label{eq:FB-SMixCP}
	\mathbb{F}^{\StMCP}_{FB}(x,u,t,\omega) =
	\begin{pmatrix}
		\psi_{FB}\lf(x_1,(\widetilde{F}_1)_1(x,u,t,\omega)\rg) \\
		\vdots \\
		\psi_{FB}\lf(x_k,(\widetilde{F}_1)_k(x,u,t,\omega)\rg) \\
		\widetilde{F}_2(x,u,t,\omega)
	\end{pmatrix}.
\end{equation}
where $\psi_{FB}(\cdot ):\R^2\rightarrow \R$ is the scalar FB C-function stated in Chapter \ref{cht:lesoc}. It should be mentioned that the FB C-function is convex, but non-smooth on $\psi_{FB}(0,0)$. According to the definition of FB C-function, a point $\bs x^*\\u^*\\t^*\es$ is a solution to the stochastic mixed complementarity problem $\StMCP(\widetilde{F}_1,\widetilde{F}_2,\R^k_+,\omega)$, if and only if
\begin{equation}\label{eq:stfb}
	\mathbb{F}^{\StMCP}_{FB}(x,u,t,\omega) = 0.
\end{equation}

Based on the results in the previous chapter, of the $\StMCP(\widetilde{F}_1,\widetilde{F}_2,\R^k_+,\omega)$ the associated merit function is:
\begin{equation}\label{eq:MF-SMixCP}
	\theta^{\StMCP}_{FB}(x,u,t,\omega)=\frac{1}{2}\mathbb{F}^{\StMCP}_{FB}(x,u,t,\omega)\tp \mathbb{F}^{\StMCP}_{FB}(x,u,t,\omega).
\end{equation}

Based on \eqref{eq:FB-SMixCP} and \eqref{eq:MF-SMixCP}, the merit function can be written as:
\[
	\theta_{FB}^{\StMCP}(x,u,t,\omega) = \frac{1}{2}\sum\limits_{i=1}^{k}\psi_{FB}^2\lf(x_i,\widetilde{F}_1^i(x,u,t,\omega)\rg) + \frac{1}{2}\sum\limits_{j=1}^{\ell}\widetilde{F}_2^j(x,u,t,\omega).
\]

By the definition of merit function, a point $\bs x^*\\u^*\\t^*\es$ is a solution to the stochastic mixed complementarity problem $\StMCP(\widetilde{F}_1,\widetilde{F}_2,\R^k_+,\omega)$, if
\[
	\theta_{FB}^{\StMCP}(x^*,u^*,t^*,\omega) = 0, \quad \omega\in \Omega \quad a.s.
\]

\bigskip \begin{proposition}
	The associated merit function $\theta_{FB}^{\StMCP}(x^*,u^*,t^*,\omega)$ is continuously differentiable on $\R^k\times\R^{\ell}\times\R$, if $\widetilde{F}_1(x^*,u^*,t^*,\omega)$ and $\widetilde{F}_2(x^*,u^*,t^*,\omega)$ are continuously differentiable on $\R^k$ and $\R^{\ell}\times\R$, respectively.
\end{proposition}
\begin{proof}
	First we prove that $\psi_{FB}^2$ is continuously differentiable. We note that $\psi_{FB}$ is continuously differentiable at every $(a, b)\neq (0,0)$. It is easy to verify that $\psi_{FB}^2$ is continuously differentiable at every $(a, b)\neq (0,0)$. Consider the following to limits at point $(a, b) = (0,0)$:
	\[
		\lim\limits_{\Delta x\rightarrow 0}\frac{\psi_{FB}^2(\Delta x, 0) - \psi_{FB}^2(0, 0)}{\Delta x} = \frac{2\lf(\Delta x^2\rg) - 2\sqrt{\Delta x^2} \cdot \Delta x}{\Delta x} = 2(\Delta x - |\Delta x|) = 0,
	\]
	and
	\[
		\lim\limits_{\Delta y\rightarrow 0}\frac{\psi_{FB}^2(0, \Delta y) - \psi_{FB}^2(0, 0)}{\Delta y} = \frac{2\lf(\Delta y^2\rg) - 2\sqrt{\Delta y^2} \cdot \Delta y}{\Delta y} = 2(\Delta y - |\Delta y|) = 0.
	\]
where $\Delta x$, $\Delta y > 0$. Both partial derivatives of $\psi_{FB}^2$ at $(0,0)$ are continuous, $\psi_{FB}^2$ is continuously differentiable. Hence, $\theta_{FB}^{\StMCP}(x^*,u^*,t^*,\omega)$ is continuously differentiable on $\R^k\times\R^{\ell}\times\R$ if and only if $\widetilde{F}_1(x^*,u^*,t^*,\omega)$ and $\widetilde{F}_2(x^*,u^*,t^*,\omega)$ are continuously differentiable on $\R^k$ and $\R^{\ell}\times\R$, respectively.
\hfill$\square$\par \end{proof}

Next we focus on the convexity of the merit function. The function $\psi_{FB}^2(a,b)$ is not convex on $\R^2$, which implies that the merit function $ \theta_{FB}^{\StMCP}(x,u,t,\omega)$ is not convex on its feasible region.


In case the merit function is not convex, Theorem \ref{thm:Sp} will be helpful for finding the solution to the minimisation problem for a non-convex merit function. By Theorem \ref{thm:Sp}, a FB regular point $\bs x^*\\u^*\\t^*\es$ is a solution to the stochastic mixed complementarity problem $\StMCP(\widetilde{F}_1,\widetilde{F}_2,\R^k_+,\omega)$, if:

\begin{equation}\label{eq:SMIXCP}
	\nabla \theta_{FB}^{\StMCP}(x^*,u^*,t^*,\omega) = 0 \quad \omega\in \Omega, \quad a.s.
\end{equation}

That is
\begin{equation}\label{eq:SMIXCP2}
	\mathcal{A}(\omega)^\top\mathbb{F}_{FB}^{\StMCP}(x^*,u^*,t^*,\omega) =  0 \quad \omega\in \Omega, \quad a.s.,
\end{equation}
where
\[
	\mathcal{A}=
	\begin{pmatrix}
		D_a+D_bJ_x\widetilde{F}_1(x^*,u^*,t^*,\omega) & D_bJ_{\bs u\\t\es}\widetilde{F}_1(x^*,u^*,t^*,\omega)\\
		J_x\widetilde{F}_2(x^*,u^*,t^*,\omega) & J_{\bs u\\t\es}\widetilde{F}_2(x^*,u^*,t^*,\omega)
	\end{pmatrix}
\]
is a nonsingular matrix. Combining equation \eqref{eq:SMIXCP2} with equation \eqref{eq:stfb} implies that equation \eqref{eq:SMIXCP} is a necessary condition for $\bs x^*\\u^*\\t^*\es$ to be a solution to $\StMCP(\widetilde{F}_1,\widetilde{F}_2,\R^k_+,\omega)$.

\begin{figure}[!ht]
  \centering
  \begin{minipage}{1\textwidth}
  \includegraphics[width=1\textwidth, height=1\textwidth]{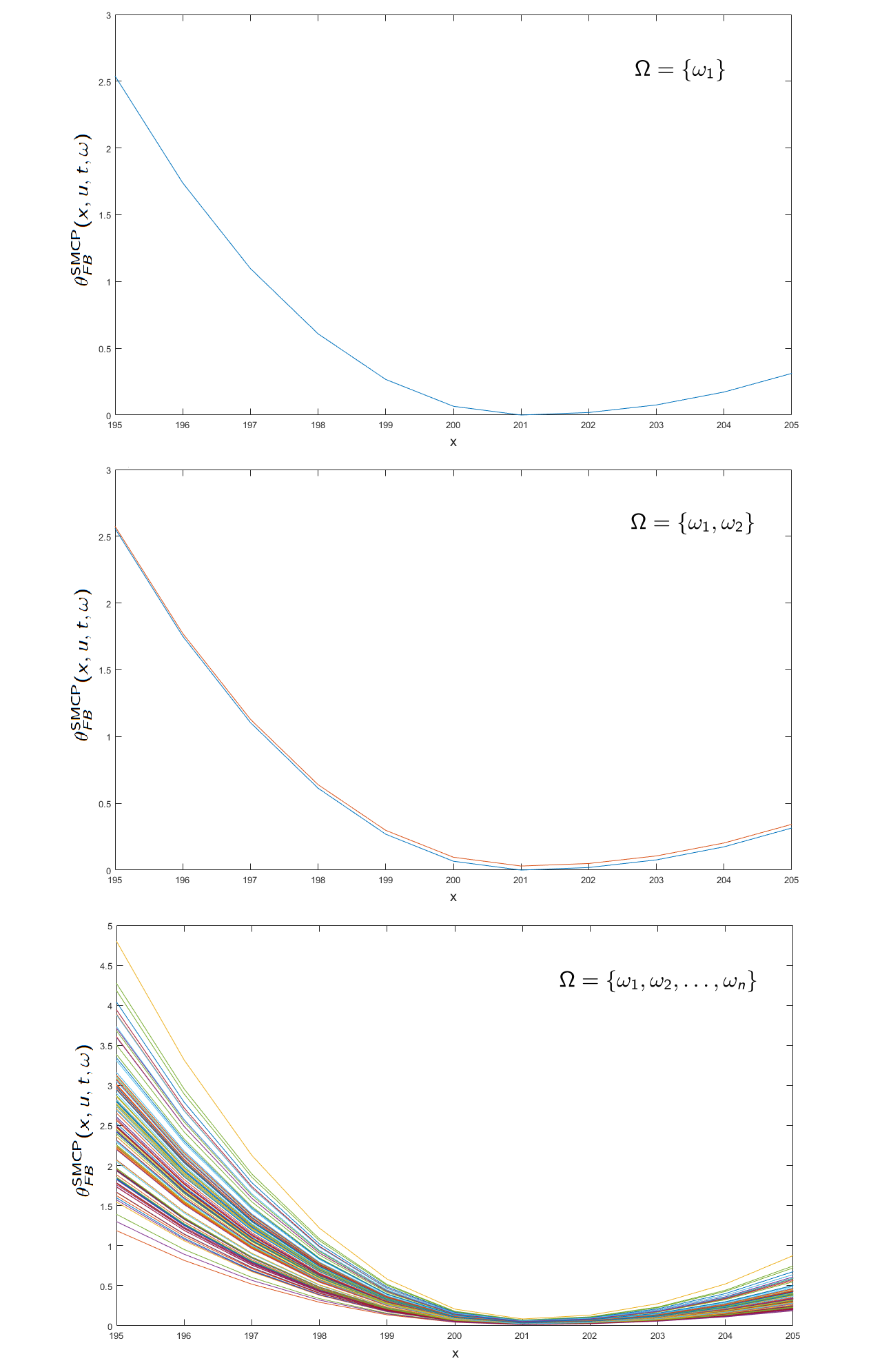}\\
  {\footnotesize Note: For a possible outcome set $\Omega$, when the size of $\Omega$ equals 1, i.e. $|\Omega|=1$ (figure 1), we can easily find a solution (the point when the merit function $\theta^{\StMCP}_{FB}(x,u,t)=0$) to the problem by using the merit function. When $|\Omega|$ increases to 2 (figure 2), the solution for the first case is not longer suitable for both outcomes. As the size of $|\Omega|$ increases (figure 3), it become almost impossible to find a solution to the problem which is suitable for any outcomes. \par}
	\end{minipage}
  \caption{The minimum point of merit function varies $\theta^{\StMCP}_{FB}$ as $|\Omega|$ increases}
  \label{fig:Omega}
\end{figure}

The feasible set of $\StMCP(\widetilde{F}_1,\widetilde{F}_2,\R^k_+,\omega)$ shrinks as $|\Omega|$ (i.e., the size of the possible outcome set $\Omega$) increases. When $|\Omega| = \infty$, we cannot generally find a solution to the problem $\StMCP(\widetilde{F}_1,\widetilde{F}_2,\R^k_+,\omega)$ such that system \eqref{eq:SMIXCP} holds almost surely for any $\omega\in \Omega$, because there will be a large number of equations in system \eqref{eq:SMIXCP}. Figure \ref{fig:Omega} shows the situation when the size of $\Omega$.

As it is introduced above, probability models provide appropriate deterministic reformulations of the stochastic complementarity problems. It can be associated with the stochastic complementarity problems to find persuasive solutions. These persuasive solutions to stochastic complementarity problems would make a proper trade-off between the satisfaction of infinite complementarity constraints and solvability of the problems.

Since $\theta^{\StMCP}_{FB}(x, u, t,\omega)\ge 0$, given a confidence level $(1-\alpha) \in (0,1)$, a point $\bs x^* \\ u^* \\t^*\es$ is a plausible solution to $\StMCP(\widetilde{F}_1,\widetilde{F}_2,\R^k_+,\omega)$ if 

\begin{equation}\label{eq:slacksmcp}
	\bs x^* \\ u^* \\t^*\es\in\arg \min_{x,u,t}\{\Theta|\mathcal{P}\{\theta^{\StMCP}_{FB}(x, u, t,\omega) \le \Theta\} \ge 1 - \alpha\},
\end{equation}

This is a relaxation of problem \eqref{eq:SMIXCP}. A small value of $\alpha$ means that  the satisfaction of the complementarity constraints is preferred to solvability of the problem. A large value of $\alpha$ means that the solvability of the problem is preferred, rather than the satisfaction of the complementarity constraints. Note that the problem \eqref{eq:slacksmcp} can be written as:
\begin{equation}\label{Ind_fun}
	\bs x^* \\ u^* \\t^*\es\in\arg \min_{x,u,t}\{\Theta|\E[\mathbf{1}_{[0,+\infty)}\lf(\theta^{\StMCP}_{FB}(x, u, t,\omega) - \Theta\rg)]\le \alpha\}.
\end{equation}

However, the indicator function $\mathbf{1}_{[0,+\infty)}(\cdot)$ is neither convex nor continuously differentiable at the point $0$. Hence, even though the function $\theta^{\StMCP}_{FB}(\cdot)$ is convex and continuously differentiable, the objective function \eqref{Ind_fun} is non-smooth. If we use the indicator function in the objective function, difficulties occur when applying algorithms which are only viable for smooth objective functions. Addressing this concern, the CVaR method will be considered, which undertakes convex and continuously differentiable objective functions. It harmonises the incompatibility between the satisfaction of infinite number of complementarity constraints and solvability of the problems, as well as inherits convexity (it is not globally convex, but at least it may be convex on some neighbourhoods) and continuous differentiability from the merit function $ \theta^{\StMCP}_{FB} (x,u,t,\omega)$. In the CVaR method, $ (\theta^{\StMCP}_{FB}(x, u, t,\omega) - \Theta)$ will be used as the ``loss function" to measure the ``loss" of complementarity. It should be emphasised that, the higher the value of the ``loss function", the more complementarity constraints of this stochastic complementarity problem are lost. We will transform \eqref{eq:slacksmcp} into CVaR based objective function and then construct the stochastic programming model in the following context.

Rewritting \eqref{Ind_fun} as Value-at-Risk (VaR) to measure of the loss of complementarity:
\[
	\bs x^* \\ u^* \\t^*\es\in \{\bs x\\ u\\ t\es)|VaR_{\alpha}\lf(\theta^{\StMCP}_{FB}(x, u, t,\omega)-\Theta\rg)\le 0\}.
\]

VaR is a measure of complementarity loss defined in \eqref{VaR}. However, the disadvantages of using VaR as the measure of complementarity loss is significant: VaR is not consistent, which means that it is neither convex nor smooth \cite{artzner1999coherent}. On the other hand, CVaR (defined in \eqref{eq:var2cvar}) has superior mathematical properties outperforming VaR, as it inherits continuous differentiability and (local) convexity from the merit function. Moreover, CVaR is a more conservative measure of complementarity loss than VaR.

\bigskip \begin{theorem}\label{thm:cvar_cont}
	If $\theta^{\StMCP}_{FB}(x,u,t,\omega)$ is continuously differentiable on $\R^k\times\R^{\ell}\times\R$, then for any $0<\alpha<1$, the measure of complementarity loss $CVaR_\alpha(\theta^{\StMCP}_{FB}(x,u,t,\omega))$ is continuously differentiable on $\R^k\times\R^{\ell}\times\R$.
\end{theorem}

\begin{proof}
	Immediate from the continuous differentiability of $\theta^{\StMCP}_{FB}(x,u,t,\omega)$ and \eqref{eq:var2cvar}.
\hfill$\square$\par \end{proof}

\bigskip \begin{theorem}\label{thm:cvar_conv}
	If $\theta^{\StMCP}_{FB}(x,u,t,\omega)$ is convex on a neighbourhood ${\cal S}\subset \R^k\times\R^{\ell}\times\R$, then for any $0<\alpha<1$, the measure of complementarity loss $CVaR_\alpha(\theta^{\StMCP}_{FB}(x,u,t,\omega))$ is also convex on ${\cal S}$.
\end{theorem}

\begin{proof}
	Denote $z$, $z'\in {\cal S}\subset \R^k\times\R^{\ell}\times\R$, suppose that $\theta^{\StMCP}_{FB}(z,\omega)$ is convex on ${\cal S}$, we have \[ \theta^{\StMCP}_{FB}(\lambda z+(1-\lambda)z',\omega) \le \lambda\theta^{\StMCP}_{FB}(z,\omega) + (1-\lambda)\theta^{\StMCP}_{FB}(z',\omega), \] where $\lambda\in[0,1]$. Noting that \begin{align*}
		CVaR_\alpha &(\theta(\lambda z+(1-\lambda)z',\omega)) \\
		& = \frac{1}{\alpha}\int_{0}^{\alpha}VaR_{\gamma}\lf(\theta(\lambda z+(1-\lambda)z',\omega)\rg)d\gamma \\
		& \le \frac{1}{\alpha}\int_{0}^{\alpha}VaR_{\gamma}\lf(\lambda\theta^{\StMCP}_{FB}(z,\omega) + (1-\lambda)\theta^{\StMCP}_{FB}(z',\omega)\rg)d\gamma \\
		& = \frac{1}{\alpha}\int_{0}^{\alpha}\lf[VaR_{\gamma}\lf(\lambda\theta^{\StMCP}_{FB}(z,\omega)\rg) + VaR_{\gamma}\lf((1-\lambda)\theta^{\StMCP}_{FB}(z',\omega)\rg)\rg]d\gamma \\
		& = \frac{\lambda}{\alpha}\int_{0}^{\alpha}VaR_{\gamma}\lf(\theta^{\StMCP}_{FB}(z,\omega)\rg)d\gamma + \frac{1-\lambda}{\alpha}\int_{0}^{\alpha} VaR_{\gamma}\lf(\theta^{\StMCP}_{FB}(z',\omega)\rg)d\gamma \\
		& = \lambda CVaR_\alpha (\theta(z,\omega)) + (1-\lambda)CVaR_\alpha (\theta(z',\omega)).
	\end{align*}

	Hence, $ CVaR_\alpha(\theta^{\StMCP}_{FB}(x,u,t,\omega))$ is convex on ${\cal S}$.
\hfill$\square$\par \end{proof}

\bigskip \begin{definition}[Conservativeness]
	Suppose $S_1\lf(\theta(x,\omega)\rg)$, $S_2\lf(\theta(x,\omega)\rg): \R^n \rightarrow {\cal S}$ are two risk measures. Given an outcome $\omega\in \Omega$, risk measure $S_1\lf(\theta(x,\omega)\rg)$ is said to be \emph{more conservative} than risk measure $S_2\lf(\theta(x,\omega)\rg)$ if
	\[
		S_1\lf(\theta(x,\omega)\rg) \ge S_2\lf(\theta(x,\omega)\rg)
	\]
	for any $x\in \R^n$.
\end{definition}

\bigskip \begin{proposition}
	For the measuring the complementarity loss of the merit function $\theta^{\StMCP}_{FB}(x,u,t,\omega)$, the measure $CVaR_\alpha(\cdot)$ is more conservative than the measure $VaR_{\alpha}(\cdot)$.
\end{proposition}

\begin{proof}
	By definition \eqref{eq:var2cvar} we have:

	\begin{align*}
		CVaR_\alpha\lf(\theta^{\StMCP}_{FB}(x,u,t,\omega)\rg)
		& = \alpha^{-1}\int^{\alpha}_{0} VaR_\tau \lf(\theta^{\StMCP}_{FB}(x,u,t,\omega)\rg)d\tau \\
		& = \E[VaR_\tau\lf(\theta^{\StMCP}_{FB}(x,u,t,\omega)\rg)|0\le\tau\le \alpha] \\
		& \ge \min\{VaR_\tau\lf(\theta^{\StMCP}_{FB}(x,u,t,\omega)\rg)|0\le\tau\le \alpha\} \\
		& = VaR_\alpha\lf(\theta^{\StMCP}_{FB}(x,u,t,\omega)\rg)
	\end{align*}

	Hence, we conclude
	\[
		CVaR_\alpha(\theta^{\StMCP}_{FB}(x,u,t,\omega)) \ge VaR_\alpha (\theta^{\StMCP}_{FB}(x,u,t,\omega)).
	\]
\hfill$\square$\par \end{proof}

Reformulate the problem \eqref{Ind_fun} to the following CVaR based minimisation problem: \begin{equation}\label{eq:mincvar}
	\min_{(x,u,t)\in\R^k\times\R^{\ell}\times\R}CVaR_{\alpha}(\theta^{\StMCP}_{FB}(x,u,t,\omega)),
\end{equation}
where
\[
	CVaR_\alpha(\theta^{\StMCP}_{FB}(x,u,t,\omega)) = \alpha^{-1}\int_{0}^{\alpha}VaR_\gamma\lf(\theta^{\StMCP}_{FB}(x,u,t,\omega)\rg)d\gamma,
\]
and
\[
	VaR_\alpha(\theta^{\StMCP}_{FB}) = min\{\Theta | \mathcal{P}[\theta^{\StMCP}_{FB}(x,u,t,\omega)\ge \Theta] \le \alpha\}.
\]

It means that a solution $\bs x^*\\u^*\\t^*\es$ to $\StMCP$ should minimise the ``loss" of complementarity from stochasticity.

Let
\[
	[t]_+ := \max\{0,t\},
\]
\[
	\nu_{(\Theta,\alpha)} (x,u,t, \omega) : = \Theta + \alpha^{-1}[\theta^{\StMCP}_{FB}(x,u,t, \omega) - \Theta]_+,
\]
and define
\[
	\mathcal{N}_{\alpha}(x,u,t, \omega,\Theta): = \E\lf[\nu_{(\Theta,\alpha)} (x,u,t, \omega)\rg] = \Theta + \alpha^{-1}\E[\theta^{\StMCP}_{FB}(x,u,t, \omega) - \Theta]_+.
\]

\bigskip \begin{lemma}\label{prop:cvareqv}
	The problem \eqref{eq:mincvar} is equivalent to the following problem:
	\begin{equation}\label{eq:cvareqv}
		\begin{array}{lcl}
			& \min\limits_{(x,u,t)\in\R^k\times\R^{\ell}\times\R} &\mathcal{N}_{\alpha}(x,u,t, \omega,\Theta^*) \\
		\end{array}
	\end{equation}
	where $\Theta^*$ is the optimal value satisfying:
	\[
		\Theta^* \in \arg\min\limits_{\Theta\in\R} \lf\{ \mathcal{N}_{\alpha}(x,u,t, \omega,\Theta)\rg\}.
	\]
\end{lemma}

\begin{proof}
	Immediate from the alternative definition of CVaR \cite{rockafellar2002conditional}:
	\[
		CVaR_{\alpha} (\theta^{\StMCP}_{FB}(x,u,t,\omega)) := \inf\limits_{\Theta\in\R}\lf\{\Theta + \alpha^{-1}\E[\theta^{\StMCP}_{FB}(x,u,t, \omega) - \Theta]_+\rg\}.
	\]
\hfill$\square$\par \end{proof}

Problem \eqref{eq:cvareqv} simplifies \eqref{eq:mincvar} because it does not contain integration, and inherits the local convexity from the merit function $\theta^{\StMCP}_{FB}(x,u,t, \omega)$. However, since the presence of the operator $[\cdot]_+$, the objective function in problem \eqref{eq:cvareqv} is not smooth at the point 0. {Using mathematical techniques to smooth the objective function can make continuation method applicable on this problem\cite{chen1996class}.} Chen and Harker \cite{chen1997smooth} summarised four palmary smoothing functions. They are provided as follows:
\begin{enumerate}
	\item[(i)] Neural network smoothing function:
	\[
		p(t,\mu) = t + \mu\log(1 + e^{-\frac{t}{\mu}}).
	\]

	\item[(ii)] Interior point smoothing function:
	\[
		p(t,\mu) = \frac{t+\sqrt{t^2 + 4\mu}}{2}.
	\]

	\item[(iii)] Auto-scaling interior point smoothing function:
	\[
		p(t,\mu) = \frac{t+\sqrt{t^2 + 4\mu^2}}{2} + \mu.
	\]

	\item[(iv)] Chen-Harker-Kanzow-Smale (CHKS) smoothing function:
	\[
		p(t,\mu) = \frac{t+\sqrt{t^2 + 4\mu^2}}{2}.
	\]

\end{enumerate}
where $\mu\ge 0$ is the parameter of the approximation function $p$. It should be noted that:
\[
	\lim_{\mu\rightarrow +0}p(t,\mu) = [t]_+.
\]

In this study, we choose \gls{chks} smoothing function and denote:
\[
	[t]_\mu =  \frac{t+\sqrt{t^2 + 4\mu^2}}{2}.
\]

We rewrite problem \eqref{eq:cvareqv} as:
\[
	\begin{array}{lcc}
		& \min\limits_{(x,u,t)\in\R^k\times\R^{\ell}\times\R,\Theta\in\R}
		& \mathcal{N}_{\alpha}(x,u,t, \omega,\Theta) = \Theta + \alpha^{-1}\E[\theta^{\StMCP}_{FB}(x,u,t, \omega) - \Theta]_\mu
	\end{array}
\]

The mathematical expectation is another difficulty that needs to be carefully treated. In many instances, the mathematical expectation $\E[\cdot]$ cannot be calculated with accuracy. A common treatment is using the \gls{saa} method, which is based on the Law of large numbers. SAA method provides a persuasive result of measuring an expectation value \cite{gurkan1999sample, jiang2008stochastic}. If the distribution of the random vector $\omega$ is known, then the Monte-Carlo approach can be used to generate a sample independently and identically distributed (i.i.d.) $\{\omega^1, \dots, \omega^N\}$ with the distribution of $\omega$. Let $\{\omega^1, \dots, \omega^N\}$ be an i.i.d. sample set. The SAA method estimates the mathematical expectation $\E[\theta^{\StMCP}_{FB}(x,u,t, \omega) - \Theta]_\mu$ using averaged value of all observations $[\theta^{\StMCP}_{FB}(x,u,t, \omega^1) - \Theta]_\mu,\, [\theta^{\StMCP}_{FB}(x,u,t, \omega^2) - \Theta]_\mu,\, \dots,\,[\theta^{\StMCP}_{FB}(x,u,t, \omega^N) - \Theta]_\mu$. That is,

\begin{align*}
	\hat{\mathcal{N}_{\alpha}}^N(x,u,t,\Theta) := & \;\frac{1}{N}\sum\limits^N_{i=1}\mathcal{N}_{\alpha}(x,u,t, \omega^i,\Theta)\\
	= & \;\Theta + \alpha^{-1}\frac{1}{N}\sum\limits^N_{i=1}[\theta^{\StMCP}_{FB}(x,u,t, \omega^i) - \Theta]_\mu.
\end{align*}

Then, problem \eqref{eq:cvareqv} becomes
\begin{equation}\label{eq:cvarsaa}
	\begin{array}{lcc}
		& \min\limits_{(x,u,t)\in\R^k\times\R^{\ell}\times\R,\Theta\in\R} & \hat{\mathcal{N}_{\alpha}}(x,u,t, \Theta) = \Theta + \alpha^{-1}\frac{1}{N}\sum\limits^N_{i=1}[\theta^{\StMCP}_{FB}(x,u,t, \omega^i) - \Theta]_\mu \\
		& s.t. & (x,u,t) \in \R^{k}\times\R^{\ell}\times\R, \quad \Theta\in \R.
	\end{array}
\end{equation}

The gradient of $ \hat{\mathcal{N}_{\alpha}}(x,u,t, \Theta) $ is:

\[
	\nabla\hat{\mathcal{N}_{\alpha}}(x,u,t, \Theta) =
	\begin{pmatrix}
		\nabla_{x,u,t}\mathcal{N}_{\alpha}(x,u,t, \Theta) \\
		\nabla_{\Theta}\mathcal{N}_{\alpha}(x,u,t, \Theta)
	\end{pmatrix}
\]
where
\small
\begin{equation}\label{eq:nu1}
	\nabla_{x,u,t}\hat{\mathcal{N}_{\alpha}}(x,u,t, \Theta) = \alpha^{-1}\frac{1}{2N}\sum\limits^N_{j=1} \lf[1+\frac{\theta^{\StMCP}_{FB}(x,u,t, \omega^j) - \Theta }{\sqrt{\lf(\theta^{\StMCP}_{FB}(x,u,t, \omega^j) - \Theta\rg)^2+4\mu}}\rg]\mathcal{A}^\top_j \mathbb{F}_{FB}^{\StMCP}(x,u,t, \omega^j),
\end{equation}
\normalsize

\begin{equation}\label{eq:nu2}
	\mathcal{A}_j =
	\begin{pmatrix}
		D_{a,j}+D_{b,j}\widetilde{A}_{j} & D_{b,j}\widetilde{B}_{j}\\
		\widetilde{C}_{j} & \widetilde{D}_{j}
	\end{pmatrix},
\end{equation}
\footnotesize
\[
	D_{a,j}= diag
	\begin{pmatrix}
		\frac{x_i}{\sqrt{(x_i)^2 + \widetilde{F}_1^i(x,u,t, \omega^j)^2}} - 1
	\end{pmatrix}
	\qquad
	D_{b,j}=diag
	\begin{pmatrix}
		\frac{\widetilde{F}_1^i(x,u,t, \omega^j)}{\sqrt{(x_i)^2 + \widetilde{F}_1^i(x,u,t, \omega^j)^2}} - 1
	\end{pmatrix},
	\qquad
	i\in\{1, \dots , k\},
\]
\normalsize
\[
	\widetilde{A}_{j} = A(\omega^j), \qquad
	\widetilde{B}_{j} =
	\begin{pmatrix}
		B(\omega^j) & A(\omega^j)e
	\end{pmatrix},
	\qquad \widetilde{C}_{j} =
	\begin{pmatrix}
		tC(\omega^j)+ue^\top A(\omega^j) \\
		0
	\end{pmatrix},
\]
\scriptsize
\[
	\widetilde{D}_{j} =
	\begin{pmatrix}
		\lf[A(\omega^j)(x+te)+B(\omega^j)u+p(\omega^j)\rg]\tp eI + ue\tp B(\omega^j) + tD(\omega^j) & C(\omega^j)x+2tC(\omega^j)e+ue^\top A(\omega^j)e+D(\omega^j)u \\ -2u^\top & 2t
	\end{pmatrix}
\]
\normalsize
and
\begin{equation}\label{eq:nu3}
	\nabla_{\Theta}\hat{\mathcal{N}_{\alpha}}(x,u,t, \Theta) = 1 - \alpha^{-1}\frac{1}{N}\sum\limits^N_{j=1} \lf[\frac{1}{2}+\frac{\theta^{\StMCP}_{FB}(x,u,t, \omega^j) - \Theta }{2\sqrt{\lf(\theta^{\StMCP}_{FB}(x,u,t, \omega^j) - \Theta\rg)^2+4\mu}}\rg].
\end{equation}

Since the objective function $ \hat{\mathcal{N}_{\alpha}}(x,u,t, \Theta) $ is continuously differentiable, Problem \eqref{eq:cvarsaa} can be solved by finding some solutions $\bs x^*\\u^*\\t^*\\ \Theta^*\es$ to
\begin{equation}\label{eq:grad_n}
	\nabla\hat{\mathcal{N}_{\alpha}}(x,u,t, \Theta)\ = 0.
\end{equation}

\section{An algorithm}

In the previous section, we have modified the $\StLCP(T(\omega),r(\omega), L,\omega)$ to the problem \eqref{eq:cvarsaa} with a convex and continuously differentiable objective function. Furthermore, the solution to the $\StLCP(T(\omega),r(\omega), L,\omega)$ can be obtained by finding some solution $\bs x^*\\u^*\\t^*\\ \Theta^*\es$ to equation \eqref{eq:grad_n}. In this section, an algorithm will be developed to solve \eqref{eq:grad_n}. Different from the algorithms in Chapter 2, the new algorithm designed in this section involves stochasticity. This algorithm contains Monte-Carlo approach to generate i.i.d. random vector sample sets. We denote $z := \bs x\\u\\t\es \in \R^{n+1}$. Given the tolerance $r >0$, stop criterion is that the following condition is satisfied:

\begin{equation}\label{eq:saatol}
	\max_i\lf\{\lf\|\frac{\partial\mathcal{N}_{\alpha}^{(N_j,\mu_t)}(x,u,t, \omega,\Theta)}{\partial z_i}\rg\|\rg\} \le r, \quad i \in\{ 1,\dots, n+1\}.
\end{equation}
It is shown as follows:


\begin{flushleft}
	\textbf{Algorithm 3} (Line search smoothing SAA):
\end{flushleft}

\textbf{Input}: initial point $z_0 := \bs x_0\\u_0\\t_0\es \in \R^{k}\times\R^{\ell}\times\R $, $\Theta_0\in \R$, confidence level $\alpha$, LM parameter $\nu$, the smoothing parameter $\mu$, maximum iteration number $j_{max}$ for $j$, $k_{max}$ for $k$, the sequence of sample set sizes $ N_1<N_2<\dots<N_{j_{max}}$, parameters of the approximation $\nu$, $\mu$, the tolerance $r>0$, $\varepsilon>0$, and parameters for Wolfe conditions $c_1$, $c_2\in (0,1)$.

\textbf{Step 1}: Set  $j = 1$.

\textbf{Step 2}: Set the sample size $N = N_j$, and generate i.i.d samples $\{\omega^1, \dots, \omega^{N}\}$.

\textbf{Step 3}: If $j>1$, and $\lf\| z^j - z^{j-1}\rg\| < \varepsilon$, \textbf{Stop}.

\textbf{Step 4}: Set  $k = 0$, and $y_0 = z_0$.

\textbf{Step 5}: If either \eqref{eq:saatol} or $k = k_{max}$, then set $j = j + 1$, $z_j = y_k$, and go to \textbf{Step 3}.

\textbf{Step 6}: Denote $\bar{\mathcal{A}_j} = \frac{1}{N}\sum^{N_j}_{i=1}\mathcal{A}_i$, and find a direction $d_k \in \R^{k}\times\R^{\ell}\times\R$ such that
\begin{equation}\label{eq:salgo}
	\bar{\mathcal{A}_j}(y_k)\tp \mathbb{F}^{\StMCP}_{FB}(y_k) + \left[\bar{\mathcal{A}_j}\tp (y_k)\bar{\mathcal{A}_j}(y_k) +\mu \mathbb{I}\right]d_k = 0.
\end{equation}

If the system \eqref{eq:salgo} is not solvable or if the condition
\[
	\nabla\hat{\mathcal{N}_{\alpha}}(y_k, \Theta)\tp d_k \le -r\|d_k\|
\]
is not satisfied, (re)set $d_k = -\nabla\hat{\mathcal{N}_{\alpha}}(y_k, \Theta)$.


\textbf{Step 7}: Find step length $s_k \in R_+$ such that
\[
	\hat{\mathcal{N}_{\alpha}}(y_k + s_k d_k, \Theta) \le \hat{\mathcal{N}_{\alpha}}(y_k + s_k d_k, \Theta) + c_1s_k\nabla\hat{\mathcal{N}_{\alpha}}(y_k, \Theta)\tp d_k,
\]
and
\[
	\hat{\mathcal{N}_{\alpha}}(y_k + s_k d_k, \Theta)\tp d_k \ge c_2 \nabla \hat{\mathcal{N}_{\alpha}}(y_k, \Theta)\tp d_k.
\]

\textbf{Step 8}: Set $y_{k+1} := y_k + s_kd_k $ and $k := k + 1$, go to \textbf{Step 5}.

\textbf{Comment}: This algorithm requires the Monte-Carlo approach to generate i.i.d. random vector samples. For a $j \in\{1,\dots, j_{max}\}$, if the value of $N_j$ is large, the algorithm is anticipated to be more accurate, but it will sacrifice time and computing power. On the other hand, if the value of $N_j$'s is small, the costs of finding result is relatively low, but the accuracy of the solution is sacrificed.



\section{A numerical example} This section illustrates a numerical example for the stochastic ESOCLCP. Denote by $L(3,2)$ an extended second order cone in $\R^3\times\R^2$. Let $x\in\R^3$ and $u\in\R^2$ be two real vectors. Denote
\[
	z=\bs x\\u\es\in\R^3\times\R^2,\quad \hat z=\bs x-\|u\|e\\u\es\in\R^3\times\R^2,\quad and \; \tilde{z}=\bs x-t\\u\\t\es\in\R^3\times\R^2\times\R.
\]

Repeat for convenience a stochastic ESOCLCP defined by the extended second order cone $L(3,2)$ and a stochastic linear function $F(x,u,\omega) = T(\omega)\bs x\\u\es + r(\omega)$ is:
\[
	SLCP(T(\omega),r(\omega),L(3,2))\left\{
	\begin{array}{l}
		Find \; x\in L(3,2),\; such \; that\\
		T(\omega)x+r(\omega) \ge 0, x\tp (T(\omega)x+r(\omega)) = 0,\; \omega\in \Omega, \quad a.s. ,
	\end{array}
	\right.
\]
where
\[
	T=\begin{pmatrix} A & B\\C & D \end{pmatrix} =
	\lf(
	\begin{array}{crcrr}
		41+\omega_1 &  -3 &   -31 &  18 & 19  \\
		 28 &  22 & -33 &  25 &  -29  \\
		-23 & -29 &  11 & -21 &  -43  \\
		- 9 & -31 & -20+2\omega_2 & -12 & 47  \\
		- 8 &  46 &  50 &  -22 & 21
	\end{array}
	\rg), \quad
	r=\lf(
	\begin{array}{r}
		p\\
		q
	\end{array}
	\rg) = \lf(
	\begin{array}{c}
		-26\\
		 4 - \omega_3\\
		 23\\
		 44\\
		-19
	\end{array}
	\rg),
\]
with $A\in\R^{3\times 3}$, $B\in\R^{3\times 2}$, $C\in\R^{2\times 3}$, $D\in\R^{2\times 2}$, $p\in\R^3$, and $q\in\R^2$. 
$\omega = \lf(\omega_1, \omega_2, \omega_3\rg)\tp \in\Omega$ is a stochastic vector with i.i.d. random variables $ \omega_i \sim N(0,1)$ for any $i \in \{1,2, 3\}$. It is easy to verify that square matrices T, A and D are nonsingular for any outcome of $\omega_i$ in  $\R $, $i \in \{ 1,2,3\} $.

By using Theorem \ref{cor:Main_thm}, we reformulate $SLCP(T(\omega),r(\omega),L(3,2))$ to a $\StMCP$ defined by $\widetilde{F}_1$, $\widetilde{F}_2$, and $\R^3_+$:
\[
	\StMCP(\widetilde{F}_1,\widetilde{F}_2, \R^3_+,\omega):\left\{
	\begin{array}{l}
		Find \; \bs x\\u\\t\es\in\R^3\times\R^2\times\R, \; such \; that\\
		\widetilde{F}_2(x,u,t,\omega)=0, \; and\; (x,\widetilde{F}_1(x,u,t,\omega))\in\C(\R^3_+),\;\omega\in \Omega, \quad a.s.
	\end{array}
	\right.
\]
where
	\[
		\widetilde{F}_1(x,u,t,\omega)=A(\omega)(x+te)+B(\omega)u+p(\omega)
	\]
	and
	\small
	\[
		\widetilde{F}_2(x,u,t,\omega) =
		\begin{pmatrix}
			\lf[tC(\omega)+ue\tp A(\omega)\rg](x+te)+ue\tp\lf[B(\omega)u+p(\omega)\rg]+t\lf[D(\omega)u+q(\omega)\rg] \\
			t^2 - \|u\|^2
		\end{pmatrix}.
	\]
	\normalsize
We will convert this to the form of \eqref{eq:cvarsaa} and then \eqref{eq:grad_n}. Given $\alpha = 0.05$, we rewrite problem \eqref{eq:cvarsaa} as:
\[
	\min\limits_{(x,u,t)\in\R^3\times\R^2\times\R,\Theta\in\R} \Theta + 0.05^{-1}\frac{1}{N}\sum\limits^{N}_{i=1}[\theta^{\StMCP}_{FB}(x,u,t, \omega^i) - \Theta]_\mu,
\]
where
\[
	\theta_{FB}^{\StMCP}(x,u,t,\omega) = \frac{1}{2}\sum\limits_{i=1}^{3}\psi_{FB}^2\lf(x_i,\widetilde{F}_1^i(x,u,t,\omega)\rg) + \frac{1}{2}\sum\limits_{j=1}^{2}\widetilde{F}_2^j(x,u,t,\omega).
\]

Since the distribution of the random vector $\omega$ is known, we use the Monte Carlo (MC) method to simulate $j_{max}$ sample sets with number of observation $N_1, N_2, \dots, N_{j_{max}}$. The solutions are shown in the following table:

%

\begin{table}[!ht]

\begin{threeparttable}
\centering
\begin{tabular}{c r c c c}
\toprule
	$j$ & $N_j$ & $\bs x\\u\es\tp$ & $F(x,u,\omega)\tp$ \\[0.4ex]
\midrule
	 1 &	  10 & (1.537, 0.273, 1.060, 0.136, -0.262) & (0.784, 29.054, -0.194, -13.466, 25.803) \\
	 2 &	 100 & (1.542, 0.263, 1.058, 0.127, -0.253) & (1.093, 28.552, -0.214, -12.609, 25.544) \\
	 3 &	1000 & (1.549, 0.257, 1.060, 0.122, -0.252) & (1.277, 28.397, -0.162, -12.418, 25.477) \\
	 4 &   10000 & (1.548, 0.262, 1.060, 0.125, -0.254) & (1.215, 28.605, -0.204, -12.701, 25.578) \\
	 5 &  100000 & (1.546, 0.261, 1.059, 0.125, -0.254) & (1.186, 28.587, -0.176, -12.643, 25.516) \\
	 6 & 1000000 & (1.546, 0.261, 1.059, 0.124, -0.254) & (1.200, 28.566, -0.177, -12.617, 25.514) \\
\bottomrule
\end{tabular}

\centering
\begin{tabular}{c r c c c}
\toprule
	$j$ & $N_j$ & Computation time (sec) & Average loss of complementarity & Threshold $\Theta$ \\[0.4ex]
\midrule
	1 &	   10 & 0.090439 & 0.347 & 0.063  \\
	2 &	  100 & 0.696431 & 0.893 & 0.095  \\
	3 &	 1000 & 5.202383 & 1.179 & 0.090  \\
	4 &	10000 & 39.39705 & 1.060 & 0.087  \\
	5 &   100000 & 553.4596 & 1.054 & 0.088  \\
	6 &  1000000 & 4759.294 & 1.073 & 0.089  \\
\bottomrule
\end{tabular}
\begin{tablenotes}
\footnotesize
\item[] Note: The first table shows the solutions to $SLCP(T(\omega),r(\omega),L(3,2))$ and the value of the function $F(x,u,\omega)\tp$ with respect to different value of $N$. The value of solution does not variate significantly, while the value of the function differs but converges to around 1.200 as the value of $N$ increase. The second table shows the computation time (in second), average loss of complementarity, and the value of threshold. The run time increases significantly along with the value $N$ increases. On the other hand, the average loss of complementarity and the value of threshold remains relative constant no matter what change to the value of $N$. 
\end{tablenotes}
\end{threeparttable}
\caption{The result of the numerical example}
\end{table}

\bigskip

The \gls{aloc} is calculated by:
\[
	ALoC = \frac{1}{N_j}\sum_{i=1}^{N_j}\|(x,u)\tp F(x,u,\omega_i)\|.
\]

As it is shown in the table, the solution converges to $(1.546, 0.261, 1.059, 0.125, -0.254)\tp$ as the value of $N_j$ increases. As the value of $N_j$ increases, the computation time increases as well. However, the Average loss of complementarity and the value of threshold $\Theta$ remains unchanged. It means that it may not be necessary to set a large $N_j$ for the algorithm to get a precise solution. 

\section{Conclusions and comments}

In this chapter, we study the stochastic linear complementarity problem on extended second order cones (stochastic ESOCLCP) which is a stochastic extension of ESOCLCP studied in Chapter \ref{cht:lesoc}. Based on Theorem \ref{thm:main} we derive Theorem \ref{cor:Main_thm}, then we can rewrite an stochastic ESOCLCP to a stochastic mixed complementarity problems (stochastic MixCP) on the nonnegative orthant. Enlightened by the idea from \cite{chen2011cvar}, we introduce the CVaR method to measure the loss of complementarity in the stochastic case. In contrast to the merit function in the deterministic case \eqref{eq:Meritque}, the merit function \eqref{eq:SMIXCP} is not required to equal zero almost surely for any $\omega \in \Omega$. Instead, a CVaR-based minimisation problem \eqref{eq:mincvar} is introduced to obtain a solution which is ``good enough" for the complementarity requirement of the original SMixCP. For solving the CVaR-based minimisation problem derived from the original SMixCP, smoothing function and sample average approximation methods are introduced and finally converted to the form in \eqref{eq:cvarsaa}. Finally, a line search smoothing SAA algorithm is provided for finding the solution to this CVaR-based minimisation problem and it is illustrated by a numerical example.

Stochastic methods on complementarity problems were pioneered by Chen and Fukushima \cite{chen2005expected}. They introduced the idea of minimising the square norm of the merit function to solve a stochastic complementarity problem (SCP). This approach is commonly used in many researches \cite{xu2014cvar,ma2013cvar,zhou2008feasible,chen2011cvar}. However, this approach led to non-convexity and consequently increased the difficulty of solving SCP by algorithms. Our algorithm introduced in this chapter only guarantees a stationary point rather than a solution to the problem. The improvement of the process of finding solutions to a stochastic ESOCLCP will be considered as a good topic of our future research.

\chapter{Application: Portfolio Optimisation Problems} \label{cht:psp}

\section{Introduction to portfolio optimisation problem}

In this chapter, we will apply the results obtained from Chapter \ref{cht:lesoc} on the portfolio optimisation problem. The foundation of the mathematical formulation of portfolio optimisation problem is established by the pioneering paper of Markowitz \cite{markowitz1952portfolio}. His \gls{mv} model is a typical quadratic optimisation programming problem. Also, the \gls{kkt} condition of this problem is a complementarity problem on nonnegative orthant. The mean-variance (MV) scheme of portfolio optimisation sets up a classic framework for the research of portfolio optimisation. Among numerous models developed based on the MV model, the mean-absolute deviation (MAD) model attracted our particular interest. The MAD model is introduced by \cite{konno1991mean} as a route to solve large-scale portfolio optimisation problems. Instead of using covariance matrix, the MAD model uses the absolute-deviation of the rate of return as a measure of the risk. From a mathematical point of view, using absolute-deviation as the measurue of the risk is almost equivalent to the way of using covariance. However, the MAD model significantly reduce the computational cost \cite{konno1999mean,konno2005mean}. The KKT condition of the MAD model is a complementarity problem on second order cone. 

Both of these two models contribute to the development of the research of portfolio optimisation problem. However, multiple papers challenge their major drawbacks. Though the MV model requires only the mean values and the covariance matrix of asset return, it still become very computationally expensive when a large number of assets are considered in the optimisation. In addition, the optimal solution derived from MV model turns out to be highly concentrated in just a few assets, which usually means they are not sufficiently diversified. Parameter sensitivity is another drawback of the MV model. The optimal solution is highly sensitive to its parameter, i.e. the asset returns' mean values and its correlation matrix \cite{laloux1999noise}. The influences of parameter sensitivity on the final result are hard to be rescued because of the inevitable appearance of estimation noises and measurement errors. Hence, optimising a portfolio with the basic MV model will be undiversified and inaccurate, which brings considerable limitations to the actual application of the model. Purely using this optimal allocation on investment decision without amendment often causes a poor out-of-sample portfolio performance. Compare to the MV model, the MAD model has a lower computational cost. However, the modulus in the absolute deviation still lead to some difficulty in computation. The Lagrange function of MAD model is semi-smooth. Unlike the MV model, the MAD model does not have an analytical solution. 

Enlightened by the theorems developed in previous chapters, we introduce a portfolio optimisation model based on the MV model and the MAD model: the \emph{Mean-Euclidean  Norm (MEN) model}. We find the analytical solution to this model. 


\subsection{Review of mean-variance model and mean-absolute deviation model}

The MV model considers a single-period investment. Assume that an investor is in a market with $n$ assets to be considered. Let $\tilde{r} \in \R^n$ denote the random column vector of asset returns in a certain period. We use $r = \E[\tilde{r}]\in \R^n$ to denote the mean return vector of the assets, where $\E[\cdot]$ represents the expected value of the random variable in the square bracket. Suppose that this investor has wealth $X$ at the beginning of the period and he would like to invest all his wealth in these $n$ assets. Let $x\in\R^n$ denote the vector of wealth the investor put in these $n$ assets and $\sum^{n}_{i=1}x_i = X$. At the end of the period, the expectation of final wealth $X'$ of the investor will be:
\[
	\E[X'] = ( e +r)\tp x,
\]
where $ e = (1,1,\dots,1)\tp \in R^n$. $w_i = \frac{x_i}{X}, i \in\{ 1,\dots, n\}$ denote the weight of wealth invested in asset $i$, so $ e \tp w = 1$. The expected rate of return $R_p$ of investor's portfolio $P$ will be:
{
\begin{equation}\label{eq:rtn}
	R_p = ( e + r)\tp  w -  1 = r\tp  w.
\end{equation}}

Since the rate of return vector on assets is random, the investor cannot be sure that how much his wealth will be at the end of the investment period. Hence, he needs to measure the risk of portfolio returns. The variance of the asset rate of return is commonly used as a surrogate for risk. Let $\sigma_{ij} = cov(R_i,R_j), i,j \in\{ 1,\dots,n\}$ denote the estimated covariance between $i$th and $j$th asset returns in a certain period. Hence, the covariance matrix is $\Sigma = (\sigma_{ij})$. The risk $ \sigma_p^2$ of the portfolio $P$ is
\begin{equation}\label{eq:risk}
	\sigma_p^2 = w\tp  \Sigma w
\end{equation}

For each unit of wealth he invested, the investor either hopes to earn at least $\gamma$ profit (i.e., the constraint of minimum rate of return), or can only tolerate a risk not exceeding $\sigma$ (i.e., the constraint of maximum risk). Hence, we obtain the following two equivalent mean variance optimisation formulations according to Markowitz's\index{Mean-variance model, MV model} \cite{markowitz1952portfolio}:
\[
	\begin{array}{rcl}
		& \max\limits_w & r\tp w \\
		&   s.t.				   & w\tp  \Sigma w \le \sigma \\
		&						  &  e \tp w = 1
	\end{array}
\]
\begin{equation}\label{eq:p_2}
	\begin{array}{rcl}
		& \min\limits_w & w\tp  \Sigma w \\
		&   s.t.				   & r\tp w \ge \gamma \\
		&						  &  e \tp w = 1.
	\end{array}
\end{equation}

The problem \eqref{eq:p_2} is a quadratic optimisation problem. The Karush-Kuhn-Tucker (KKT) conditions of it can be written as:\index{Mean-variance model, MV model}
{
\[
	\begin{cases}
		 2\Sigma w - \lambda r - \mu  e = 0\\
		r\tp  w -\gamma \ge 0 ,  \lambda \ge 0,  e \tp  w = 1,\\
		\lambda\tp  (r\tp  w - \gamma) = 0
	\end{cases}
\]}
for some $\lambda$ and $\mu$. Noting that $(\lambda,r\tp  w - \gamma)\in {\cal C}(\R_+)$. So this KKT condition is a mixed complementarity problem on nonnegative orthant. Since the covariance matrix $\Sigma$ is a symmetric and positive semi-definite matrix, if there is a vector $(w^*, \lambda^*, \mu^*)\tp $ satisfies the KKT condition above, $w^*$ will be a solution to problem \eqref{eq:p_2}.

An alternative formulation of \eqref{eq:p_2} is to include both risk and return in the objective function by using the Arrow-Pratt absolute risk -aversion index $c_0 >0$ \cite{kallberg1984mis}. Such inclusion can be deemed as trading risks against return. Problem \eqref{eq:p_2} is reformulated as:
\begin{equation}\label{eq:p2_p}
	\begin{array}{rcl}
		& \min\limits_{y,w} & c_0 y - r\tp  w\\
		& s.t.					& y \ge w\tp  \Sigma w \\
		&						  &  e \tp w = 1.
	\end{array}
\end{equation}

%

The solution to problem \eqref{eq:p2_p} is:
\begin{equation}\label{eq:mv_sol}
	w = \lf(2c_0 \Sigma\rg)^{-1}(r- \frac{e\tp\Sigma ^{-1}r}{e\tp\Sigma^{-1} e}e)+ \frac{\Sigma^{-1} e}{e\tp\Sigma^{-1} e},
\end{equation}
{
\[
	y = w\tp \Sigma w.
\]}

Kallberg and Ziemba \cite{kallberg1984mis} showed that different coefficients $c_0$ can reflects different risk-preference attitudes of an investor. When $c_0 \ge 6$, the objective function reveals strong risk-aversion; when $2\le c_0 < 6$, it shows a moderate risk-aversion; whilst $0<c_0<2$, the function reflects a risk-seeking situation.


One major dispute about the MV model is its computational inefficiency. It requires $n(n+1)/2$ covariance coefficients $\sigma_{ij}$ to be calculated based on the historical data or some results of scenario models. The calculation will be tedious when solving a large-scale portfolio optimisation problem, say, a portfolio with 500 securities included. This is a reason why MV model has not been extensively applied in practises. A good way to alleviate the computation difficulty is to switch the risk measure from covariance matrix to absolute deviation. Konno and Yamazaki \cite{konno1991mean} introduced the mean-absolute deviation (MAD) model to reduce the computational cost of the MV model. According to Konno and Yamazaki, if the return is multivariate normally distributed, the MAD model provides similar results with the MV model. The introduction of risk aversion coefficient $c_0$ makes the models equivalent \cite{rudolf1999linear}. The MAD model demonstrated a stronger computability than the MV model as the computational time for solving a linear programming problem will not be substantially increased comparing with that of a quadratic programming problem.

The MAD model outperforms the MV model in many other aspects. The MV model quantifies the portfolio selection into a form with only two criteria: expected returns measured by means and risks measured by covariance. This simple quantification is also criticized as not consistent with any degrees of stochastic dominance \cite{whitmore1978stochastic,levy1992stochastic}. On the other hand, the MAD model depends on a relation of \gls{ssd} (A portfolio is said to be second-order stochastic dominant of another if this portfolio involves less risk and has at least as high return).

Assume that the rates of returns of assets $\tilde{r} = (\tilde{r}_1,\tilde{r}_2, \dots,\tilde{r}_n)\tp\in\R^{n}$ are distributed over a finite (discrete) sequence of points $\{R_j\}=\{(R_{1j},R_{2j},\dots,R_{nj})\tp\}\in\R^{n},\; j \in\{ 1, 2, \dots, T\}$. That is, there are $T$ different scenarios leads to different outcomes of asset returns. Let $f_j,\; j \in\{ 1, 2, \dots, T\}$ denote the probability distribution of the outcomes of the rates of returns of assets:
\begin{equation}\label{eq:f}
	f_j = Pr\{(\tilde{r}_1,\tilde{r}_2, \dots,\tilde{r}_n)\tp=(R_{1j},R_{2j},\dots,R_{nj})\tp\},\; t=1,2,\dots,T.
\end{equation}

The sequences $\{R_j\}$ and $\{f_j\}$ are acquired through historical data or some techniques of future projection. By definition, it is clear that { $\sum_{j=1}^{T}f_j = 1$} and $f_j\ge 0$ for any $j \in\{ 1, \dots, T\}$. In particular,
\begin{equation}\label{eq:mean_r}
	r = \E[\tilde{r}] =\sum\limits_{j=1}^{T}f_j R_j \in\R^n.
\end{equation}

Denote $U = (U_1, U_2, \dots, U_T)\tp$, where $U_j = (R_j-r),\; j \in\{ 1,\dots, T\}$. The MAD model is the following linear programming problem:\index{Mean-absolute deviation model, MAD model}
\begin{equation}\label{eq:MAD}
	\begin{array}{rcl}
		& \min\limits_{y,w} & c_0 f\tp y - r\tp w\\
		& s.t. & y_j \ge |U_j\tp w|, \quad j \in\{ 1,\dots,T\} \\
		&						  &  e \tp w = 1.
	\end{array}
\end{equation}

%
%
Denote by ${\cal L}$ the second order cone: 
\[
	{\cal L} := \lf\{ (x,y)\in \R\times\R^n:x\geq \|y\| \rg\}.
\]

The KKT condition of Problem \eqref{eq:MAD} can be written as the following complementarity problem on second order cone:

\[
	\mathcal{L} \ni
	\begin{pmatrix}
		y_j \\ U_j\tp w
	\end{pmatrix}
	\perp
	\begin{pmatrix}
		c_0f_j - \theta_j  \\
		\left(U_j^{-1}\right)\tp \left(- r + \mu e + \sum\limits^T_{i=1}\theta_i\frac{U_i U_i\tp }{|U_i\tp w|}w\right)
	\end{pmatrix}
	\in \mathcal{L}, \quad j \in\{ 1,\dots, T\},
\]
\[
	e \tp w - 1 = 0,
\]
where $\mu\in \R$, and $\theta\in \R^T$ are Lagrangian multipliers. Noting that this KKT condition is a mixed complementarity problem on second order cone ${\cal L}$. 

Unlike the MV model, the MAD model does not have an analytical solution \cite{bower2005portfolio}. If we try to solve for $w$, we have:
\begin{equation}\label{eq:smad}
	w = c_0^{-1}B\lf(r - \frac{e\tp B r}{e\tp B e}e\rg) + \frac{B e}{e\tp B e}
\end{equation}
where
\[
	B = \lf( \sum_{j=1}^{T} \frac{f_j}{|U_j\tp w|}U_jU_j\tp\rg)^{-1}
\]

The existence of the modulus in the term $|U_j\tp w|$ implies that \eqref{eq:smad} is not an analytic solution the MAD model. We emphasis that $T \ge n$ is a necessary condition for the matrix $B$ to be nonsingular. It can be easily proved by some basic linear algebra knowledge.

\subsection{Formulation of the mean-Euclidean  norm model}
~\\
If we slightly modify the constraint of problem \eqref{eq:MAD}, we get the \gls{men}:\index{Mean-Euclidean  norm model, MEN model}

\begin{equation}\label{eq:men}
	\begin{array}{rcll}
		& \min\limits_{y,w} & c_0 f\tp y - r\tp w \\
		& s.t. & y_j \ge\|U_j\|\|w\| & j \in\{ 1,\dots, T\}\\
		&						  &  e \tp w = 1,
	\end{array}
\end{equation}
where $ \|\cdot\|$ is the Euclidean  norm and $\|w\| = \sqrt{\langle w,w\rangle}$.
 As a modification of problem \eqref{eq:MAD}, problem \eqref{eq:men} has a different feasible set comparing to problem \eqref{eq:MAD}. The feasible set of problem \eqref{eq:MAD} is:

\begin{equation}\label{eq:msdm1}
	\mathcal{F}_3 = \lf\{(y,w) : y_j \ge|\lf(R_j-r\rg)\tp w|,\; and \; e \tp w = 1 ,\; j \in\{ 1,\dots, T\} \rg\},
\end{equation}
whereas the feasible set of problem \eqref{eq:men} is

\begin{equation}\label{eq:mdmark}
	\mathcal{F}'_3 = \lf\{(y, w) : y_j \ge \|U_j\|\|w\| ,\; and \; e \tp w = 1,\; j \in\{ 1,\dots, T\}\rg\}.
\end{equation}

The following corollary shows the relationship between these two feasible sets.

\bigskip \begin{proposition}\label{cor:wsubset}
	The feasible set $\mathcal{F}'_3 $ \eqref{eq:mdmark} is a subset of $\mathcal{F}_3$  \eqref{eq:msdm1}.
\end{proposition}

\begin{proof}
	For any $j \in\{ 1, 2, \dots, T\}$, we have by Cauchy's inequality
	\[
			|U_j\tp w|  \le \|U_j\|\|w\|.
	\]

	Hence, we have $\mathcal{F}' \subseteq \mathcal{F}$. The equation holds ($\mathcal{F}' = \mathcal{F}$) only if $ U_j$ and $w$ are linearly dependent.
\hfill$\square$\par \end{proof}

Since $\|U_j\| > 0$ for any $j \in\{ 1,\dots, T\}$, $\mathcal{F}_3' $ can be written as:
\[
	\mathcal{F}_3' = \lf\{(y,w) : \frac{y_j}{\|U_j\|} \ge\|w\|, ~ e \tp w = 1 ,~ j \in\{ 1,\dots, T \}\rg\}
\]

{
\textbf{Comment:} The MEN model \eqref{eq:men} and the MAD model \eqref{eq:MAD} have the same objective function and similar feasible sets, but the feasible set of MEN model is the subset of that of MAD model. Unfortunately, the optimal solution to MAD model may be excluded from the feasible set of MEN model. On the other hand, the advantage of MEN model over MAD model is that the former, by applying the Proposition \ref{prop:cs-esoc}, provides possibility of finding analytical solution. In subsequent we will show in Proposition \ref{prop:msckkt} how can we use Proposition \ref{prop:cs-esoc} to work out an analytical solution to the MEN model. }

Recall the definitions of the mutually dual extended second order cone $L(T,n), M(T,n) \in \R^T\times\R^n$ introduced in \eqref{esoc}, \eqref{desoc}:

\[
	L(T, n) = \{\bs y\\w\es \in \R^T\times \R^n : y \geq \|w\| e ,\, y \ge 0 \},
\]

\[
	M(T, n) = \{\bs y\\w\es \in \R^T\times \R^n : e \tp y\geq \|w\|,\, y\ge 0 \},
\]


Denoting $U_{\|\cdot\|} = (\|U_1\|,\|U_2\|,\dots,\|U_T\|)\tp $. We use ``$\circ$" to represent the Hadamard product \cite{horn1990hadamard}. The KKT condition of Problem \eqref{eq:men} is:
\begin{equation}\label{eq:menkkt}
	L \ni
	\begin{pmatrix}
		y\circ U_{\|\cdot\|}^{-1} \\ w
	\end{pmatrix}
	\perp
	\begin{pmatrix}
		c_0U_{\|\cdot\|}\circ f - \theta  \\
		- r + \mu e  + \frac{w}{\|w\|}\sum\limits^T_{j=1}\theta_j
	\end{pmatrix}
	\in M,
\end{equation}

\begin{equation}\label{eq:sum1}
	e \tp w - 1 = 0,
\end{equation}
where, $\mu\in \R$, and $\theta\in \R^T$ are Lagrangian multipliers.

KKT condition of the M2LN model is a nonlinear complementarity problem on ESOC. Applying the item (iv) of Proposition \ref{prop:cs-esoc} on condition \eqref{eq:menkkt} with $x = y\circ U_{\|\cdot\|}^{-1}$, $ u = w$, $ z = c_0U_{\|\cdot\|}\circ f - \theta$, and $ v = - r + \mu e + \frac{w}{\|w\|}\sum\limits^T_{j=1}\theta_j$, we obtain the following proposition:

\bigskip \begin{proposition}\label{prop:msckkt}
	If $- r + \mu e + \frac{w}{\|w\|}\sum\limits^T_{j=1}\theta_j \neq 0$, then there exists a parameter $\lambda\in\R_+$ such that
	\[
		- r + \mu e + \frac{w}{\|w\|}\sum\limits^T_{j=1}\theta_j = -\lambda w,
	\]
	\[
		e \tp \lf(c_0U_{\|\cdot\|}\circ f - \theta\rg) = \lf\|- r + \mu e + \frac{w}{\|w\|}\sum\limits^T_{j=1}\theta_j\rg\|, \] and \[ \R^T_+\ni \lf( y\circ U_{\|\cdot\|}^{-1} - \|w\| e \rg) \perp \lf(c_0U_{\|\cdot\|}\circ f - \theta\rg) \in\R^T_+.
	\]
\end{proposition}
\bigskip \begin{remark}\label{rmk: msckkt}
	Item (i) and item (ii), and item (iii) of Proposition \ref{prop:cs-esoc} are inapplicable in the circumstance of finding the solution to problem \eqref{eq:men}. In Proposition \ref{prop:cs-esoc}, item (i) and item (ii) state that $w = 0$, which contradicts to condition \eqref{eq:sum1} as the later requires $w\neq 0$. Therefore, item (i) and (ii) are inapplicable.

	Item (iii) and item (iv) are applicable in finding the solution to problem \eqref{eq:men}. However, both items have its own limitations.
In Proposition \ref{prop:cs-esoc}, item (iii) assume that
	\begin{equation}\label{eq:lem5eq0}
		- r + \mu e  + \frac{w}{\|w\|}\sum\limits^T_{j=1}\theta_j = 0
	\end{equation}

	We will use the following proposition to show that the conjecture in \eqref{eq:lem5eq0} is not always appropriate.

\end{remark}

\bigskip \begin{proposition}\label{prop:item_iii}
	 With conjecture $ - r + \mu e  + \frac{w}{\|w\|}\sum\limits^T_{j=1}\theta_j = 0$ we cannot always find a solution to problem \eqref{eq:men} for any $\{R_j\}\in \R^n$, $j \in\{ 1, \dots, T\}$.
\end{proposition}

\begin{proof}
	Suppose that \eqref{eq:lem5eq0} holds, it can be rewritten to:
	\begin{equation}\label{eq:lem5eq1}
		w = \lf( r - \mu e \rg)\frac{ \|w\|}{\sum^T_{j=1}\theta_j},
	\end{equation}
	adding \eqref{eq:sum1}, we get
	\begin{equation}\label{eq:lem5eq2}
	   1 = e \tp w = \lf(e \tp r -n \mu\rg) \frac{\|w\|}{\sum^T_{j=1}\theta_j}.
	\end{equation}

	Combine \eqref{eq:lem5eq1} and \eqref{eq:lem5eq2}, we conclude that:

	\begin{equation}\label{eq:lem5eq3}
		w = \frac{ r - \mu e }{ e \tp r - n \mu}.
	\end{equation}

	By KKT condition \eqref{eq:menkkt}, we have:
	\begin{equation}\label{eq:lem5eq13}
	\R^T_+ \ni \lf(y\circ U_{\|\cdot\|}^{-1}\rg) \perp \lf(c_0U_{\|\cdot\|}\circ f - \theta\rg) \in\R^T_+.
	\end{equation}

	On the other hand, by the complementarity in \eqref{eq:lem5eq13}, and $ y\circ U^{-1}_{\|\cdot\|} > 0$ implied by \eqref{eq:menkkt} and \eqref{eq:sum1}, we have
	\[
		0 =  c_0\|U_j\|f_j - \theta_j, \quad j \in\{ 1,\dots, T\},
	\]
	that is
	\begin{equation}\label{eq:lem5eq4}
		\theta_j =  c_0\|U_j\|f_j , \quad j \in\{ 1,\dots, T\}.
	\end{equation}

	Substitute $w$ and $\theta_j$ in \eqref{eq:lem5eq1} by \eqref{eq:lem5eq3} and \eqref{eq:lem5eq4}, respectively, we have
	\begin{equation}\label{eq:lem5eq5}
		\lf( r - \mu e \rg)\lf(1 - \frac{c_0 U_{\|\cdot\|}\tp f}{\|r - \mu e \|}\rg) = 0.
	\end{equation}

	Apparently, $ r \neq \mu e $ because of the nature of asset returns. Let the term in the right bracket of \eqref{eq:lem5eq5} equals zero. Recalling $ U\tp_{\|\cdot\|} = \lf(\|R_1 -r\|, \dots, \|R_T - r\|\rg)$, we get:
	\[
		\|r - \mu e\| = c_0 U\tp_{\|\cdot\|}f = c_0 \sum_{j=1}^{T}\|R_j - r\|f_j.
	\]

	Solving above equation for $\mu$, we have
	\[
		\|r\|^2 - 2\bar{r}\mu + n\mu^2 = \lf( c_0 \sum_{j=1}^{T}\|R_j - r\|f_j\rg)^2,
	\]
	\begin{equation}\label{eq:lem5eq6}
		\mu = \bar r \pm \sqrt{\frac{\bar r - \|r\|^2 + \lf( c_0 \sum_{j=1}^{T}\|R_j - r\|f_j\rg)^2}{n}}.
	\end{equation}

	Recall that $\bar{r} = \frac{1}{n}e\tp r$. From \eqref{eq:lem5eq6} we can observe that if we pick up some $\{R_j\}\in \R^n$, $j \in\{ 1, \dots, T\}$ such that
	\begin{equation}\label{eq:le5eq7}
		\bar r + \lf( c_0 \sum_{j=1}^{T}\|R_j - r\|f_j\rg)^2\le \|r\|^2
	\end{equation}
	then $\mu$ is not a real number.
\hfill$\square$\par \end{proof}

\bigskip \begin{example}\label{exp:item_iii}
	This example gives a numerical case to show that inequality \eqref{eq:le5eq7} does not hold for some $\{R_j\}\in \R^n$, $j \in\{ 1, \dots, T\}$. Given the absolute risk-preference index $c_0 = 4$, number of asset class $n = 3$, and the number of scenarios $ T = 5$:
	\[
		f =
		\begin{pmatrix}
			0.01 \\
			0.14 \\
			0.27 \\
			0.12 \\
			0.46
		\end{pmatrix}, \qquad
			R =
		\begin{pmatrix}
			0.10 & 0.70 &  0.80 &  0.80 & 1.00 \\
			0.30 & 0.80 &  0.60 &  0.40 & 0.70 \\
			0.50 & 0.60 &  0.50 &  0.00 & 0.60
		\end{pmatrix},
	\]
	then we can calculate $r = (0.8710, 0.6470, 0.5000)\tp$ and $\bar r = 0.6727$. Hence,
	\begin{align*}
		  &\; \bar r + \lf( c_0 \sum_{j=1}^{T}\|R_j - r\|f_j\rg)^2 - \|r\|^2 \\
		= &\, 0.6727 + 0.7251 - 1.4273 \\
		= &\, -0.0294 < 0
	\end{align*}
\end{example}

Proposition \ref{prop:item_iii} and Example \ref{exp:item_iii} reveal that the assumptions in Proposition \ref{prop:cs-esoc} item (iii) does not always hold. Using this item of Proposition \ref{prop:cs-esoc} is appropriate only if:
\begin{equation}\label{eq:le5eq8}
		\bar r + \lf( c_0 \sum_{j=1}^{T}\|R_j - r\|f_j\rg)^2 - \|r\|^2 \ge 0
\end{equation}

\begin{figure}[!ht]
	\centering
	\begin{minipage}{1\textwidth}
	\includegraphics[width=1\textwidth, height=0.85\textwidth]{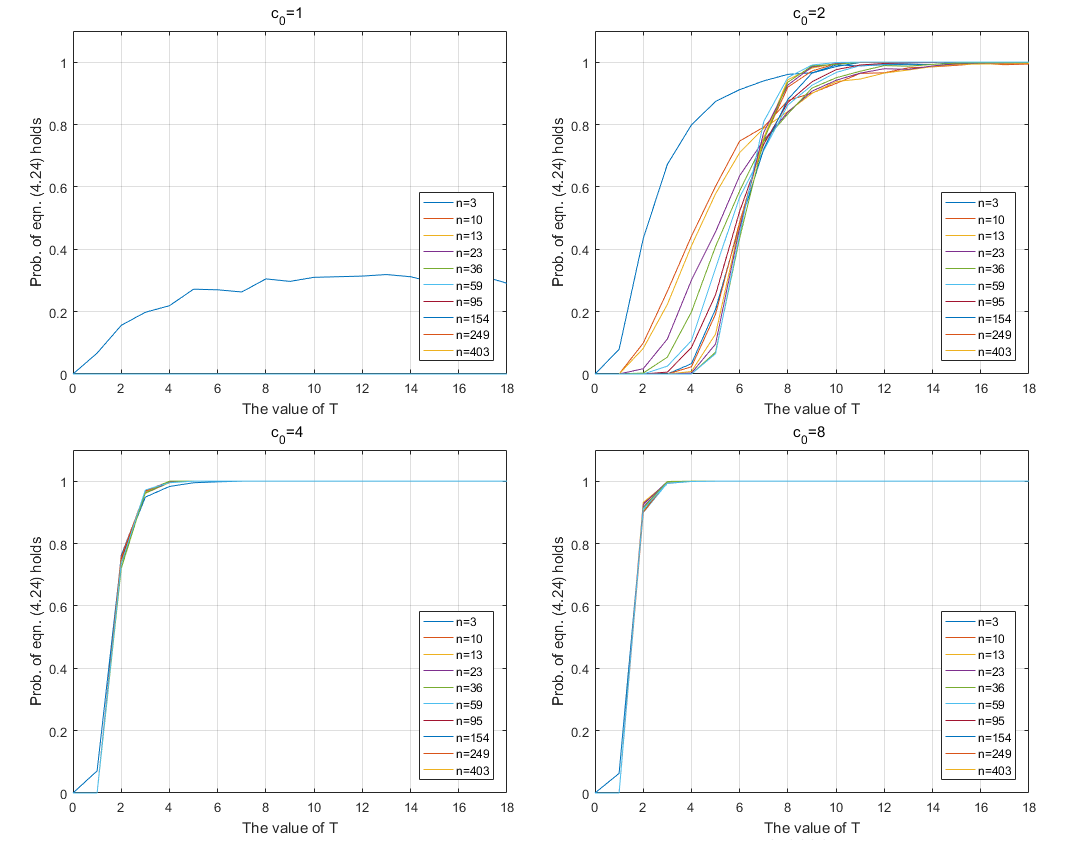}\\
	{\footnotesize Note: This figure shows the correlation between probability and the size of sample set. When $c_0 = 1$, no matter what the values of $n$ and $T$ are, the inequality \eqref{eq:le5eq8} is less likely to hold. Excluding the case when $c_0=1$, the probability is approaching 1 as the size of sample set increases. It can be observed that the parameter $n$ does not have a linear correlation with the probability. The correlation between them seems more likely quadratic. In the case when $c_0 > 0$, as the value of $n$ increases, the probability decreases first and then increases after $n \ge 154$. Curves converge to 1 as $T$ increases. \par}
	\end{minipage}
	\caption{The graph of the probability of inequality \eqref{eq:le5eq8} to hold relates to its parameters}
	\label{fig:item_iii}
\end{figure}

For any $\{R_j\}\in \R^n$, $j \in\{ 1, \dots, T\}$, whether inequality \eqref{eq:le5eq8} holds or not depends on the parameters such as the number of assets ($n$), the size of sample set ($T$), and absolute risk-preference index ($c_0$). The numerical experiment is implemented by generating a series of random return data $R$ and its distribution $f$, then test whether \eqref{eq:le5eq8} holds with the data. Hence, we can calculate the probability of inequality \eqref{eq:le5eq8} to hold. By doing some numerical experiments, we discover that the probability of inequality \eqref{eq:le5eq8} to hold is positively correlated to $n$ and $T$. Also, a large absolute risk-preference index $c_0$ also leads to a higher probability that \eqref{eq:le5eq8} holds. If the value of $n$ and $T$ are large enough, inequality \eqref{eq:le5eq8} holds with a probability almost equals 1. Figure \ref{fig:item_iii} can give an impression to readers about the correlations between the probability of inequality \eqref{eq:le5eq8} to hold and its parameters. In Figure \ref{fig:item_iii}, each curves represent a result with different number of asset (n). The sample number of asset (n) is selected based on a Fibonacci sequence.

Proposition \ref{prop:msckkt} item (iv) is also an appropriate case for problem \eqref{eq:men}. Suppose $- r + \mu + \frac{w}{\|w\|}\sum\limits^T_{j=1}\theta_j \neq 0$, by item (iv) of Proposition \ref{prop:msckkt}, there exists a parameter $\lambda\in\R_+$ such that
\begin{equation}\label{eq:lem6eq1}
	- r + \mu e + \frac{w}{\|w\|}\sum\limits^T_{j=1}\theta_j = -\lambda w,
\end{equation}

\begin{equation}\label{eq:lem6eq2}
	e\tp \lf(c_0U_{\|\cdot\|}\circ f - \theta\rg) = \lf\|- r + \mu e + \frac{w}{\|w\|}\sum\limits^T_{j=1}\theta_j\rg\|,
\end{equation}

and
\begin{equation}\label{eq:lem6eq3}
	\R^T_+\ni \lf( y\circ U_{\|\cdot\|}^{-1} - \|w\| e \rg) \perp \lf(c_0U_{\|\cdot\|}\circ f - \theta\rg) \in\R^T_+.
\end{equation}

From \eqref{eq:lem6eq1} and \eqref{eq:lem6eq2} , we get
\begin{align}
	w & = \frac{\|w\|\lf( r - \mu e \rg)}{\sum\limits_{j=1}^T\theta_j + \lambda\|w\|},\label{eq:lem6eq4}\\
	\lambda \|w\| & = e \tp \lf(c_0U_{\|\cdot\|}\circ f - \theta\rg). \label{eq:lem6eq5}
\end{align}

Combine \eqref{eq:lem6eq4}, \eqref{eq:lem6eq5} and \eqref{eq:sum1}, we have
\[
	\mu = \bar r - \frac{c_0\sum_{j=1}^{T}\|R_j - r\|f_j}{n\|w\|}
\]

Substitute $\mu$ in \eqref{eq:lem6eq4}, we obtain
\begin{equation}\label{eq:smenraw}
	w = \frac{\lf( r - \bar{r} e \rg)}{c_0 \sum_{j=1}^{T}\|R_j - r\|f_j}\|w\| + \frac{ e }{n}.
\end{equation}

The modulus $|U_j\tp w|$ in \eqref{eq:smad} is an obstacle against finding analytical solution to problem \eqref{eq:MAD}, but the norm $\|w\|$ in \eqref{eq:smenraw} is removable therefore leading to the analytical solution to problem \eqref{eq:men}. Since
\begin{equation}\label{eq:wnorm}
	\|w\|^2 = \langle w,w\rangle,
\end{equation}
noting that $\langle r- \bar{r}e, e\rangle = 0$, we substitute $w$ in the right-hand side of the equation \eqref{eq:wnorm} by \eqref{eq:smenraw}, then we have:

\begin{align*}
	\|w\|^2  & = \lf\langle \frac{\lf( r - \bar{r} e \rg)}{c_0 \sum_{j=1}^{T}\|R_j - r\|f_j}\|w\| + \frac{ e }{n}, \frac{\lf( r - \bar{r} e \rg)}{c_0 \sum_{j=1}^{T}\|R_j - r\|f_j}\|w\| + \frac{ e }{n} \rg\rangle \\
	& = \frac{\|w\|^2}{\lf(c_0 \sum_{j=1}^{T}\|R_j - r\|f_j\rg)^2} \lf\langle \lf( r - \bar{r} e \rg) - \frac{ e }{n}, \lf( r - \bar{r} e \rg) - \frac{ e }{n}\rg\rangle \\
	& = \frac{\|w\|^2}{\lf(c_0 \sum_{j=1}^{T}\|R_j - r\|f_j\rg)^2} \|r-\bar{r}e\|^2 + 1,
\end{align*}
then make some transformations, we have the following equation
\begin{equation}\label{eq:smentrn}
	\lf(1 - \frac{\|r-\bar{r}e\|^2}{\lf(c_0 \sum_{j=1}^{T}\|R_j - r\|f_j\rg)^2}\rg)\|w\|^2 = \frac{1}{n}.
\end{equation}

Since $\|w\| >0$, $\|w\|$ in \eqref{eq:smentrn} is a real number only if:
\begin{equation}\label{eq:smenieq}
	1 - \frac{\|r-\bar{r}e\|^2}{\lf(c_0 \sum_{j=1}^{T}\|R_j - r\|f_j\rg)^2} > 0.
\end{equation}
%
%
%

\begin{figure}[!ht]
	\centering
	\begin{minipage}{1\textwidth}
	\includegraphics[width=1\textwidth, height=0.85\textwidth]{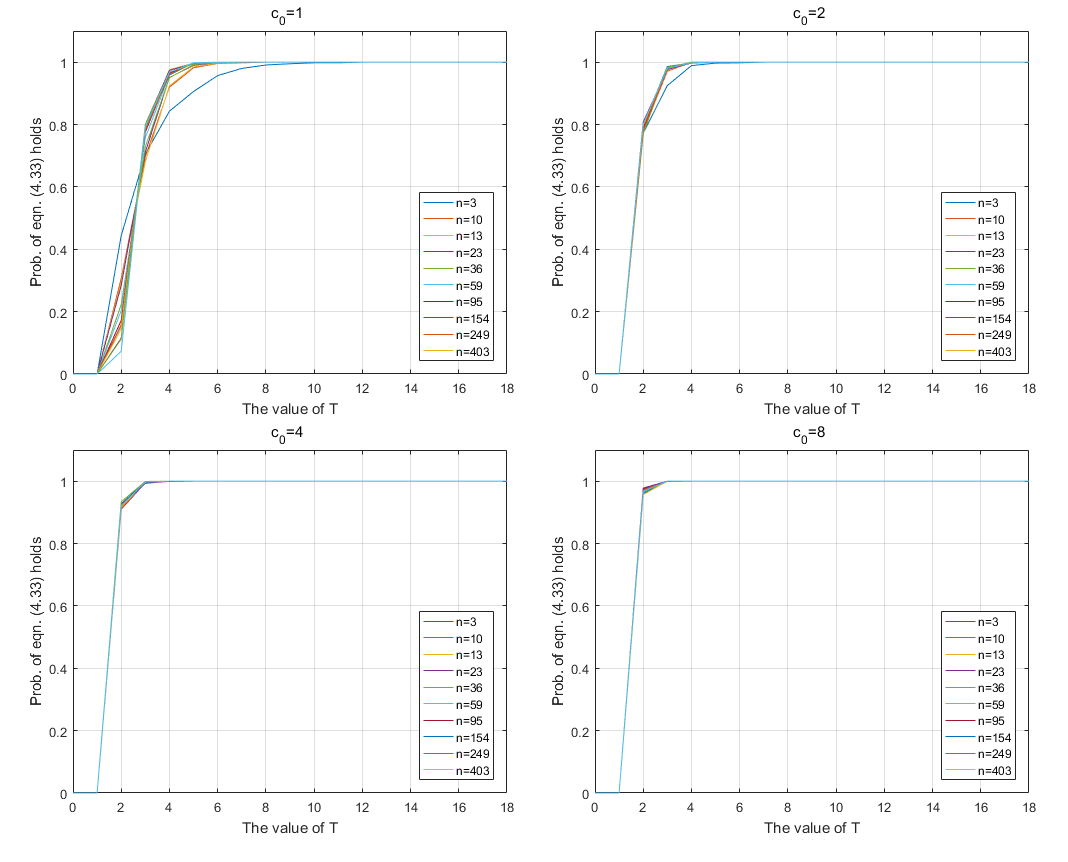}\\
	{\footnotesize Note: Compare to Proposition \ref{prop:cs-esoc} item (iii), item (iv) is more applicable. Unlike the case in Figure \ref{fig:item_iii}, when $c_0 = 1$, the probability of inequality \eqref{eq:smenieq} to hold will be very close to 1 if $T \ge 12$. For $c_0 \ge 1$, curves converge to 1 as $T$ increases. It also shows a quadratic correlations between probability and parameter $n$. \par}
	\end{minipage}
	\caption{The graph of the probability of inequality \eqref{eq:smenieq} to hold relates to its parameters}
	\label{fig:item_iv}
\end{figure}

Figure \ref{fig:item_iv} shows the correlations between inequality \eqref{eq:smenieq} and its parameters. Similar to inequality \eqref{eq:le5eq8}, inequality \eqref{eq:smenieq} is not always hold. However, inequality \eqref{eq:smenieq} is somehow more plausible than inequality \eqref{eq:le5eq8}. When the size of sample set ($T$) is large enough (greater than 12 in the $c_0=1$ case), the inequality \eqref{eq:smenieq} has very high probability to hold. Unlike inequality \eqref{eq:le5eq8}, the probability is less related to the risk-preference index $c_0$. It means that it is suitable for more scenarios.

\bigskip

Suppose that inequality \eqref{eq:smenieq} holds, we can solve for $w$ by \eqref{eq:smentrn} and \eqref{eq:smenraw}:
\begin{equation}\label{eq:smenfin}
	w = \frac{\lf( r - \bar{r} e \rg)}{\sqrt{n\lf(\lf(c_0\sum\limits_{j=1}^{T}\|R_j - r \|f_j\rg)^2 - \|r-\bar{r} e \|^2\rg)}} + \frac{e}{n}.
\end{equation}

The result \eqref{eq:smenfin} indicates that the weight of asset $i$ relates to its excess return $\lf( r_i - \bar{r} e \rg)$, total risk of the market $U_{\|\cdot\|}\tp  f$, and the absolute risk-preference index $c_0$. MEN conservatively consider the total market risk instead of the risk of single asset in the process of determining the weight of each asset.

From \eqref{eq:lem6eq3}, we get
\begin{equation}
	\lf( y_j \|U_j\|^{-1} - \|w\|\rg)\lf(c_0\|U_j\| f_j - \theta_j\rg) = 0, \quad j \in\{ 1, \dots, T\},
\end{equation}
and
\[
	y_j = \|U_j\|\|w\|, \quad j \in\{ 1, \dots, T\},
\]
always hold if the objective function is minimised.

\section{Conclusions and comments}

In this chapter, we introduced a modified version of portfolio selection model based on the mean-variance model (MV) and the mean-absolute deviation model (MAD): the mean-Euclidean  norm (MEN) model. The MV model has its analytical solution. However, this analytical solution requires the estimation of the covariance matrix, and the estimation of the covariance matrix of the MV model is computationally expensive. In addition, the result of such estimation is always negatively influenced by the estimation noises. Based on the MV model, Konno and Yamazaki \cite{konno1991mean} introduced the MAD model. The advantage of MAD model is that does not require to calculate covariance matrix of asset returns, therefore its computational costs is less than that of the MV model. However, given the modulus term in the constraint of MAD model, it is almost impossible to work out an analytical solution to the MAD model. Hence, when the number of assets is large, it also becomes computationally expensive in the process of finding the optimal weight of assets. Hence, the results from both the MV model and MAD model are criticised by many scholars \cite{michaud2008efficient, laloux1999noise, pafka2003noisy, laloux2000random}.

The MEN model considers a similar way of finding the optimal portfolio to both the mean-variance (MV) model and the MAD model: maximising the risk adjusted return. The objective of introducing this model is to find a new method based on the historical return data of assets. This model is designed for calculating the optimal weight of a portfolio with low computational cost. It is also designed to help in reducing the negative impacts of the inaccuracy of optimal solutions caused by estimation noises. The idea of introducing the MEN model is based on Corollary \ref{cor:wsubset}, which conservatively tighten the feasible set of the MAD model from \eqref{eq:msdm1} to \eqref{eq:mdmark}. Given the fact that the KKT condition of the MEN model is a nonlinear complementarity problem on extended second order cone, we innovatively obtained the analytical solution to the MEN model based on Proposition \ref{prop:msckkt} in Chapter \ref{cht:lesoc}. 

The major advantage of the MEN model is has analytical solutions. 
However, due to the insufficient time and effort, this study fails to provide an empirical evaluations for comparing these three portfolio optimisation approaches. The empirical evaluations usually contain constraints such as short-sales constraints, lower-bound/upper-bound-on-individual-asset constraints, sector constraints, etc. These constraints are commonly occurred in practices. Without considering these constraints, the empirical evaluation results always looks absurd. Hence, our future work is to explore the MEN models with linear constraints and evaluate this new model empirically. 

Our current study can be extended and improved from practical aspects. The empirical evaluations for the MEN model is the suggested direction. Also, the complementarity problem on extended second order cone has potential on the application of solving the asset allocation problem for the defined-contribution pension plan \cite{cairns2000optimal}.

\chapter{Spherically quasi-convex quadratic functions}\label{cht:sphere}

\section{Motivation of this study}

In this chapter, we study the spherical quasi-convexity of quadratic functions on spherically convex sets, which is related to the problem of finding
their minimiser. The spherically convex set is a natural extension of the concepts and techniques of convexity of mathematical programming problem. The original intention of this study is to explicitise certain fixed point theorems, surjectivity theorems, and existence theorems for complementarity problems and variational inequalities. Parts of the results of this chapter are published in \cite{FerreiraNemethXiao2018}, but we provide more detailed definitions and proofs in this thesis.

Recall the definition of the nonlinear complementarity problem:
\index{Complementarity problem!Nonlinear complementarity problem}
\bigskip \begin{definition}[Nonlinear complementarity problem]
	Let $F:\R^n\to\R^n$ be a mapping. Let ${\cal K}\subseteq\R^n$ be a nonempty closed convex cone and ${\cal K}^*$ its dual.
	Defined by ${\cal K}$ and $F$ the \emph{\gls{ncp}} \cite{FacchineiPang2003} is:
	\begin{equation}\label{eq:nlcp}
		NCP(F,{\cal K})\left\{
		\begin{array}{l}
			Find \; x\in \R^n, \; such \; that\\
			\left(x, F(x)\right)\in\C({\cal K}).
		\end{array}
		\right.
	\end{equation}
	The solution set of $\NCP(F,{\cal K})$ is denoted by $\SNCP(F,{\cal K})$:
	\[
		\SNCP(F,{\cal K}) =  \{x\in \R^n: (x, F(x)) \in \C({\cal K})\}.
	\]
\end{definition}

To explain the motivation of this study in detail, we start with presenting the following three definitions:
\index{Inversion}
\begin{definition}[Inversion]\cite[Definition~5]{Nemeth2006}
	The operator
	\[
		i: \R^n\setminus\{0\} \rightarrow \R^n\setminus \{0\};~ i(x):= \frac{x}{\|x\|^2}
	\]
	is called {\it inversion} (of pole 0).
\end{definition}
It is easy to see that $i$ is a one-to-one mapping, and $i^{-1} = i$.
\begin{definition}[Inversion of a mapping]\cite[Definition~6]{Nemeth2006}
	The inversion (of pole 0) of the mapping $F: \R^n\rightarrow \R^n$, is the mapping $\mathcal{I}(F): {\cal K}\rightarrow \R^n$ defined by:
	\[
		\mathcal{I}(F)(x):=\left\{
		\begin{array}{l}
			\|x\|^2(F\circ i)(x) \;\; if~x\neq 0,\\
			0\;\; if~ x= 0.
		\end{array}
		\right.
	\]
\end{definition}
\index{Inversion of a mapping}

\begin{definition}[Lower scalar derivative]\cite[Definition~1.6]{IsacNemeth2008}
	Consider the mapping $F: {\cal K} \rightarrow \R^n$. The limit
	\[
		\underline{F}^{\#}(x_0, {\cal K}):=\liminf_{x\rightarrow x_0, x-x_0\in {\cal K}}\frac{\langle F(x)-F(x_0),x-x_0\rangle}{\|x-x_0\|^2}
	\]
	is called the lower scalar derivative of $F$ at $x_0$.
\end{definition}
\index{Lower scalar derivative}

The lower scalar derivative is highly related to the minimising questions on spherically convex sets. We can observe the clues from \cite[Theorem~18]{Nemeth2006}. This theorem gives a more explicit expression of the lower scalar derivative if the mapping $F$ is Fr{\'e}chet differentiable (see definition \ref{def:Frechetdif}):
\begin{theorem}\cite[Theorem~18]{Nemeth2006}\label{thm:lower_s_d}
	Let ${\cal K} \subseteq \R^n$ be a closed convex cone with non-empty interior and $x$ an interior point of ${\cal K}$. If $F:{\cal} \rightarrow \R^n$ is Fr{\'e}chet differentiable in $x$, with the differential $JF(x)$, then
	\[
		\underline{F}^{\#}(x, {\cal K}) = \min_{\|u\|=1,u\in{\cal K}}\langle JF(x)(u),u\rangle.
	\]
\end{theorem}

Theorem \ref{thm:lower_s_d} leads to the study about minimising questions on the intersection between a cone and the sphere. The Corollary \cite[Corollary~8.1]{IsacNemeth2006} is the trigger of this study. We state this corollary here for convenience.
\begin{corollary}\cite[Corollary~8.1]{IsacNemeth2006}
	Let $\mathcal{K} \subseteq \R^n$ be a closed convex cone, and $F:{\cal K} \rightarrow \R^n$ be a continuous mapping such that its inversion (of pole 0) is differentiable at 0. Denote by $J {\cal I}(F)(0)$ the Jacobian matrix of the inversion of mapping $F$ at 0. Then if
	\[
		\underline{{\cal I}(F)}^{\#}(0) := \min_{\|x\|= 1, u\in {\cal K}}\langle J {\cal I}(F)(0) u, u\rangle > 0.
	\]
	then the nonlinear complementarity problem $\NCP$ has a solution.
\end{corollary}

By \cite[Theorem~18]{Nemeth2006} and \cite[Corollary~8.1]{IsacNemeth2006}, the question of the existence of the solution to a nonlinear complementarity problem can be converted to a problem of minimising a quadratic function on the intersection between a cone and the sphere. These sets are exactly the spherically convex sets (see \cite{FerreiraIusemNemeth2014}), which leads to the problem of minimising quadratic functions on spherically convex sets.

Apart from  the above, the motivation of this study is much wider. For instance, consider the quadratic constrained optimisation problem on the sphere
\begin{equation}\label{eq:gp}
	\min \{\langle Qx, x\rangle~:~x\in {\cal S}\cap {\cal K}\},
\end{equation}
where ${\cal K}$ is a cone on the sphere, and ${\cal S}\subseteq\SP^{n-1}:=\left\{ x\in \R^{n}~:~\|x\|=1\right\}$ is a sphere. Suppose the matrix $Q$ is a symmetric matrix, then the problem \eqref{eq:gp} is a minimum eigenvalue problem in ${\cal S}$. In particular, this problem includes the problem of finding the spectral norm of the matrix \(-Q\) when ${\cal S}=\SP^{n-1}$ (see, e.g.,~\cite{Smith1994}). We focus on the cases when ${\cal S}$ is an intersection of a subdual convex set with the sphere. Also, it is important to highlight that the special case when ${\cal S}$ is the intersection of the nonnegative orthant with the shpere is of particular interest because the nonnegativity of the minimum value is equivalent to the copositivity of the matrix $Q$ \cite[Proposition~1.3]{UrrutySeeger2010} and to the nonnegativity of all Pareto eigenvalues of $Q$ \cite[Theorem~4.3]{UrrutySeeger2010}. As far as we are aware there are no methods for finding the Pareto spectra by using the intrinsic geometrical properties of the sphere, hence our study is expected to open new perspectives for detecting the copositivity of a symmetric matrix. Another important special case is when ${\cal S}$ is the intersection of the Lorentz cone with the sphere. We pay attention to this case because the minimum eigenvalue of $Q$ in ${\cal S}$ is nonnegative if and only if the matrix $Q$ is Lorentz copositive, see \cite{GajardoSeeger2013,LoewySchneider1975}. In general, changing the Lorentz cone by an arbitrary closed convex cone ${\cal K}$ would lead to a more general concept of ${\cal K}$ copositivity, thus our study is anticipated to initialise new perspectives of investigating the general copositivity of a symmetric matrix. 
More problems that deals with ``spherical" constraint can be found in \cite{Malick2007}.

Optimisation problems posed on the sphere have a specific underlying algebraic structure that could be exploited to greatly reduce
the cost of obtaining the solutions; see \cite{Hager2001, HagerPark2005, Smith1994,So2011,Zhang2003,Zhang2012}. It is
worth to point out that when a quadratic function is spherically quasi-convex, then a spherical strict local minimiser is equal to a spherical strict
global minimiser. Therefore, it is natural to consider the problem of determining the spherically quasi-convex quadratic functions on spherically convex sets. The goal of the study is to
present necessary conditions and sufficient conditions for quadratic functions which are spherically quasi-convex on spherical convex sets. As a particular case, we exhibit several such results for both the spherical positive orthant and, more general, the spherical subdual convex set. 

Apart from the questions about the existence of a solution to the nonlinear complementarity problem, this study also related to many other questions. For example, minimising a quadratic function defined on spherical nonnegative orthant is equivalent to finding the minimum Pareto eigenvalues of the quadratic function. Hence, if the minimum value of this quadratic function is nonnegative, then the minimum Pareto eigenvalues of the quadratic function will also be nonnegative. In this chapter, we present several conditions that characterise the spherical quasi-convexity of quadratic functions.
The study can be considered as a first spherical analogue for the study of quasi-convexity of quadratic functions. Without the aim of completeness, we list
here some of the main papers about the quasi-convexity of quadratic functions: \cite{Martos1969,Ferland1972,Schaible1981,Komlosi1984,Karamardian1993}.

\section{Spherically quasi-convex quadratic functions  on spherically convex sets} \label{sec:qcqfcs}
In this section our aim is  to  present  some conditions   characterising   quadratic spherically quasi-convex functions on a general  spherically convex set. We assume for convenience that {\it from now on the cone ${\cal K}\subseteq \R^{n}$ is a proper subdual cone}. Define
\begin{equation}\label{eq:ccb}
	{\cal S}=\SP^{n-1}\cap\inte({\cal K}), \quad \bar{\cal S} = \SP^{n-1}\cap {\cal K},
\end{equation}
and assume that ${\cal S}$ is an open spherically convex set.

\bigskip \begin{definition}[Quadratic function]
	The associated quadratic function  $q_A: {\cal S}\to\R$ defined by the symmetric matrix $A=A^T\in\R^{n\times n}$ is
	\begin{equation} \label{eq:QuadFunc}
		q_A(x):=\langle Ax,x\rangle.
	\end{equation}
\end{definition}\index{Quadratic function}

We remark that $q_A$ can be extended to $\bar{\cal S}$. For the simplicity of notations we will denote the extended values by $q_A(x)$ too, but the spherical quasi-convexity of $q_A$ will always be understood as a function defined on
${\cal S}$. To proceed we need the following definition:

\bigskip \begin{definition}[Rayleigh quotient function]
	The Rayleigh quotient function $\varphi_A : \inte({\cal K}) \to \R$ restricted on $\inte({\cal K}) $ defined by matrix $A$ is
	\begin{equation} \label{eq:RayleighFunction}
		\varphi_A(x):=\frac{\langle Ax, x \rangle}{\|x\|^2}.
	\end{equation}
\end{definition}\index{Rayleigh quotient function}

In the following proposition we present some  equivalent characterisations of the convexity of the associated quadratic function on spherically convex sets $q_A$ defined in \eqref{eq:QuadFunc}. 
\bigskip \begin{proposition}\label{pr:spher-quasiconv}
Let $q_A$ and $\varphi_A$ be the functions defined in \eqref{eq:QuadFunc} and \eqref{eq:RayleighFunction}, respectively. The following statements are equivalent:
	\begin{enumerate}
		\item[(a)] The quadratic function $q_A$ is spherically quasi-convex;
		\item [(b)] $\langle Ax,y\rangle\leq\langle x,y\rangle\max \left\{q_A(x), ~q_A(y)\right\}$ for any $ x,y\in\SP^{n-1}\cap {\cal K}$;
		\item [(c)] $ \displaystyle \frac{\langle Ax,y\rangle}{\langle x,y\rangle}\leq\max\left\{\varphi_A(x), ~\varphi_A(y)\right\}, $ for any $x,y\in {\cal K}$ with $\langle x,y\rangle\ne 0$.
	\end{enumerate}
\end{proposition}

\begin{proof}
(a)$\Rightarrow$(b): First of all, we assume that item (a) holds. Arbitrarily take $x,y\in {\cal S}$. Thus, either $ q_A(x)\leq q_A(y)$ or $q_A(y)\leq q_A(x)$ holds. By using Proposition~\ref{pr:CharDiff} we conclude that 
\begin{align*}
	& ~ q_A(x)\leq q_A(y) \\
	\Rightarrow & ~\langle Dq_A(y), x\rangle - \langle x, y\rangle \langle Dq_A(y), y\rangle \le 0 \\
	\Rightarrow & ~\langle Ay,x\rangle\leq\langle x,y\rangle \langle Ay, y\rangle \\
	\Rightarrow & ~\langle Ay,x\rangle\leq\langle x,y\rangle q_A(y)
\end{align*}
Similarly, 
\[
	q_A(y)\leq q_A(x) \Rightarrow \langle Ax,y\rangle\leq\langle x,y\rangle q_A(x)
\]
 Thus, the symmetric matrix $A$ implies $\langle Ax,y\rangle= \langle Ay,x\rangle$, taking into account that ${\cal S} = \SP^{n-1}\cap\inte({\cal K})$ and ${\cal K}$ is a subdual cone and hence $\lng x,y\rng > 0$, we have
\[
	\langle Ax,y\rangle\leq\max\{\langle x,y\rangle q_A(x),\langle x,y\rangle q_A(y)\}=\langle x,y\rangle\max\{q_A(x), q_A(y)\}, \, \forall ~ x,y\in {\cal S}.
\]
Therefore, by continuity the above inequality can be extended to all $ x,y\in\SP^{n-1}\cap {\cal K}$ and, then item (b) holds. 

(b)$\Rightarrow$(a): Conversely, we assume that item (b) holds. Take $x,y\in {\cal S}$ satisfying $q_A(x)\leq q_A(y)$. Then, by the inequality in item (b) and the fact that ${\cal K}$ is a subdual cone, we have
\[
	q_A(x)\leq q_A(y)~\Rightarrow~\langle Ax,y\rangle - \langle x,y\rangle q_A(y)\leq 0.
\]

Hence, by using Proposition \ref{pr:CharDiff} we conclude that $q_A$ is a spherically quasi-convex function.

(b) $\Rightarrow$ (c): To establish the equivalence between (b) and (c), we firstly assume that item (b) holds. Let $x,y\in {\cal K}$ with $\langle x,y\rangle\ne 0$. Then, $x\ne 0$ and $y\ne 0$. Moreover, we have by the property of a cone: 
\[
	u:=\frac{x}{\|x\|} \in \SP^{n-1}\cap{\cal K}, \qquad v:=\frac{y}{\|y\|} \in \SP^{n-1}\cap {\cal K}.
\]
Hence, by using the inequality in item (b) with $x=u$ and $y=v$, we obtain the inequality in item (c). 

(c) $\Rightarrow$ (b): Conversely, suppose that (c) holds. 
Take $x,y\in \SP^{n-1}\cap {\cal K}$ with $\langle x,y\rangle\ne 0$. We have $\|x\|=\|y\|=1$ as $x, y\in \SP^{n-1}$. From the inequality in item (c) we conclude that
\[
	\frac{\langle Ax,y\rangle}{\langle x,y\rangle}\leq\max\left\{ q_A(x), ~ q_A(y) \right\}.
\]
Due to ${\cal K}$ being a subdual cone, we have $\langle x,y\rangle\geq 0$, and hence the last inequality is equivalent to the inequality in item (b). 

Now, assume that $\langle x,y\rangle=0$. Then, take two sequences $\{x_k\}, \{y_k\}\subseteq {\cal S}$ such that
$\lim_{k\to + \infty}x_k= x$, $\lim_{k\to + \infty}y_k= y$ and $\lng x_k,y_k\rng\ne 0$. Since ${\cal K}$ is a subdual cone, we have $\langle x_k,y_k\rangle>0$
for any $k\in\{1, 2, \ldots\}$. Therefore, considering that $\|x_k\|=\|y_k\|=1$ for any $k=1, 2, \ldots$, we can apply again the inequality in item
(c) to conclude
$$
\langle A x_k, y_k \rangle\leq\langle x_k, y_k\rangle\max \left\{q_A(x_k), ~ q_A(y_k)\right\}, \qquad k=1, 2, \ldots.
$$
By tending with $k$ to infinity, we conclude that the inequality in item (b) also holds for $\langle x,y\rangle=0$ and the proof of the equivalence between (b) and (c) is complete.
\hfill$\square$\par \end{proof}
 
\bigskip \begin{corollary}\label{pr:qpcf}
Assume that ${\cal K}$ is a self-dual cone. If the quadratic function $q_A$ is spherically quasi-convex, then $A$ has the $\cal K$-Z-property (see Definition \ref{def:k_z_p}).
\end{corollary}
\begin{proof}
To prove $A$ has the $\cal K$-Z-property, we need to prove that 
\[
	\langle Ax, y\rangle \le 0
\]
for any $(x,y)\in {\cal C}({\cal K})$. Take  $x,y \in\R^n$ such that $(x,y)\in {\cal C}({\cal K})$. If either $x=0$ or $y=0$, we have $\langle Ax,y\rangle=0$. Thus, assume that $x\neq 0$ and $y\neq 0$. Considering that ${\cal K}$ is a self-dual cone, we have $
\frac{x}{\|x\|},~\frac{y}{\|y\|}\in\SP^{n-1}\cap {\cal K}$. Suppose that $q_A$ is spherically quasi-convex and by the items (a) and (b) of Proposition~\ref{pr:spher-quasiconv}, we have
\[
	\Big\langle A\frac{x}{\|x\|}, \frac{y}{\|y\|}\Big\rangle \le \Big\langle \frac{x}{\|x\|}, \frac{y}{\|y\|} \Big\rangle \max \lf\{ q_A\lf(\frac{x}{\|x\|}\rg), q_A\lf(\frac{y}{\|y\|}\rg)\rg\}, \quad \forall x, y \in {\cal C}({\cal K}).
\]
By fact that $\big\langle \frac{x}{\|x\|}, ~\frac{y}{\|y\|} \big\rangle=0$, we obtain $\langle Ax,y\rangle\leq 0$. 
\hfill$\square$\par \end{proof}

\bigskip \begin{theorem}\label{th:quasiconv-iff}
The function $q_A$ defined in \eqref{eq:QuadFunc} is spherically quasi-convex if and only if $\varphi_A$ defined in \eqref{eq:RayleighFunction} is quasi-convex.
\end{theorem}
\begin{proof}
For any $c\in\R$, let $[q_A\leq c]:=\left\{x\in {\cal S}~:~q_A(x)\leq c\right\}$ and $ [\varphi_A\leq c]:=\{x\in\inte(\mathcal K):\varphi_A(x)\leq c\}$ be the sublevel sets of $q_A$ and $\varphi_A$, respectively. Let ${\cal K}_{[q_A\leq c]}$ be the cone spanned by $[q_A \leq c]$.
Since ${\cal S}=\SP^{n-1}\cap\inte({\cal K})$, we conclude that $x\in \inte \mathcal K $ if and only if $x/\|x\| \in {\cal S}$. Hence, by the definitions of $[q_A\leq c]$ and $ [\varphi_A\leq c]$ we obtain:
\begin{align*}
	{\cal K}_{[q_A\leq c]}& = \Big\{tx~:~x\in {\cal S}, q_A(x)\leq c, t\in[0,+\infty)\Big\} \\
	& = \lf\{x\in\inte(\mathcal K)~:~ q_A\lf(\frac{x}{\|x\|}\rg)\leq c\rg\} \\
	& = \{x\in\inte(\mathcal K)~:~ \varphi_A(x)\leq c\}. 
\end{align*}
That is:
\begin{equation}\label{eq:eqmain}
{\cal K}_{[q_A\leq c]} = [\varphi_A\leq c].
\end{equation}
Suppose that the quadratic function $q_A$ is spherically quasi-convex. Thus, from Poposition~\ref{pr:charb1} we conclude that $[q_A\leq c]$ is spherically convex for any $c\in \R$. Hence, it follows from Proposition~\ref{pr:ccs} that the cone ${\cal K}_{[q_A\leq c]}$ is convex and pointed, which implies from \eqref{eq:eqmain} that $[\varphi_A\leq c]$  is convex for any $c\in \R$. Therefore, again by Poposition \ref{pr:charb1} we conclude that $\varphi_A$ is quasi-convex. 
 
Conversely, suppose that $\varphi_A$ is quasi-convex. Thus,
 $[\varphi_A\leq c]$ is convex for any $c\in \R$. On the other hand, given ${\cal K}$ is a proper subdual cone, $\inte \mathcal K$ is pointed. Thus, we conclude that $[\varphi_A\leq c]\subseteq \inte \mathcal K$ is also a pointed cone. Hence, from \eqref{eq:eqmain} it follows that ${\cal K}_{[q_A\leq c]}$ is, again, a pointed convex cone. Hence, Proposition \ref{pr:ccs} implies that $[q_A\leq c]$ is spherically convex for any $c\in \R$. Therefore, by using Proposition~\ref{pr:charb1}, we conclude that $q_A$ is a spherically quasi-convex function.
\hfill$\square$\par \end{proof}

Let $c\in \R$, recall the definition \eqref{eq:sls} with $f\equiv \varphi_A$: 
\begin{align*}
	[\varphi_A\leq c]:=& \{x\in {\cal S} :\; \varphi_A(x)\leq c\}  \\
	 = & \{x\in {\cal S} :\; \langle A_cx,x\rangle\leq 0\}, \qquad A_c:=A - cI_n.
\end{align*}
	
\bigskip \begin{corollary}\label{cor:cor}
 The function $q_A$ is spherically quasi-convex if and only if for any $c\in\R$ the set $[\varphi_A\leq c]$ is convex.
\end{corollary}
\begin{proof}
Suppose that the quadratic function $q_A$ is spherically quasi-convex. Hence Theorem \ref{th:quasiconv-iff} implies that $\varphi_A$ is quasi-convex, and the sub-level set $[\varphi_A\leq c]$ is convex for any $c\in\R$ by Proposition \ref{pr:charb1}. Since $\{x\in\inte({\cal K}):~\langle A_{c}x,x\rangle< 0\} \neq \emptyset$, we conclude that
$$
\textrm{cl}\Big(\big\{x\in\inte({\cal K})~:~\langle A_{c}x,x\rangle\leq 0\big\}\Big)= \{x\in{\cal K}~:~\langle A_{c}x,x\rangle\leq 0\}.
$$
where ``cl$(\cdot)$'' is the topological closure operator of a set. Thus, considering that 
\begin{align*}
	[\varphi_{A}\leq c] & = \Big\{x\in\inte({\cal K})~:~\frac{\langle Ax,x\rangle}{\|x\|^2}\leq c\Big\} \\
	& =  \{x\in\inte({\cal K})~:~\langle Ax,x\rangle - c\|x\|^2 \leq 0\} \\
	& =  \{x\in\inte({\cal K})~:~\langle Ax,x\rangle - c\langle x, x\rangle \leq 0\} \\
	& =  \{x\in\inte({\cal K})~:~\langle Ax - cx, x\rangle \leq 0\} \\
	& =  \Big\{x\in\inte({\cal K})~:~\big\langle (A - cI_n)x, x\big\rangle \leq 0\Big\},
\end{align*}

we obtain that
$$
		\textrm{cl}\left([\varphi_A\leq c]\right)=\{x\in{\cal K}~:~\langle A_{c}x,x\rangle\leq 0\},
$$
Taking into account that $[\varphi_A\leq c]$ is convex, the set $\textrm{cl}\big([\varphi_A\leq c]\big)$ is also convex. 
\hfill$\square$\par \end{proof}

\section{Spherically quasi-convex quadratic functions on the spherical positive orthant} \label{sec:qcqfpo}

In this section we present some properties of a quadratic function defined in the spherical positive orthant, which corresponds to ${\cal K}=\R^n_+$ (therefore $\inte K = \R^n_{++}$). If $A$ has only one eigenvalue, it is easy to conclude that $q_A$ is spherically quasi-convex. However, suppose that the only eigenvalue of $A$ is $\lambda$, then $q_A(x) = \lambda$ for any $x\in {\cal S}$ (we say $q_A$ is constant in this case), which is meaningless to discuss. Therefore, {\it throughout this section we assume
that $A$ has at least two distinct eigenvalues}. By the definitions in \eqref{eq:QuadFunc} and \eqref{eq:RayleighFunction}, the domains of $q_A$ and $\varphi_A$ (${\cal S}$ and $\inte({\cal K})$, respectively) are given by
\begin{equation} \label{eq:KR++}
	 {\cal S}:=\SP^{n-1}\cap\R^n_{++}, \qquad \inte({\cal K}):=\R^n_{++},
\end{equation}
Next we present a technical lemma which will be useful in the sequel.
\bigskip \begin{lemma}\label{Lem:Basic}
Let $n\geq 2$ and $V=[v^1~ v^2~v^3~\cdots ~v^n] \in\R^{n\times n}$ be an orthogonal matrix, $A=V\Lambda V\tp$ and $\Lambda=\diag(\lambda_1, \ldots, \lambda_n)$. Assume that $\lambda_1 < \lambda_2 \leq \ldots \leq \lambda_n$. If $ v^1 \in \R^n_{+}$, then the sublevel set $[\varphi_A\leq c]$ is convex for any $c\notin (\lambda_2,\lambda_n)$.
\end{lemma}
\begin{proof}
 By using that $V\tp V=I_n$ and $A=V\Lambda V\tp $ we obtain from the definition \eqref{eq:RayleighFunction} that

\begin{align}
	[\varphi_A\leq c] & = \left\{x\in\R^n_{++}:~ \frac{\langle Ax, x\rangle}{\|x\|^2} \leq c\right\}\nonumber\\
	& = \left\{x\in\R^n_{++}:~ \langle Ax, x\rangle - c\|x\|^2 \leq 0\right\}\nonumber\\
	& = \left\{x\in\R^n_{++}:~ \langle (A - cI_n)x, x\rangle\leq 0\right\}\nonumber\\
	& = \Big\{x\in\R^n_{++}:~ x\tp V\lf(\Lambda - cI_n\rg)V\tp x\leq 0\Big\}\nonumber\\
	& = \Big\{x\in\R^n_{++}:~ \sum_{i=1}^n(\lambda_i-c)\lf(x\tp v^i\rg)^2 \leq 0\Big\}\nonumber\\
	& = \left\{x\in\R^n_{++}:~ \sum_{i=1}^n(\lambda_i-c)\langle v^i, x\rangle^2 \leq 0\right\}\label{eq:bl1}
\end{align}

In the following we will show that $[\varphi_A\leq c]$ is convex for any $c\notin (\lambda_2,\lambda_n)$. 

If $c<\lambda_1$, then since $v^1, v^2, \ldots, v^n$ are linearly independent, we conclude from \eqref{eq:bl1} that $[\varphi_A\leq c]=\{0\}$ and therefore it is convex. 

If $c=\lambda_1$, then from \eqref{eq:bl1} we conclude that $[\varphi_A\leq c]={\cal O}\cap \R^n_{++},$ where ${\cal O}: =\{x\in\R^n ~: \langle v^2, x\rangle=0, ~\ldots,~ \langle v^n, x\rangle=0\}$ is a convex cone. Hence $[\varphi_A\leq c]$ is convex. 

If $\lambda_1<c\le \lambda_2$, letting $y=V^\top x$, i.e., $y_i=\langle v^i, x\rangle$, for $i=1, \ldots, n$. Since $v^1\in {\mathbb R^n_{++}}$ and $x\in\R^n_{++}$, we have $y_1 = \langle v^1, x\rangle >0$. From \eqref{eq:bl1}
we obtain
\begin{align*}
	[\varphi_A\leq c] & = \left\{x\in\R^n_{++}:~ (c - \lambda_1)\langle v^1, x\rangle^2 \ge \sum_{i=2}^n(\lambda_i-c)\langle v^i, x\rangle^2 \right\}\\
	& = \left\{x\in\R^n_{++}:~ \langle v^1, x\rangle^2 \ge \sum_{i=2}^n\theta_i\langle v^i, x\rangle^2 \right\}
\end{align*}
where $\theta_i=\frac{\lambda_i-c}{c-\lambda_1}$ for $i=2, \ldots, n$. Denote the cone
\[
	{\cal L}:=\left\{y=(y_1, \ldots, y_n)\in\R^n:~y_1\geq \sqrt{\theta_2y_2^2+ \ldots + \theta_{n}y_n^2}\right\},
\]
we have $[\varphi_A\leq c]= {\cal L}\cap V\tp\R^n_{++}$. Since both ${\cal L}$ and $V\tp\R^n_{++}$ are convex sets, we conclude that $[\varphi_A\leq c]$ is convex. 

If $c\ge \lambda_n$, then $[\varphi_A\leq c]=\R^n_{++}$ is convex.
\hfill$\square$\par \end{proof}

It should be mentioned that if $ \lambda_1< \lambda_2=\dots=\lambda_n$, then the sublevel set $[\varphi_A\leq c]$ is convect for any $c\in\R$. This will be a useful fact in a proof of a following theorem. 

\bigskip \begin{lemma} \label{lem:Copositive}
	Let $\lambda $ be an eigenvalue of $A$. If $\lambda I_n-A$ is copositive and $\lambda\leq c$, then 
	\[
		[\varphi_A\leq c]=\R^n_{++}
	\] 
	and consequently it is a convex set.
\end{lemma}
\begin{proof}
Let $c\in \R$ and $[\varphi_A\leq c]=\{x\in \R^n_{++}:~ \langle A x,x \rangle-c \|x\|^2\leq 0\}$. Suppose that $\lambda\leq c$, for any $x\in \R^n_{++}$ we have 
\[
	\langle A x,x \rangle-c\|x\|^2\leq \langle A x,x \rangle-\lambda \|x\|^2= \langle (A-\lambda I_n) x,x \rangle,
\]
and $\lambda I_n-A$ is copositive, that is
\[
	\langle A x,x \rangle-c\|x\|^2\leq\langle (A-\lambda I_n) x,x \rangle\leq 0,
\]
hence $\langle A x,x \rangle-c\|x\|^2\leq 0$ holds for any $x\in \R^n_{++}$, which implies that $[\varphi_A\leq c]=\R^n_{++}$.
\hfill$\square$\par \end{proof}
The next theorem exhibits a series of implications and, in particular, conditions which imply that the quadratic function $q_A$ is spherically quasi-convex.

\bigskip \begin{theorem}\label{Thm:Copositive}
	Let $A\in\mathbb R^{n\times n}$ be a symmetric matrix and let $\lambda_1\le\lambda_2\leq\dots\leq\lambda_n$ its eigenvalues. Consider the following statements:
	\begin{enumerate}
		\item[(i)] $q_A$ is a spherically quasi-convex function.
		\item[(ii)] $A$ is a Z-matrix.
		\item[(iii)] $A$ is a Z-matrix and $\lambda_2\geq a_{ii}$ for any $i\in\{1,2,\ldots,n\}$.
		\item[(iv)] $A$ is a Z-matrix, $\lambda_1<\lambda_2$ and $\lambda_2\geq a_{ii}$ for any $i \in \{1,2,\ldots,n\}$.
		\item[(v)] $A$ is an irreducible Z-matrix and $\lambda_2\geq a_{ii}$ for any $i\in\{1,2,\ldots,n\}$.
		\item[(vi)] $\lambda_2I_n-A$ is copositive and there exists an eigenvector $v^1\in\R^n_{++}$ corresponding to the eigenvalue $\lambda_1$ of $A$.
   \end{enumerate}

	Then the following implications hold:
  $$
  \begin{array}{ccccccc}
	  &&\textrm{(iv)}		&&\\
	&&\Downarrow && \\
	\textrm{(iii)} & \Leftarrow & \textrm{(vi)} & \Rightarrow &\textrm{(i)} & \Rightarrow & \textrm{(ii)} \\
	&&\Uparrow && \\
	&&\textrm{(v)}		&&\\
  \end{array}
  $$
\end{theorem}

\begin{proof}

$\,$
\vspace{2mm}

(iv)$\Rightarrow$(vi):
Suppose that  $A$ is a Z-matrix, $\lambda_1<\lambda_2$ and $\lambda_2\geq a_{ii}$ for any $i \in \{1,2,\ldots,n\}$. It is easy to verify that $\lambda_2I_n-A$ is a nonnegative matrix, and hence $\lambda_2I_n-A$ is copositive. Moreover, applying item (i) of Perron-Frobenius theorem (see Theorem \ref{Perron-Frobenius theorem1}) to the matrix
$\lambda_2I_n-A$, we obtain that there exists an eigenvector $v^1\in\mathbb R^n_{++}$ of $\lambda_2I_n-A$ corresponding to the largest eigenvalue $\lambda_2-\lambda_1$ ($v^1$ is also the eigenvector of $A$ corresponding to $\lambda_1$).

(v)$\Rightarrow$(vi):
Suppose that $A$ is an irreducible Z-matrix and $\lambda_2\geq a_{ii}$ for any $i\in\{1,2,\ldots,n\}$. Again, by applying item (i) of Perron-Frobenius theorem we conclude that there exists an eigenvector $v^1\in\mathbb R^n_{++}$ of $\lambda_2I_n-A$ corresponding to the largest eigenvalue $\lambda_2-\lambda_1$ ($v^1$ is also the eigenvector of $A$ corresponding to $\lambda_1$).

(vi)$\Rightarrow$(i): If $c\le\lambda_2$, then Lemma~\ref{Lem:Basic} implies that $[\varphi_A\leq c]$ is convex. If $c\geq \lambda_2$, then from Lemma~\ref{lem:Copositive} we have $[\varphi_A\leq c]=\R^n_{++}$, which is convex. Hence, $[\varphi_A\leq c]$ is convex for any $c\in \R$. Therefore, by using Theorem~\ref{th:quasiconv-iff}, we conclude that $q_A$ is spherically quasi-convex function.

(i)$\Rightarrow$(ii): Suppose that $q_A$ is spherically quasi-convex. From Corollary \ref{pr:qpcf}, it follows that $A$ has the $\R^n_+$-Z-property. By applying Theorem \ref{thm:k_z_z_m}, we obtain that $A$ is a Z-matrix.

(vi)$\Rightarrow$(iii): Suppose (vi) holds, by above proof (vi)$\implies$(i)$\implies$(ii), it follows that $A$ is a Z-matrix. Arbitrarily take $x = e^i$ with $i\in \{1,2,\dots, n\}$, $e^i$ are canonical vectors of $\R^n$. Since $\lambda_2I_n-A$ is copositive, it follows that
\[
	\langle (\lambda_2 I_n - A)x, x \rangle \ge 0 \Rightarrow \lambda_2  - a_{ii} \ge 0
\]
for any  $i\in \{1,2,\dots, n\}$. 

\hfill$\square$\par \end{proof}

\bigskip \begin{corollary}\label{Cor:NegMatrix}
 Let $n\ge 2$ and  $ \lambda_1, \ldots , \lambda_n \in \R$ be the eigenvalues of $A$.  Assume that  $-A$ is  an entrywise positive matrix, $ \lambda_1 < \lambda_2 \leq \ldots \leq \lambda_n$ and $0< \lambda_2$ . Then $q_A$ is spherically quasi-convex.
\end{corollary}
\begin{proof}
Suppose that  $-A$ is an entrywise positive matrix, $ \lambda_1 < \lambda_2 \leq \ldots \leq \lambda_n$ and $0< \lambda_2$, then the matrix $\lambda_2I_n - A$ is an entrywise positive matrix with $\lambda_2-\lambda_1>0$ to be its largest eigenvalue. Thus, Theorem~\ref{Perron-Frobenius theorem1} implies that the eigenvalue $\lambda_2-\lambda_1$ has the associated eigenvector $v^1\in \R^n_{++}$. Since by the definition of eigenvector
\[
	(\lambda_2I_n - A)v^1=(\lambda_2-\lambda_1)v^1 \Rightarrow Av^1=\lambda_1v^1,
\]
we conclude that $ v^1$ is also an eigenvector of $A$ associated to $\lambda_1$. Therefore, considering that $A$ is a entrywise negative matrix, then $A$ is also a Z-matrix. Since $ v^1 \in \R^n_{++}$, $\lambda_1 < \lambda_2$ and $\lambda_2\ge a_{ii}$ for any $i\in\{1,2\dots,n\}$, it follows from Theorem~\ref{Thm:Copositive} (iv)$\Rightarrow$(i) that $q_A$ is spherically quasi-convex.
\hfill$\square$\par \end{proof}

In the following two examples we use Theorem~\ref{Thm:Copositive} (vi)$\Rightarrow$(i) to illustrate a class of quadratic quasi-convex functions defined in the spherical positive orthant.
\bigskip \begin{example}\label{Ex:countergg}
Let $n\geq 3$ and $V=[v^1~ v^2~v^3~\cdots ~v^n] \in\R^{n\times n}$ be an orthogonal matrix, $A=V\Lambda V\tp$ and $\Lambda:=\diag(\lambda,\mu, \ldots,\mu ,\nu)$, where $\lambda,\mu,\nu\in\R$. Then $q_A$ is a spherically quasi-convex if
\begin{equation} \label{eq:ciqc}
v^1-\sqrt{\frac{\nu-\mu}{\mu-\lambda}}|v^{n}|\in {\mathbb R^n_{+}}, \qquad \lambda<\mu<\nu,
\end{equation}
where $|v^{n}|:=(|v_1^{n}|, \ldots, |v_n^{n}|)$. Indeed, by using that $V^\top V=I_n$ and $A=V\Lambda V\tp$, then 
\begin{align}
	\mu\|x\|^2 - \langle A x,x \rangle& = \mu\|x\|^2 - \langle V\Lambda V\tp x, x\rangle \nonumber\\
	& = \mu \|x\|^2 - \lf(x\tp v^1\rg)^2 \lambda + \sum_{i=2}^{n -1} \lf(x\tp v^i\rg)^2 \mu + \lf(x\tp v^n\rg)^2 \nu\nonumber\\
	& = \mu \langle I_nx,x\rangle - \lambda \langle v^1, x\rangle^2 + \sum_{i=2}^{n-1}\mu \langle v^i, x\rangle^2 + \nu \langle v^n, x\rangle ^2 \nonumber \\
	& = \mu \langle V\tp V x,x\rangle - \lambda \langle v^1, x\rangle^2 + \sum_{i=2}^{n-1}\mu \langle v^i, x\rangle^2 + \nu \langle v^n, x\rangle ^2 \nonumber \\
	& = \sum_{i=1}^{n} \mu \langle v^i,x\rangle^2 - \lambda \langle v^1, x\rangle^2 + \sum_{i=2}^{n-1}\mu \langle v^i, x\rangle^2 + \nu \langle v^n, x\rangle ^2 \nonumber \\
	& = (\mu-\lambda) \left[\langle v^1, x\rangle^2 - \frac{\nu-\mu}{\mu-\lambda}\langle v^n, x\rangle^2\right].\label{eq:ceq}
\end{align}

 Thus, using the condition in \eqref{eq:ciqc} and arbitrarily take $x\in\R^n_{++}$, we have
\begin{equation*}
\langle v^1, x\rangle^2 - \frac{\nu-\mu}{\mu-\lambda}\langle v^n, x\rangle^2 \geq \frac{\nu-\mu}{\mu-\lambda}\left[\langle |v^n|, x\rangle^2 - \langle v^n, x\rangle^2\right] \geq 0.
\end{equation*}
Hence, by combining the last inequality with  \eqref{eq:ceq}, we conclude that 
\[
	\langle (\mu I_n-A)x,x\rangle \ge 0
\]
for any $x\in\R^n_{++}$. Hence $\mu I_n-A$ is copositive. Therefore, since $ v^1 \in \R^n_{+}$ we can apply
Theorem~\ref{Thm:Copositive} (vi)$\Rightarrow$(i) with $\lambda_2=\mu$ to conclude that $q_A$ is a spherically quasi-convex function. 

For instance, taking $\lambda<(\lambda+\nu)/2<\mu<\nu$, and the vectors $v^1=(e^1+e^n)/\sqrt{2}, v^2=e^2, ~ \ldots, ~ v^{n-1}=e^{n-1}, v^n=(e^1-e^n)/\sqrt2$, satisfy \eqref{eq:ciqc}. We can conclude that $q_A$ is spherically quasi-convex. 
\end{example}

\bigskip \begin{example}\label{Ex:CountExMulEngfm}
Let $n\geq 3$ and $V=[v^1~ v^2~v^3~\cdots ~v^n] \in\R^{n\times n}$ be an orthogonal matrix, {$\Lambda=\diag(\lambda_1, \ldots, \lambda_n)$ and $A=V\Lambda V\tp$ }. Then $q_A$ is a spherically quasi-convex function, if
\begin{equation} \label{eq:CondEngfm}
 v^1=(v^1_1, \ldots, v^1_n)\tp\in {\mathbb R^n_{++}}, \qquad \lambda_1<\lambda_2\leq \cdots \leq \lambda_n \leq\lambda_2 +\frac{\alpha^2}{(n-2)}(\lambda_2-\lambda_1),
\end{equation}
where $\alpha:=\min\big\{ v^1_i :~ i\in\{1, \ldots, n\}\big\}$ is the minimum (nonzero) element of eigenvector $v^1$.
Indeed, by using $V\tp V = I_n$ and the definition of the matrix $A$, we obtain 
\begin{align}
	\lambda_{2}\|x\|^2 - \langle A x,x \rangle & = \lambda_2\langle V\tp V x,x\rangle - \langle V\Lambda V\tp x, x\rangle \nonumber \\
	& = \sum_{i=1}^{n} \lambda_{2} \langle v^i,x\rangle^2 - \sum_{i=1}^{n}\lambda_{i} \langle v^i, x\rangle^2 \nonumber \\
	& = \lf(\lambda_{2} -\lambda_1\rg)\langle v^1,x\rangle^2 + \sum_{i=3}^{n} \lf(\lambda_{2} -\lambda_i\rg)\langle v^i,x\rangle^2. \label{eq:ExCopo}
\end{align}

Since $\lambda_{2}-\lambda_1>0$ and $\lambda_2-\lambda_{n}\leq \lambda_2 - \lambda_j\leq 0$, for any $j\in\{3, \ldots, n\}$, from \eqref{eq:ExCopo} we have
\begin{align}
	\lambda_{2}\|x\|^2 - \langle A x,x \rangle & = \lf(\lambda_{2} -\lambda_1\rg)\langle v^1,x\rangle^2 + \sum_{i=3}^{n} \lf(\lambda_{2} -\lambda_i\rg)\langle v^i,x\rangle^2 \nonumber \\
	& = \lf(\lambda_{2} -\lambda_1\rg)\lf[\langle v^1,x\rangle^2 + \sum_{i=3}^{n} \frac{\lf(\lambda_{2} -\lambda_i\rg)}{\lf(\lambda_{2} -\lambda_1\rg)}\langle v^i,x\rangle^2\rg] \nonumber \\
	& \geq (\lambda_{2}-\lambda_1)\left[\langle v^1, x\rangle^2 + \sum_{i= 3}^{n} \frac{\lambda_2-\lambda_{n}}{\lambda_{2}-\lambda_1} \langle v^{i}, x\rangle^2\right]. \label{eq:coprfm}
\end{align} 

On the other hand, by using that $v^1_i\in {\mathbb R_{++}}$ and $v^1_i \geq \alpha$ for any $i\in\{1, \ldots, n\}$, we conclude that
\begin{align}\label{eq:CopFi}
 \langle v^1, x\rangle^2 & = (v^1_1 x_1+ \cdots + v^1_n x_n)^2 \nonumber \\
 & \geq \alpha^2 ( x_1+ \cdots + x_n)^2\geq \alpha^2 ( x_1^2+ \cdots + x_n^2)=\alpha^2\|x\|^2,
 \end{align}

for any $x\in\R^n_{+}$. Moreover, taking into account that $\|v^j\|=1$ for any $j\in\{3, \ldots, n\}$, applying Cauchy-Schwarz ineqauality, it follows that
$$
 \langle v^{3}, x\rangle^2+ \cdots + \langle v^n, x\rangle^2\leq \| v^{3}\|^2\|x\|^2+\cdots+ \| v^{n}\|^2\|x\|^2 \leq (n-2)\|x\|^2,
$$
for any $x\in\R^n_{+}$. Thus, combining the last inequalities with \eqref{eq:coprfm} and \eqref{eq:CopFi} and considering that the last inequality in \eqref{eq:CondEngfm} is equivalent to  $ -\alpha^2+ (n-2)(\lambda_n-\lambda_{2})/(\lambda_{2}-\lambda_1) \leq 0$, we have
\begin{align*}
	\lambda_{2}\|x\|^2 - \langle A x,x \rangle & \geq (\lambda_{2}-\lambda_1)\left[\langle v^1, x\rangle^2 + \sum_{i= 3}^{n} \frac{\lambda_2-\lambda_{n}}{\lambda_{2}-\lambda_1} \langle v^{i}, x\rangle^2\right] \\
	& \geq (\lambda_{2}-\lambda_1)\left[\alpha^2\|x\|^2 + \sum_{i= 3}^{n} \frac{\lambda_2-\lambda_{n}}{\lambda_{2}-\lambda_1} \langle v^{i}, x\rangle^2\right] \\
	& \geq (\lambda_{2}-\lambda_1)\left[\alpha^2 + (n-2) \frac{\lambda_2-\lambda_{n}}{\lambda_{2}-\lambda_1}\right]\|x\|^2\geq 0,
\end{align*}
for any $x\in\R^n_{+}$. Hence, we conclude that $\lambda_{2} I_n-A$ is copositive. Therefore, since $v^1\in{\mathbb R}^n_{++}$ is the eigenvector of $A$ corresponding to the eigenvalue $\lambda_1$, we apply Theorem~\ref{Thm:Copositive} (vi)$\Rightarrow$(i), to conclude that $q_A$ is a spherically quasi-convex function. 

For instance, $n\geq 3$, $A=V\Lambda V\tp$, $\Lambda=\diag(\lambda_1, \ldots, \lambda_n)$, $V=[v^1~ v^2~v^3~\cdots ~v^n] \in\R^{n\times n}$, and $\alpha=1/\sqrt{n}$,
\small
$$
v^1:= \frac{1}{\sqrt{n}} \sum_{i=1}^ne^i, \quad v^j:= \frac{1}{\sqrt{(n+1-j)+(n+1-j)^2}}\left[e^1- (n+1-j)e^j + \sum_{i>j}^ne^i \right],
$$
\normalsize
for $j\in\{2, \ldots, n\}$ and $\lambda_1<\lambda_2\leq \ldots \leq \lambda_n <\lambda_2 +(1/[n(n-2)])(\lambda_2-\lambda_1)$, satisfy the
orthogonality of $V$ and the condition \eqref{eq:CondEngfm}, therefore we conclude that $q_A$ is a spherically quasi-convex.
\end{example}

In the next theorem we establish the characterisation for a quasi-convex quadratic function $q_{A}$ on the spherical positive orthant, where $A$ is a symmetric matrix having only two distinct eigenvalues.

\bigskip \begin{theorem}\label{thm:2neg-v}
Let $n\ge 3$ and $A\in\R^{n\times n}$ be a symmetric matrix with only two distinct eigenvalues, such that its smallest eigenvalue has multiplicity one. Then, $q_{A}$ is spherically quasi-convex if and only if there is an eigenvector of $A$ corresponding to the smallest eigenvalue with all components nonnegative.
\end{theorem}
\begin{proof}
Let $A:=(a_{ij})\in\R^{n\times n}$, $\lambda_1, \lambda_2,\dots, \lambda_n$ be the eigenvalues of $A$ corresponding to an orthonormal set of eigenvectors $v^1, v^2,\dots,v^n$, respectively. Then, we can assume with no loss of generality that $\lambda_1=:\lambda <\mu:=\lambda_2=\dots=\lambda_n. $
Thus, we have
\begin{equation} \label{eq:VLVT}
A=V\Lambda V^T, \qquad V:=[v^1\textrm{ }v^2\textrm{ }\dots\textrm{ }v^n]\in\R^{n\times n}, \qquad \Lambda:= \diag(\lambda,\mu, \ldots, \mu) \in\R^{n\times n}.
\end{equation}
First we suppose that $q_{A}$ is a spherically quasi-convex function. The matrix $\Lambda$ can be equivalently written as follows
\begin{equation} \label{eq:eqtr}
 \Lambda=\mu I_n+(\lambda-\mu)D,
\end{equation}
where $D:=\diag(1,0, \dots, 0)\in\R^{n\times n}$. Then \eqref{eq:eqtr} and \eqref{eq:VLVT} imply
\begin{equation} \label{eq:neg-v}
a_{ij}=(\lambda-\mu)v^1_iv^1_j,\, \qquad i\ne j.
\end{equation}

Since $q_{A}$ is spherically quasi-convex and $e^i \in {\cal S}=\SP^{n-1}\cap\R^n_{++}$ for any $i\in\{1, \ldots, n\}$, by using Theorem~\ref{Thm:Copositive} (ii)$\Rightarrow$(i) we conclude that the matrix $A$ is a Z-matrix, that is $a_{ij}\leq 0$ for any $i, j\in\{1, \ldots, n\}$ with $i\ne j$. Thus, owing that $\lambda <\mu
$, we obtain form \eqref{eq:neg-v} that $0\leq v^1_iv^1_j $ for any $i\ne j$, which implies $v^1\in \R^n_+$ (or $-v^1\in \R^n_+$, they are the same because $v^1$ is an eigenvector). Therefore, there is an eigenvector corresponding to the smallest eigenvalue with all components nonnegative. 

Conversely, suppose that $v^1\in \R^n_{+}$. Then, applying Lemma~\ref{Lem:Basic} with $ \lambda =\lambda_1<\mu=\lambda_2=\dots=\lambda_n$ we conclude that $[\varphi_A\leq c]$ is convex for any $c\in\R$, and hence $\varphi_A$ is quasi-convex by Proposition \ref{pr:charb1}. Therefore, by using Theorem~\ref{th:quasiconv-iff}, we conclude that $q_{A}$ is spherically quasi-convex.
\hfill$\square$\par \end{proof}

In the following examples we present two classes of matrices satisfying the assumptions of Theorem~\ref{thm:2neg-v}.

\bigskip \begin{example}
Let $v\in \R^n_{+}$ and define the Householder matrix $H:=I_n-2vv\tp/\|v\|^2$. The matrix $H$ is nonsingular and symmetric. Moreover, the Householder matrix has a property that $Hv=-v$. Given that $\frac{vv\tp}{\|v\|^2}$ is a Rank-1 matrix, the characteristic polynomial of $H$ is
\begin{align*}
	p(\lambda) & = |\lambda I_n - H| = (\lambda - 1)^n + Tr\lf(2\frac{vv\tp}{\|v\|^2}\rg)(\lambda - 1)^{n-1} \\
	& = |\lambda I_n - H| = (\lambda - 1)^n + Tr\lf(2\frac{v\tp v}{\|v\|^2}\rg)(\lambda - 1)^{n-1} \\
	& = |\lambda I_n - H| = (\lambda - 1)^n + 2(\lambda - 1)^{n-1} \\
	& = |\lambda I_n - H| = (\lambda + 1)(\lambda - 1)^{n-1},
\end{align*}
we conclude that $-1$ and $1$ are eigenvalues of $H$ with multiplicities one and $n-1$, respectively.

Furthermore, the eigenvector corresponding to the smallest eigenvalue of $H$ has all components nonnegative. Therefore, Theorem~\ref{thm:2neg-v} implies that the quadratic function $q_{H}(x)=\langle Hx,x\rangle$ is spherically quasi-convex.
\end{example}\index{Householder matrix}

{
\bigskip \begin{example}
	Let $u\in \R^n$, $\alpha \in \R_{++}$. Define the matrix $A:= I_n -\alpha uu\tp$. The eigenvalues of $A$ are $\{\lambda_1, \lambda_2,\dots, \lambda_n\} = \{\alpha\langle u, u\rangle + 1, 1, \cdots, 1\}$. Thus, we have
\begin{equation}\label{eq:VLVTex}
	A=V\Lambda V^T, \qquad V:=[v^1\textrm{ }v^2\textrm{ }\dots\textrm{ }v^n]\in\R^{n\times n}, \qquad \Lambda:= \diag(1-\alpha\langle u, u, \rangle ,\mu, \ldots, \mu) \in\R^{n\times n}.
\end{equation}

The matrix $\Lambda$ can be equivalently written as follows
\begin{equation} \label{eq:eqtrex}
 \Lambda= I_n - \alpha \langle u, u \rangle D,
\end{equation}
where $D:=\diag(1,0, \dots, 0)\in\R^{n\times n}$. Then \eqref{eq:eqtrex} and \eqref{eq:VLVTex} imply
\[
a_{ij}=\lf\{\begin{array}{ll}
	- \alpha \langle u, u \rangle v^1_iv^1_j, & \qquad i\ne j, \\
	- \alpha \langle u, u \rangle (v^1_i)^2 +1, & \qquad i= j. 
	\end{array}
	\rg.
\]

Therefore, if $v^1\in\R^n_+$, Theorem \ref{thm:2neg-v} implies that the function $q_A =\langle Ax, x\rangle$ is spherically quasi-convex.

Also, if $v^1\in\R^n_+$, then $A$ is a $Z$-matrix, we have $\lambda_1<\lambda_2$ and $\lambda_2\ge a_{ii}$ for any $i\in\{1, 2, \dots, n\}$. Theorem \ref{Thm:Copositive} (iv)$\Rightarrow$ (i) implies that the function $q_A =\langle Ax, x\rangle$ is spherically quasi-convex.
\end{example}
}

In order to give a {\it complete characterisation of the spherical quasi-convexity of $q_A$ for the case when $A$ is diagonal}, in the following result we
start with a necessary condition for $q_A$ to be spherically quasi-convex on the spherical positive orthant.
\bigskip \begin{lemma}\label{Lem:DiagPosCond}
Let $n\ge 3$, ${\cal S}=\SP^{n-1}\cap \R^n_{++}$ and $A\in\R^{n\times n}$ be a nonsingular diagonal matrix. If $q_{A}$ is spherically
quasi-convex, then $A$ has only two distinct eigenvalues, such that its smallest one has multiplicity one.
\end{lemma}
\begin{proof}
The proof will be made by contradiction. First we suppose that $A$ has at least three distinct eigenvalues, among which exactly two are negative, or at least two distinct eigenvalues, among which exactly one is negative and has multiplicity greater than one, i.e.,
\begin{equation} \label{eq:DefVec}
Ae^1= -\lambda_1 e^1, \qquad Ae^2= -\lambda_2 e^2, \quad Ae^3= \lambda_3 e^3, \qquad \lambda_1, \lambda_2, \lambda_3 >0
\end{equation}
with either $-\lambda_1< -\lambda_2<0< \lambda_3$ or $-\lambda_1= -\lambda_2<0<\lambda_3$ and $e^1, e^2,e^3$ are canonical vectors of $\R^n$.
Define the following two auxiliary vectors
\begin{equation} \label{eq:AuxVec}
v^1:=e^1+t_1e^3, \qquad v^2:=e^2+t_2e^3, \qquad t_i=\sqrt{\frac{\lambda_i}{\lambda_3}}, \qquad i\in\{1,2\}.
\end{equation}
Hence, \eqref{eq:DefVec} and \eqref{eq:AuxVec} implies that 
\[
	\langle Av^1, v^1\rangle = \langle Ae^1, e^1\rangle + t_1^2\langle Ae^3, e^3\rangle = -\lambda_1\|e^1\| + \frac{\lambda_1}{\lambda_3}\lambda_3\|e^3\| = 0,
\]
and 
\[
	\langle Av^2, v^2\rangle = \langle Ae^2, e^2\rangle + t_2^2\langle Ae^3, e^3\rangle = -\lambda_2\|e^2\| + \frac{\lambda_2}{\lambda_3}\lambda_3\|e^3\| = 0.
\]
Since $v^1, v^2\in \R^n_{+}$, we conclude that $v^1, v^2 \in \left\{x\in\R^n_{+}~:~\langle Ax, x\rangle \leq 0\right\}$. However, using again \eqref{eq:DefVec} and \eqref{eq:AuxVec} we obtain that
$$
\langle A(v^1+v^2), v^1+v^2 \rangle =2\langle Av^1, v^2\rangle = 2\frac{\sqrt{\lambda_1\lambda_2}}{\lambda_3}\lambda _3 = 2\sqrt{\lambda_1 \lambda_2} >0,
$$
and therefore $v^1+v^2\notin \left\{x\in\R^n_{+}~:~\langle Ax, x\rangle \leq 0\right\}$. We conclude that $\left\{x\in\R^n_{+}~:~\langle Ax, x\rangle \leq 0\right\}$ is not a convex cone. 

Finally, suppose that $A$ has at least three distinct eigenvalues or at least two distinct ones with the smallest one having multiplicity greater than one. Let $\lambda, \mu, \nu$ be eigenvalues of $A$ such that either 
$\lambda < \mu < \nu$ or $\lambda = \mu < \nu$. Take a constant $c\in \R$ such that $\mu<c<\nu$. Letting $A_{c}:=A-cI_n$ we conclude
that $\lambda-c, \mu-c, \nu-c$ are eigenvalues of $A_{c}$ and satisfying
\[
\lambda -c<\mu-c<0<\nu-c
\]
or
\[
	\lambda -c=\mu-c<0<\nu-c.
\]

Thus, by the first part of the proof, with $A_c$ in the role of $A$, we conclude that 
\[
	\left\{x\in\R^n_{+}~:~\langle A_cx, x\rangle \leq 0\right\}
\]
is not a convex cone. On the other hand, due to $e^i\in \R^n_{+}$ and $ \langle Ae^i, e^i\rangle =\lambda-c<0$, for some $i$, we obtain that 
\[
	\left\{x\in\R^n_{++}~:~\langle A_cx, x\rangle < 0\right\}\neq \emptyset.
\]

Henceforth, applying Corollary~\ref{cor:cor} with ${\cal K}=\R^n_{+}$ and taking into account that both 
\[
	\left\{x\in\R^n_{+}~:~\langle Ax, x\rangle \leq 0\right\}
\]
 in the first part of the proof, and
\[
	\left\{x\in\R^n_{+}~:~\langle A_cx, x\rangle \leq 0\right\}
\]
in the second part of the proof are not convex, we conclude that $q_{A}$ is not spherically quasi-convex.
\hfill$\square$\par \end{proof}

To make the study self-contained we state the result of \cite[Theorem 1]{FerreiraNemeth2017} explicitly here:
\bigskip \begin{theorem} \label{eq:cpoth}
	\cite[Theorem 1]{FerreiraNemeth2017} Let ${\cal S}=\SP^{n-1} \cap \R^n_{++}$ and $A\in\R^{n\times n}$ be a symmetric matrix. Then, $q_A$ is spherically convex if and only if
	there exists $\lambda \in\R$ such that $A=\lambda I_{n}$. In this case $q_A$ is a constant function.
\end{theorem}

The next result gives a full characterisation for $q_A$ to be spherically quasi-convex quadratic function on the spherical positive orthant, where {\it $A$ is a diagonal matrix}. The proof of this result is a combination of Theorem~\ref{thm:2neg-v}, Lemma~\ref{Lem:DiagPosCond} and Theorem \ref{eq:cpoth}. Before presenting the result we need the following definition: 

\bigskip \begin{definition}[Merely spherically quasi-convex] 
	A function is called \emph{merely spherically quasi-convex} if it is spherically quasi-convex, but it is not spherically convex.
\end{definition}\index{Spherically convex!Merely spherically quasi-convex}

\bigskip \begin{theorem}\label{cr:main}
Let $n\ge 3$ and $A\in\R^{n\times n}$ be a nonsingular diagonal matrix. Then $q_A$ is merely spherically quasi-convex if and only if $A$ has only two eigenvalues, such that its smallest one has multiplicity one and has a corresponding eigenvector with all components nonnegative.
\end{theorem}
\begin{proof}
Given that $A\in\R^{n\times n}$ is a nonsingular diagonal matrix. We suppose $q_A$ is a merely spherically quasi-convex function defined on ${\cal S} = \SP^{n-1}\cap \R^n_{++}$. By applying Lemma \ref{Lem:DiagPosCond}, we conclude that $A$ has only two eigenvalues, such that its smallest one has multiplicity one and has a corresponding eigenvector with all components nonnegative.

Conversely, suppose that $A$ is a nonsingular diagonal (therefore, symmetric) matrix and has only two eigenvalues, such that its smallest one has multiplicity one and has a corresponding eigenvector with all components nonnegative. By Theorem \ref{eq:cpoth} we conclude that $q_A$ is not spherically convex because $A\neq \lambda I_n$ for any $\lambda\in \R$. Using Theorem~\ref{thm:2neg-v}, we obtain that $q_A$ is spherically quasi-convex. Therefore $q_A$ is merely spherically quesi-convex.
\hfill$\square$\par \end{proof}

We end this section by showing that, if a symmetric matrix $A$ has three eigenvectors in the nonnegative orthant associated to at least two distinct
eigenvalues, then the associated quadratic function $q_A$ cannot be spherically quasi-convex.
\bigskip \begin{lemma}\label{cor:DiagNegCond}
Let $n\ge 3$ and $v^1, v^2, v^3 \in \R^n$ be distinct eigenvectors of a symmetric matrix $A$ associated to the eigenvalues $\lambda_1,
\lambda_2, \lambda_3 \in \R$, respectively, among which at least two are distinct. If $q_A$ is spherically quasi-convex, then $v^i \notin \R^n_{+}$
for some $i\in\{1, 2, 3\}$.
\end{lemma}
\begin{proof}
Assume by contradiction that $v^i\in \R^n_{+}$ for any $i \in\{1, 2, 3\}$. Without loss of generality we can also assume that $\|v^i\|=1$, for $i\in\{1,2,3\}$. Given that at least two eigenvalues are distinct, we have three possibilities: $\lambda_1< \lambda_2< \lambda_3$, $\lambda_1= \lambda_2<\lambda_3$ or $\lambda_1< \lambda_2=\lambda_3$. We start by analysing the possibilities $\lambda_1< \lambda_2< \lambda_3$ or $\lambda_1= \lambda_2<\lambda_3$. First we assume that either $\lambda_1< \lambda_2<0<\lambda_3$ or $\lambda_1= \lambda_2<0<\lambda_3$. Define the following auxiliary vectors
\begin{equation} \label{eq:AuxVecg}
	w^1:=v^1+t_1v^3, \qquad w^2:=v^2+t_2v^3, \qquad t_1:=\sqrt{\frac{-\lambda_1}{\lambda_3}}, \qquad t_2:=\sqrt{\frac{-\lambda_2}{\lambda_3}}.
\end{equation}
We have $\langle v^i, v^j\rangle =0$ for any $i, j\in\{1,2,3\}$ with $i\neq j$, and since
\begin{equation} \label{eq:DefVecg}
Av^1= \lambda_1 v^1, \qquad Av^2= \lambda_2 v^2, \quad Av^3= \lambda_3 v^3, \qquad v^1, v^2, v^3 \in \R^n_{+},
\end{equation}
we conclude from \eqref{eq:AuxVecg} that 
\[
	\langle Aw^1, w^1\rangle = \langle Av^1, v^1\rangle + \frac{-\lambda_1}{\lambda_3}\langle Av^3, v^3\rangle = \lambda_1\|v^1\| + \frac{-\lambda_1}{\lambda_3}\lambda_3\|v^3\| = 0, 
\]
and
\[
	\langle Aw^2, w^2\rangle = \langle Av^2, v^2\rangle + \frac{-\lambda_2}{\lambda_3}\langle Av^3, v^3\rangle = \lambda_2\|v^2\| + \frac{-\lambda_2}{\lambda_3}\lambda_3\|v^3\| = 0.
\]
Moreover, since $v^1, v^2, v^3 \in \R^n_{+}$ we conclude that $w^1, w^2 \in \left\{x\in\R^n_{+}~:~\langle Ax, x\rangle \leq 0\right\}$. On the other hand, by using \eqref{eq:DefVecg} and
\eqref{eq:AuxVecg}, we obtain that
$$
\langle A(w^1+w^2), w^1+w^2 \rangle =2\langle Aw^1, w^2\rangle = 2t_1t_2\langle A v^3, v^3\rangle = 2\frac{\sqrt{\lambda_1\lambda_2}}{\lambda_3}\lambda _3 = 2\sqrt{\lambda_1 \lambda_2} >0,
$$
hence $w^1+w^2\notin \left\{x\in\R^n_{+}~:~\langle Ax, x\rangle \leq 0\right\}$. Thus, $\left\{x\in\R^n_{+}~:~\langle Ax, x\rangle \leq 0\right\}$ is not a convex cone. 

For the general case, take $c\in \R$ such that $\lambda_2<c< \lambda_3$.
Letting $A_{c}:=A-cI_n$ we conclude that $\lambda_1-c, \lambda_2-c, \lambda_3-c$ are eigenvalues of $A_{c}$ and satisfying 
\[
	\lambda_1 -c<\lambda_2-c<0<\lambda_3-c
\]
or 
\[
	\lambda_1 -c=\lambda_2-c<0<\lambda_3-c
\]
with the three corresponding orthonormal eigenvectors $v^1, v^2, v^3 \in \R^n_{+}$. Hence, by the first part of the proof, with $A_c$ in the role of $A$, we conclude that the cone $\left\{x\in\R^n_{+}~:~\langle A_cx, x\rangle \leq 0\right\}$ is not convex. On the other hand, due to $v^1\in \R^n_{+}$ and $ \langle Av^1, v^1\rangle =\lambda_1-c<0$, we have $ \left\{x\in\R^n_{++}~:~\langle A_cx, x\rangle < 0\right\}\neq~\emptyset$. Thus, applying Corollary~\ref{cor:cor} with ${\cal K}=\R^n_{+}$ and taking into account that $\left\{x\in\R^n_{+}~:~\langle A_cx, x\rangle \leq 0\right\}$ is not convex, we conclude that $q_{A}$ is not spherically quasi-convex. 

To analyse the possibility $\lambda_1< \lambda_2=\lambda_3$, first assume that $\lambda_1<0 < \lambda_2= \lambda_3$ and define the vectors
$$
w^1:=t_1v^1+v^3, \qquad w^2:=t_2v^1+v^3, \qquad t_1=\sqrt{\frac{\lambda_2}{-\lambda_1}}, \qquad t_2=\sqrt{\frac{\lambda_3}{-\lambda_1}},
$$
and then proceed as above to obtain again a contradiction. Therefore, $v^i \notin \R^n_{+}$ for some $i\in\{1,2,3\}$.
\hfill$\square$\par \end{proof}

\section{Spherically quasi-convex quadratic functions on the subdual convex sets} \label{sec:qcqfsdcs}

In this chapter we present a  condition characterising  the  spherical   quasi-convexity of quadratic
functions on  spherically subdual convex
sets associated to subdual cones.  The results obtained  generalise the corresponding ones obtained in previous chapter (or in our published paper \cite[Section 4.1]{FerreiraNemethXiao2018}). We also summarised the results of this chapter in \cite{ferreira2019spherical}.  {\it Throughout this chapter we follow the conventions used before and assume that the cone ${\cal K}$ is a subdual ( i.e., ${\cal K}\subseteq{\cal K}^*$}) and proper cone.  A  closed set \( {\cal S} \subseteq \SP^{n-1}\) is called a {\it spherically subdual convex set} if the associated cone ${\cal K}_{\cal S}$ (defined in \eqref{eq:pccone}) is subdual. It is clear that if $A=A\tp\in\R^{n\times n}$ has only one eigenvalue, then  $q_A$ is constant and, consequently,    it is   spherically quasi-convex. Henceforth, throughout this chapter {\it we assume that $A$  has at least two distinct  eigenvalues}.  We remind that   $q_A$ and $\varphi_A$  are defined in \eqref{eq:QuadFunc} and \eqref{eq:RayleighFunction}, respectively. Two technical lemmas, which are useful in the following text, will be presented. They are generalisations of Lemma \ref{Lem:Basic} and \ref{lem:Copositive} (or, Lemmas~14~and~15~of~\cite{FerreiraNemethXiao2018}), respectively. For stating the next lemma, 
denote by $\{v^1, v^2,\dots,v^n\}$ a orthonormal system of eigenvectors of $A$  corresponding to the eigenvalues $\lambda_1 < \lambda_2\leq  \ldots \leq \lambda_n$. Given $c\in(\lambda_1,\lambda_2]$, we define the following convex cone 
\begin{equation}\label{eq-ltheta}
	{\cal L}_c:=\left\{x\in\R^n:~\lng v^1,x\rng\geq  \sqrt{ \sum_{i=2}^{n}\theta_i(c)\lng v^i,x\rng^2}\right\},  \qquad \theta_i(c):=\frac{\lambda_i-c}{c-\lambda_1}, 
\end{equation}
for $i\in\{2, \ldots, n\}$.  Note that if $\lambda_1<c<\lambda_2$, then $\theta_i(c)> 0$, for $i\in\{2, \ldots, n\}$, and both $\mathcal L_{c}$ and $-\mathcal L_{c}$ are {proper cones} (recall that a proper cone is a closed, convex, pointed cone with nonempty interior). We also need to consider the following cone
\begin{equation} \label{eq:cw}
{\mathcal W}:=(\mathcal L_{\lambda_2}\cup-\mathcal L_{\lambda_2})\cap\inte({\mathcal K}).
\end{equation} 

Considering that ${\mathcal K}$, $\mathcal L_{c}$, and $-\mathcal L_{c}$ are proper cones, as a conclusion the cone ${\mathcal W}$ is  also a proper cone, and $\inte({\cal W})\neq \emptyset$. The following lemma is a general version of Lemma \ref{Lem:Basic}:
\begin{lemma}\label{cor:Basic}
Let $n\geq 2$,     $A=A\tp\in\R^{n\times n}$ and   $\{v^1, v^2,\dots,v^n\}$ be an orthonormal system of eigenvectors of $A$  corresponding to the
eigenvalues $\lambda_1 < \lambda_2 \leq  \ldots \leq \lambda_n$, respectively.  Then, the sublevel set $[\varphi_A\leq c]$ is convex for any $c\notin(\lambda_2,\lambda_n)$ if and only if
$v^1 \in {\cal W}^*\cup-{\cal W}^*$. In particular if $v^1 \in {\cal K}^*$, then $[\varphi_A\leq c]$ is convex for any 
$c\notin(\lambda_2,\lambda_n)$.
\end{lemma}
\begin{proof} 
	By using the spectral decomposition of $A$, we have  $A=V\Lambda V\tp = \sum_{i=1}^n\lambda_iv^i(v^i)\tp$. From the definition \eqref{eq:RayleighFunction} we have 
	\begin{align}
		[\varphi_A\leq c] & = \left\{x\in\inte({\cal K}):~ \frac{\langle Ax, x\rangle}{\|x\|^2} \leq c\right\}\nonumber\\
		& = \left\{x\in\inte({\cal K}):~ \langle (A - cI_n)x, x\rangle\leq 0\right\}\nonumber\\
		& = \left\{x\in\inte({\cal K}):~ \sum_{i=1}^n(\lambda_i-c)\langle v^i, x\rangle^2 \leq 0\right\}\label{eq:bl1-1}
	\end{align}
	If $\lambda_1<c\le\lambda_2$, then by using \eqref{eq-ltheta}, the equality  \eqref{eq:bl1-1} can be completed as follows
	\begin{align}\label{eq:bl1-2} 
		{\mathcal W}  =[\varphi_A\leq\lambda_2]&\supseteq[\varphi_A\leq c] = ({\mathcal L}_c\cup-{\mathcal
		L}_c)\cap\inte({\mathcal K}) \notag \\
		&=\left\{x\in\inte(\mathcal K):~\lng v^1,x\rng^2\geq \sum_{i=2}^{n}\theta_i(c)\lng v^i,x\rng^2\right\}.
	\end{align}
\noindent			
 \textbf{\it Sufficiency of the first statement:} 
		
	Let $v^1\in {\cal W}^*$ (a similar argument holds for $v^1\in-{\cal W}^*$).  
	
	If $c<\lambda_1$, then considering that  $v^1, v^2, \ldots , v^n$ are linearly independent and $0\notin \inte ({\cal K})$, we obtain from \eqref{eq:bl1-1} that $\sum_{i=1}^n(\lambda_i-c)\langle v^i, x\rangle^2 > 0$ for any $x\in \inte({\cal K})$ and hence $[\varphi_A\leq c]=\emptyset $ is convex.  
	
	If $c=\lambda_1$, then  \eqref{eq:bl1-1} implies that $[\varphi_A\leq c]={\cal S}\cap \inte ({\cal K}),$ where 
	\begin{align*}
		{\cal S}& := \left\{x\in\R^n:~ \sum_{i=2}^n(\lambda_i-c)\langle v^i, x\rangle^2 = 0\right\} \\
		&=\lf\{x\in\R^n ~: \langle v^i, x\rangle=0,~ i\in \{2,\dots, n\}\rg\}.
	\end{align*}
	Thus, due to $\inte ({\cal K})$ and ${\cal S}$ being convex, we conclude that $[\varphi_A\leq c]$ is also convex. 
	
	If $\lambda_1<c\le\lambda_2$. Since $v^1\in {\cal W}^*$, for any $x\in\mathcal W$ we obtain that $\lng v^1, x\rng\ge0$ and from \eqref{eq:bl1-2} we have $-{\cal L}_c\cap \inte ({\cal K})=\emptyset$ and hence $[\varphi_A\leq c]= {\cal L}_c\cap \inte ({\cal K})$. Due to the convexity of the cones ${\cal L}_c$ and $\inte ({\cal K})$, we obtain  that  $[\varphi_A\leq c]$ is convex.
	
	Finally, if $c\ge \lambda_n$, then \eqref{eq:bl1-1} implies that $\sum_{i=1}^n(\lambda_i-c)\langle v^i, x\rangle^2 \leq 0$ for any $x\in\inte({\cal K})$ and hence $[\varphi_A\leq c]=\inte ({\cal K})$ is convex.\\
\noindent	
\textbf{\it Necessity of the first statement:}  

	We will show that  $v^1\notin {\cal W}^*\cup-{\cal W}^*$ implies $[\varphi_A\leq c]$ is not convex, for some $c\in (\lambda_1,\lambda_2)$.	Suppose that $v^1\notin {\cal W}^*\cup-{\cal W}^*$.  Thus, considering that $\inte(\mathcal W)\neq \emptyset$, there exist $y,z\in\inte(\mathcal W)$ such that $\lng v^1,y\rng>0$ and $\lng v^1,z\rng<0$. Thus, \eqref{eq-ltheta} and \eqref{eq:cw} implies that 
	\begin{equation} \label{eq:yz}
		y\in\inte(\mathcal K)\cap\inte(\mathcal L_{\lambda_2}), \qquad \quad   z\in\inte(\mathcal K)\cap\inte(-\mathcal L_{\lambda_2}).
	\end{equation}
	We claim that there exists a ${\bar c}\in (\lambda_1,\lambda_2)$ such that $y\in\inte(\mathcal K)\cap \inte({\mathcal L}_{{\bar c}})$ and  $z\in\inte(\mathcal K)\cap \inte(-\mathcal L_{{\bar c}})$.  In order to simplify the notations, for  $ x\in\R^n$ and   $c\in (\lambda_1,\lambda_2]$, we define the following function   
	\begin{equation} \label{eq:fpsy}
		\xi(x,c):=\sqrt{\sum_{i=2}^n\theta_i(c)\lng v^i,x\rng^2}.
	\end{equation}
	Note that $\xi$ is  a  continuous function  and,  from   the definition of $\theta_i$ in   \eqref{eq-ltheta},  it is also
			decreasing with respect to the second variable $c$. By using \eqref{eq-ltheta} and \eqref{eq:fpsy} we have 
	\begin{equation} \label{eq:edef}
		\inte(\mathcal K)\cap \inte({\cal L}_c)=\left\{x\in\inte{\mathcal K}:~\lng v^1,x\rng>\xi(x,c)\right\}, \qquad \forall ~c\in (\lambda_1,\lambda_2].
	\end{equation}		
	Thus,  taking into account  the  first inclusion in \eqref{eq:yz}  we conclude,  by setting   $c=\lambda_2$ in   \eqref{eq:edef},    that 
			\[
			\lim_{c\to\lambda_2}\xi(y,c)=\xi(y,\lambda_2)<\lng v^1,y\rng. 
			\] 
			Hence,  there exists a ${\hat c}\in(\lambda_1,\lambda_2)$ sufficiently close to $\lambda_2$ such that $\xi(y,{\hat c})< \lng v^1,y\rng$. Similarly, we can also prove that   there exists  a ${\tilde c}\in(\lambda_1,\lambda_2)$ sufficiently close to $\lambda_2$ such that $\xi(z,{\tilde c})< -\lng v^1,z\rng$. Thus, letting ${\bar c}=\max\{{\hat c}, {\tilde c} \}$ we conclude that $\xi(y,{\bar c})< \lng v^1,y\rng$ and  $\xi(z,{\bar c})< -\lng v^1,z\rng$, which by   \eqref{eq:fpsy} and \eqref{eq:edef} yields
			\begin{equation} \label{eq:icba}
			y\in \inte({\mathcal L}_{{\bar c}}), \qquad \quad z\in\inte(-\mathcal L_{{\bar c}}). 
			\end{equation} 
We know by 	 \eqref{eq:yz} that 	$y\in\inte(\mathcal K)$ and  $z\in\inte(\mathcal K)$, which together  with \eqref{eq:icba} yields
$y\in\inte(\mathcal K)\cap \inte({\mathcal L}_{{\bar c}})$ and  $z\in\inte(\mathcal K)\cap \inte(-\mathcal L_{{\bar c}})$ and the  claim is
concluded.  Therefore,  there exist  $r_y>0$ and  $r_z>0$ such that $\textbf{B}(y,r_y)\subset \inte(\mathcal K)\cap \inte({\mathcal L}_{{\bar c}})$ and
$\textbf{B}(z,r_z)\subset \inte(\mathcal K)\cap \inte(-\mathcal L_{{\bar c}})$, where   $\textbf{B}(y,r_y)$ and  $\textbf{B}(z,r_z)$ denote  the open balls  with  centers $y$, $z$ and radius $r_y>0$, $r_z>0$, respectively. Hence, by  dimensionality reasons, we can take  $u_y\in \inte(\mathcal K)\cap \inte({\mathcal L}_{{\bar
c}})$ and  $u_z \in \inte(\mathcal K)\cap \inte(-\mathcal L_{{\bar c}})$ such that $v^1$,  $u_y$ and $u_z$ are \gls{li}.  Thus,
in particular, we have  $0\notin [u_y,u_z]$,  where  $[u_y,u_z]$ denotes the straight line segment joining $u_y$ to $u_z$. Since $ \inte({\mathcal
L}_{{\bar c}})\cap \inte(-\mathcal L_{{\bar c}})=~\emptyset$ and $0\notin [u_y,u_z]$, the  segment $[u_y,u_z]$ is intersecting, at the distinct
points
$w_y\neq 0$ and $w_z\neq 0$,   the boundaries of the sets  $ \inte({\mathcal L}_{{\bar c}})$ and $  \inte(-\mathcal L_{{\bar c}})$, respectively.
Moreover, due to   $u_y$ and $u_z$  being  l.i.,  $0\notin [u_y,u_z]$ and    $w_y, w_z\in  [u_y,u_z]$, we conclude that the vectors $v^1$,  $w_y$
and $w_z$ are  also l.i.. Our next task is to  prove  that
\begin{equation} \label{eq:mwywz}
\frac{1}{2}(w_y+w_z)\notin  {\mathcal L_{\bar c}}\cup-{\mathcal L_{\bar c}}.
\end{equation}
First, due to $w_y$ and $ w_z$ belonging  to  the boundaries  of  $ {\mathcal L}_{{\bar c}}$  and $  -\mathcal L_{{\bar c}}$, respectively,  we obtain from \eqref{eq-ltheta}  that
\begin{equation} \label{eq:ewywz}
\lng v^1,w_y\rng=  \sqrt{\sum_{i=2}^n\theta_i({\bar c})\lng v^i,w_y\rng^2}, \qquad \lng v^1,w_z\rng= - \sqrt{\sum_{i=2}^n\theta_i({\bar c})\lng v^i,w_z\rng^2}.
\end{equation}
On the other hand, by using the two equalities in  \eqref{eq:ewywz},  we obtain after some algebraic manipulations that 
\footnotesize
\begin{align*}
	\sum_{i=2}^n\theta_i({\bar c})\left\lng v^i,\frac{1}{2}(w_y+w_z) \right\rng^2 &= \sum_{i=2}^n\theta_i({\bar c}) \lf(\lf\lng v^i,\frac{1}{2}w_y\rg\rng + \lf\lng v^i,\frac{1}{2}w_z\rg\rng \rg)^2 \\
	& = \left\lng v^1,\frac{1}{2}w_y\right\rng^2+ \left\lng v^1,\frac{1}{2}w_z \right\rng^2+ 2\sum_{i=2}^n\theta_i({\bar c})\left\lng v^i,\frac{1}{2}w_y\right\rng\left\lng v^i,\frac{1}{2}w_z\right\rng.
\end{align*}
\normalsize
Thus, considering that  
{\footnotesize 
\[
	\left\lng v^1,\frac{1}{2}(w_y+w_z) \right\rng^2= \left\lng v^1,\frac{1}{2}w_y\right\rng^2+ \left\lng v^1,\frac{1}{2}w_z \right\rng^2+ 2 \left\lng v^1,\frac{1}{2}w_y\right\rng\left\lng v^1,\frac{1}{2}w_z \right\rng,
\]
}
we have 
{\footnotesize
\begin{align} \label{eq:fil3}
\sum_{i=2}^n\theta_i({\bar c})\left\lng v^i,\frac{1}{2}(w_y+w_z) \right\rng^2=\left\lng v^1,\frac{1}{2}(w_y+w_z) \right\rng^2 &- 2 \left\lng v^1,\frac{1}{2}w_y\right\rng\left\lng v^1,\frac{1}{2}w_z \right\rng \nonumber\\
	&+ 2\sum_{i=2}^n\theta_i({\bar c})\left\lng v^i,\frac{1}{2}w_y\right\rng\left\lng v^i,\frac{1}{2}w_z\right\rng.
\end{align}
}
Applying   Cauchy-Schwarz   inequality and then,  using again both equalities  in  \eqref{eq:ewywz}, we conclude that  
\begin{align} \label{eq:esil3}
-\sum_{i=2}^n\theta_i({\bar c})\left\lng v^i,\frac{1}{2}w_y\right\rng\left\lng v^i,\frac{1}{2}w_z\right\rng &\leq    \sqrt{\sum_{i=2}^n\theta_i({\bar c})\lng v^i,w_y\rng^2}   \sqrt{\sum_{i=2}^n\theta_i({\bar c})\lng v^i,w_z\rng^2} \\                                                                                                                                                                                                  &=-\left\lng v^1,\frac{1}{2}w_y\right\rng\left\lng v^1,\frac{1}{2}w_z \right\rng \notag. 
\end{align}

We are going to prove that the inequality \eqref{eq:esil3} is  strict. For that, assume the contrary, i.e., that the last inequality holds as equality. In this case,  there exists  $\alpha \neq 0$ such that 
\begin{align*}
	\Big(\sqrt{\theta_2({\bar c})}\left\lng v^2,\frac{1}{2}w_y\right\rng,& \dots, \sqrt{\theta_n({\bar c})}\left\lng v^n,\frac{1}{2}w_y\right\rng \Big) \\
	& = \alpha \Big(\sqrt{\theta_2({\bar c})}\left\lng v^2,-\frac{1}{2}w_z\right\rng, \dots, \sqrt{\theta_n({\bar c})}\left\lng v^n,-\frac{1}{2}w_z\right\rng \Big),
\end{align*}

which implies  that $w_y+\alpha w_z$ is orthogonal to the set of vectors $\{v^2, \ldots, v^n \}$. Thus, since  the set $\{v^1, v^2,\dots,v^n\}$ is an orthonormal system,  $w_y+\alpha w_z$ is parallel to the vector $v^1$, which is absurd due to vectors $v^1$,  $w_y$ and  $w_z$ being   l.i.. Hence, \eqref{eq:esil3} holds strictly and combining it with \eqref{eq:fil3} we conclude that
$$
\sum_{i=2}^n\theta_i({\bar c})\left\lng v^i,\frac{1}{2}(w_y+w_z) \right\rng^2> \left\lng v^1,\frac{1}{2}(w_y+w_z) \right\rng^2, 
$$
and  \eqref{eq:mwywz} holds.  Therefore,  considering that $\frac{1}{2}(w_y+w_z) \in  (u_y,u_z)$, we conclude that   $(u_y,u_z) \not\subset {\mathcal L_{\bar c}}\cup-{\mathcal L_{\bar c}}$. Thus,  using notation \eqref{eq:bl1-2},  we also have   $(u_y,u_z) \not\subset ({\mathcal L}_{\bar c}\cup-{\mathcal L}_{\bar c})\cap\inte({\mathcal K})=[\varphi_A\le {\bar c}]$, and due to $u_y, u_z\in ({\mathcal L}_{\bar c}\cup-{\mathcal L}_{\bar c})\cap\inte({\mathcal K})=[\varphi_A\le {\bar c}]$, it follows that $[\varphi_A\le {\bar c}]$ is not 
			convex.\\
\noindent		
 \textbf{\it Proof of second statement:} It follows  from $\mathcal K^*\subseteq\mathcal W^*$.  
\hfill$\square$\par\end{proof}

\begin{remark}
	It is easy to check that for any two cones ${\cal A}\subset \R^n$ and ${\cal B}\subset \R^n$, we have $\lf({\cal A}\cap {\cal B}\rg) = \lf({\cal A}^* + {\cal B}^*\rg)$. The dual of $\mathcal W$ in \eqref{eq:cw}can be expressed as 
	\begin{eqnarray}\label{eq-o-w}
		\begin{array}{rcl}
			\mathcal W^*=[(\mathcal K\cap\mathcal L_{\lambda_2})\cup (\mathcal K\cap-\mathcal
			L_{\lambda_2})]^*
			&=&(\mathcal K\cap\mathcal L_{\lambda_2})^*\cap (\mathcal K\cap-\mathcal L_{\lambda_2})^*\\
			&=&(\mathcal K^*+\mathcal L_{\lambda_2}^*)\cap(\mathcal K^*-\mathcal L_{\lambda_2}^*).
		\end{array}
	\end{eqnarray} 
\end{remark}

\begin{corollary}\label{cor-conv-sl}
	Suppose that $n\ge 3$ and $\lambda_2\le
	(\lambda_1+\lambda_3)/2$. If either $\mathcal K\cap -\mathcal L_{\lambda_2}=\{0\}$ or $\mathcal
	K\cap \mathcal L_{\lambda_2}=\{0\}$, then $[\varphi_A\le c]$ is convex for any $c\notin (\lambda_2,\lambda_n)$.
\end{corollary}

\begin{proof}
	First note that if $n\ge 3$ and $\lambda_2\le (\lambda_1+\lambda_3)/2$, then $\theta_i(\lambda_2)\ge 1$ for any $i\ge 3$.  Define  the cone 
	\[
		\mathcal L_{[v^2]^\perp}:=\lf\{x\in\R^n:~\lng v^1,x\rng\ge\sqrt{\sum_{i=3}^n\lng v^i,x\rng^2}\rg\}. 
	\]
 Note that $ \mathcal L_{[v^2]^\perp}$ is a self-dual Lorentz cone  as a subset of the subspace $[v^2]^\perp$.   Moreover,   considering that  $\theta_i(\lambda_2)\ge 1$ for any $i\ge 3$, we conclude 
	\[
		\mathcal L_{\lambda_2}\cap [v^2]^\perp = \lf\{x\in\R^n:~\lng v^1,x\rng\ge\sqrt{\sum_{i=3}^n\theta(\lambda_2) \lng v^i,x\rng^2}\rg\}\subseteq \mathcal L_{[v^2]^\perp}.
	\]
Consequently,  taking into account  that  $ \mathcal L_{[v^2]^\perp}$ is a self-dual cone,  the cone  $\mathcal L_{\lambda_2}\cap [v^2]^\perp$ is subdual as a subset of the subspace $[v^2]^\perp$. To simplify the notation, denote by upper star (i.e., $^*$)
	the dual of a cone in $\mathbb R^n$ and by lower star (i.e., $_*$) the dual of a cone in $[v^2]^\perp$. Thus, using this notation we will prove
	\begin{equation} \label{eq:ssls}
	\mathcal
	L_{\lambda_2}^*=(\mathcal L_{\lambda_2}\cap [v^2]^\perp)_*
	\end{equation}
	 Indeed, since $v^2,-v^2\in\mathcal L_{\lambda_2}$, for 
	any $z\in\mathcal L_{\lambda_2}^*$, we have $\lng z,v^2\rng=0$ and hence $\mathcal L_{\lambda_2}^*\subseteq [v^2]^\perp$, which implies
	$\mathcal L_{\lambda_2}^*\subseteq (\mathcal L_{\lambda_2}\cap [v^2]^\perp)_*$. 
	
	Conversely,  arbitrarily take $u\in (\mathcal L_{\lambda_2}\cap [v^2]^\perp)_*$, and take $w\in\mathcal L_{\lambda_2}\cap [v^2]^\perp$ then for any $t\in\R$ we have $v=w+tv^2\in {\cal L}_{\lambda_2}$. Hence,   $\lng u,v\rng=\lng u,w\rng\ge 0$, which implies that $u\in \mathcal L_{\lambda_2}^*$. Hence,  we conclude that 
	$(\mathcal L_{\lambda_2}\cap [v^2]^\perp)_*\subseteq\mathcal L_{\lambda_2}^*$, and \eqref{eq:ssls} is proved.  Suppose $\mathcal K\cap -\mathcal L_{\lambda_2}=\{0\}$,  by using the first equality in \eqref{eq-o-w} we obtain $\mathcal W^*=(\mathcal K\cap\mathcal L_{\lambda_2})^*$.  Therefore, considering that  $\mathcal L_{\lambda_2}\cap [v^2]^\perp$ is subdual  and \eqref{eq:ssls}, we obtain 
	 \[v^1\in\mathcal L_{\lambda_2}\cap[v^2]^\perp\subseteq (\mathcal L_{\lambda_2}\cap[v^2]^\perp)_*=\mathcal L_{\lambda_2}^*\subseteq
	(\mathcal K\cap\mathcal L_{\lambda_2})^*=\mathcal W^*.\] Hence, following Lemma \ref{cor:Basic} we conclude that $[\varphi_A\le c]$ is convex for any $c\notin (\lambda_2,\lambda_n)$.  The case $\mathcal K\cap {\cal L}_{\lambda_2}=\{0\}$ can be proved similarly. 
\hfill$\square$\par\end{proof}
\begin{lemma}\label{lem:locnonconv}
	Let $n\ge 3$ and $B=B\tp\in\R^{n\times n}$. Let $\mu_1 \leq \mu_2 \leq  \ldots \leq \mu_n$ be  eigenvalues of the matrix $B$.  Assume   that one of the  following two conditions holds:
	\begin{enumerate}
		\item[(a)] $\mu_1=\mu_ 2<0  < \mu_n$;
		\item[(b)]  $ \mu_1<\mu_ 2<0<\mu_n$. 
	\end{enumerate}
	 Then, for any ${\bar x}\in\R^n\setminus\{0\}$ such that $B{\bar x}\ne 0$ and  $\lf\lng B{\bar x},{\bar x}\rg\rng=0$,  and any number  $\delta>0$, the set 
	$
	\Xi\lf(B,{\bar x},\delta\rg):=\lf\{x\in\R^n:\lf\|x-{\bar x}\rg\|\le\delta,\textrm{ }\lng
	Bx,x\rng\le0\rg\}
	$
	is not convex.
\end{lemma}
\begin{proof}
Since  $\mu_1=\min_{x\in \SP^{n-1}} q_B(x) <\max_{x\in \SP^{n-1}} q_B(x)= \mu_ n$, we can take   ${\bar x}\in\R^n\setminus\{0\}$  such that $B{\bar x}\ne 0$ and  $\lf\lng B{\bar x},{\bar x}\rg\rng=0$. Define the following vector subspace of $\R^{n}$:
\[
	{\mc N}:=[\{u\in \R^n:~Bu=\mu  u,\textrm{ for some }\mu<0\}].
\]

It follows from assumption (a) or (b)  that $\dim({\mc N})\ge2$.   For simplifying  the notation we set 
\begin{equation}\label{eq:bar_y} 
	{\bar y}:=B{\bar x}\neq 0.
\end{equation}

To proceed with the proof, we first need to prove that ${\mc N}\neq [\bar y]^\perp$.  Assume to the contrary that ${\mc N}= [\bar y]^\perp$. In this case, due to \eqref{eq:bar_y} and $B=B\tp$, the definition  of  $[\bar y]^\perp$ implies that 
\begin{equation}\label{eq:sub_sp_y}
	\lf\lng Bv, {\bar x}\rg\rng=\lf\lng B{\bar x}, v\rg\rng=\lf\lng {\bar y}, v\rg\rng=0, \quad \forall v\in {\mc N}.
\end{equation}

Thus, it follows from the definition  of ${\mc N}$ that   $ \lf\lng Bv, {\bar x}\rg\rng=\lf\lng v, {\bar x}\rg\rng=0$, for any $v\in {\mc N}$,  which implies   ${\mc N}\subset [\bar x]\pr:=\{v\in \R^{n}:~  \lf\lng v,{\bar x}\rg\rng=0\}$. Moreover,  considering that $\lf\lng {\bar y},{\bar x}\rg\rng=0$,  we also have ${\bar y} \in [\bar x]\pr$. Hence,  we conclude that $[\bar{y}]+\mathcal{N}\subset [\bar x]\pr$. Since by definition \eqref{eq:bar_y} we have${\bar y}\neq 0$, then we conclude  ${\bar y}\notin [\bar{y}]\pr = {\mc N}$. Due to ${\bar y}\neq0$ and  ${\mc N}= [\bar y]^\perp$ we have $\dim ([\bar{y}]+\mathcal{N})=n$. Combining with the fact that $[\bar{y}]+\mathcal{N}\subset [\bar x]\pr$, we obtain  ${\bar x} = 0$, which contradicts the assumption ${\bar x} \neq 0$. Therefore, ${\mc N}\neq [\bar y]^\perp$. Thus, we have  
	\[
		\dim(\mc N\cap [{\bar y}]^\perp)\geq \dim{\mc N}+ \dim{ [{\bar y}]^\perp} -  \dim\R^n \ge 2+(n-1)-n=1.
	\]
Hence, there exist a unit vector $a\in\mc N\cap [{\bar y}]\pr$, so that $\lng a,{\bar y}\rng=0$.   Since $\mathcal{N}\ne [\bar{y}]^\perp$, we can choose a sequence
of vectors  $\{a^n\}\subset \mc N$  such that $\lim_{n\to\infty}a^n=a$ and  $\lng a^n,{\bar y}\rng\ne 0$.  Let $\{u^1, u^2,\dots,u^n\}$ be an orthonormal system of eigenvectors of $B$  corresponding to the eigenvalues $\mu_1 , \mu_2 ,  \ldots , \mu_n$, respectively. Note that the spectral decomposition of $B$ implies
$
B=\sum_{i=1}^n\mu_i u^i(u^i)\tp.
$
Since  $\{a^n\}\subset \mc N$, we can write $a^n=\sum_{i=1}^\ell\alpha_{n, i}  u^i$, where $2\leq \ell=\dim(\mc N)<n$ and $\mu_1,  \ldots , \mu_\ell$ are the negative eigenvalues  of $B$. Thus, 
\[
	\lng Ba^n,a^n\rng=\sum_{i=1}^\ell\sum_{j=1}^\ell\lng B\alpha_{n,i}u^i,\alpha_{n,j}u^j\rng =\sum_{i=1}^\ell\alpha^2_{n, i}(u^i)\tp Bu^i=\sum_{i=1}^\ell\alpha^2_{n, i}\mu_i<0.
\]
For proceeding with the proof,  we define 
\[
 p^n:={\bar x}+t_na^n, \qquad t_n:=-2\f{\lng a^n,{\bar y}\rng}{\lng Ba^n,a^n\rng}.
\] 
Then,
	$\lng Bp^n,p^n\rng=0$ and,  due to $\lng a,{\bar y}\rng=0$ and $\lim_{n\to\infty}a^n=a$, we have  $\lim_{n\to\infty}p^n={\bar x}$. Hence, if $n$ is sufficiently large, then for any  $\delta>0$ arbitrary but fixed, we have  $p^n\in\Xi\lf(B,{\bar x},\delta\rg)$. 
	For such an $n$,  after some simple algebraic manipulations we conclude 
	\[\lf\lng B\lf(\f{{\bar x}+p^n}2\rg),\f{{\bar x}+p^n}2\rg\rng=-\f{\lf\lng a^n,{\bar y}\rg\rng^2}{\lng Ba^n,a^n\rng}>0.
	\]
	Hence, ${\bar x},p^n\in\Xi\lf(B,{\bar x},\delta\rg)$, but $({\bar x}+p^n)/2\notin\Xi\lf(B,{\bar x},\delta\rg)$. Therefore, $\Xi\lf(B,{\bar x},\delta\rg)$ is not
	convex.   
\hfill$\square$\par\end{proof}
\begin{proposition}\label{pr:multwo}
	Let $n\ge 3$ and $A=A\tp\in\R^{n\times n}$  is a nonsingular matrix.   Suppose that $q_A$  is    not constant  and   $\lambda_1 \leq \lambda_2 \leq  \ldots \leq \lambda_n$ are eigenvalues of $A$.  If $q_A$  is  quasi-convex, then the following conditions hold:
\begin{enumerate}
		\item[(i)] $\lambda_1 < \lambda_2$;
		\item[(ii)] either   $\lambda_2\leq \min_{x\in {\bar {\cal C}}} q_A(x)$  or  $\max_{x\in {\bar {\cal C}}} q_A(x)\leq\lambda_2$. 
	\end{enumerate}	
\end{proposition}
\begin{proof}
 Suppose by contradiction that one of the  following two conditions holds:
	\begin{enumerate}
		\item[(a)] $\lambda_1= \lambda_2$;
		\item[(b)]  $\min_{x\in {\bar {\cal C}}} q_A(x) <\lambda_2<\max_{x\in {\bar {\cal C}}} q_A(x)$. 
	\end{enumerate}
First of all, note that due to  $q_A$ not being constant, we have   $\lambda_1\leq\min_{x\in {\bar {\cal C}}} q_A(x) <\max_{x\in {\bar {\cal C}}} q_A(x)\leq \lambda_ n$, where ${\bar {\cal C}}$ is defined in \eqref{eq:ccb}. If the   condition (a) holds, we can  take a scalar $\mu \in \R$ such that $\mu \neq \lambda_i$ for any $i\in\{1, \dots, n\}$ and  satisfying 
\begin{equation} \label{eq:cdfff}
 \lambda_1=\lambda_ 2\leq\min_{x\in {\bar {\cal C}}} q_A(x)< \mu  <\max_{x\in {\bar {\cal C}}} q_A(x)\leq \lambda_n. 
\end{equation}
Otherwise,  if the condition (b) holds, we take $\mu \in \R$  satisfying 
\begin{equation} \label{eq:cdsss}
\lambda_1\leq\min_{x\in {\bar {\cal C}}} q_A(x)<\lambda_ 2<\mu< \max_{x\in {\bar {\cal C}}} q_A(x)\leq \lambda_n.
\end{equation}
	Then, either the conditions \eqref{eq:cdfff} or \eqref{eq:cdsss} implies that  $\pm(A-\mu I_n)$ is not $\mc K$-copositive. Since the matrix $A-\mu I_n$ is not $\mc K$-copositive, we can find a point $p\in\mc K$ such that $\lng Ap,p\rng<\mu\|p\|^2$. Hence, we can find $u\in\inte(\mc K)$ sufficiently close to $p$ such that $\lng Au,u\rng<\mu\|u\|^2$. 
	
	Similarly, since $-(A-\mu I_n)=\mu I_n-A$ is not $\mc K$-copositive, we can find $v\in\inte(\mc K)$ such that $\lng Av,v\rng>\mu\|v\|^2$. Therefore, we take $t\in (0,1)$ and define
	\[
		\inte({\cal K})\ni {\bar x}:=(1-t)u+tv.
	\]

By continuity, we have $\lng A{\bar x},{\bar x}\rng=\mu\|{\bar x}\|^2$.  Denoting $B=A-\mu I$, the eigenvalues of $B$ are given by $\mu_i:=\lambda_i-\mu$, for $i\in\{1, 2, \dots, n\}$. Thus,  we conclude from \eqref{eq:cdfff} and \eqref{eq:cdsss} that either
\begin{equation} \label{eq:cdfffa}
 \mu_1=\mu_ 2<0 < \mu_n,   \quad \mbox{or} \qquad \mu_1<\mu_ 2<0<\mu_n,  
\end{equation}
if  either the condition (a) or (b) holds, respectively.  Considering that  $B{\bar x}\ne 0$ and  $\lng B{\bar x},{\bar x}\rng=0$, we conclude from Lemma~\ref{lem:locnonconv} that,  for any $\delta>0$,  the set 
	\[\Xi\lf(B,{\bar x},\delta\rg):=\lf\{x\in\R^n:\lf\|x-{\bar x}\rg\|\le\delta,\textrm{ }\lng Bx,x\rng\le0\rg\}, \] 
	  is not convex. Hence, there exists an $s\in (0,1)$ and
	$a^0,a^1\in\Xi\lf(B,{\bar x},\delta\rg)$ such that $a^s:=(1-s)a^0+sa^1\notin\Xi\lf(B,{\bar x},\delta\rg)$. Thus, since the closed ball centered at ${\bar x}$ and radius $\delta$ is convex, $a^s\notin \Xi\lf(B,{\bar x},\delta\rg)$ implies 
	$\lng Aa^s,a^s\rng-\mu\|a^s\|^2=\lng Ba^s,a^s\rng>0$. On the other hand, since $a^0,a^1\in\Xi\lf(B,{\bar x},\delta\rg)$, we
	have $\lng Aa^i,a^i\rng-\mu\|a^i\|^2=\lng Ba^i,a^i\rng\le0$, for $i\in\{0,1\}$. Furthermore, if $\delta$ is sufficiently
	small, then since ${\bar x}\in\inte(\mc K)$, we have $a^0,a^1\in\inte\mc K$. Hence,	$a^0,a^1\in[\varphi_A\le\mu]$ and
	$a^s\notin[\varphi_A\le\mu]$. By using Corollary~\ref{cor:cor}, this contradicts the spherical quasi-convexity of $A$.  
 \hfill$\square$\par\end{proof}

The proof of following lemma is based on Lemma \ref{lem:Copositive}. 

\begin{lemma} \label{lem:k_Copositive}
 Let  $A \in\R^{n\times n}$ and $\lambda,c\in\R$ such that $\lambda\leq c$.  If  $\lambda  I_n-A$ is ${\cal K}$-copositive,   then  $[\varphi_A\leq c]=\inte ({\cal K})$. As a consequence, the set $[\varphi_A\leq c]$ is convex.
\end{lemma}

\begin{proof}
Let $c\in \R$ and $[\varphi_A\leq c]=\{x\in \inte({\cal K}):~ \langle A x,x \rangle-c \|x\|^2\leq 0\}$. Suppose that $\lambda\leq c$, for any $x\in \inte({\cal K})$ we have 
\[
	\langle A x,x \rangle-c\|x\|^2\leq \langle A x,x \rangle-\lambda \|x\|^2= \langle (A-\lambda I_n) x,x \rangle,
\]
and considering that $\lambda I_n-A$ is ${\cal K}$-copositive, that is
\[
	\langle A x,x \rangle-c\|x\|^2\leq\langle (A-\lambda I_n) x,x \rangle\leq 0,\quad \forall x\in \inte({\cal K}),
\]
hence $\langle A x,x \rangle-c\|x\|^2\leq 0$ holds for any $x\in {\cal K}$, which implies that $[\varphi_A\leq c]=\inte ({\cal K})$.
\hfill$\square$\par \end{proof}

The following theorem combines the results from Proposition \ref{pr:multwo}, Lemma \ref{cor:Basic},  Lemma \ref{lem:k_Copositive} and Corollary \ref{cor:cor}:
\begin{theorem}\label{th:best}
Let $n\geq 3$, $k\ge1$, $A=A\tp\in\R^{n\times n}$ and  $\{v^1, v^2,\dots,v^n\}$ be an orthonormal system of eigenvectors of $A$  corresponding to
the eigenvalues $\lambda_1=\dots=\lambda_j<\lambda_{k+1}\leq \ldots \leq \lambda_n$, respectively.  Then, we have the following
statements:
\begin{enumerate}
	\item[(i)] If   $q_A$  is  quasi-convex and not constant, then $k=1$. 
	\item[(ii)]  If   $q_A$  is  quasi-convex and not constant, then either   $\lambda_2\leq \min_{x\in {\bar {\cal C}}} q_A(x)$  or  $\max_{x\in {\bar {\cal C}}} q_A(x)\leq\lambda_2$. 
	\item[(iii)] Suppose that $k=1$ and $\lambda_2I_n-A$ is $\mc K$-copositive. Then, $q_A$ is spherically quasi-convex if and only if 
		$v^1 \in {\cal W}^*\cup-{\cal W}^*$. In particular if $v^1 \in {\cal K}^*$, then $q_A$ is spherically quasi-convex.
\end{enumerate}
\end{theorem}

\begin{proof}
	Items (i) and (ii) follow from Proposition \ref{pr:multwo}. Item (iii) follows from Lemma \ref{cor:Basic},  Lemma \ref{lem:k_Copositive} and Corollary~\ref{cor:cor}.  
 \hfill$\square$\par\end{proof}
The next corollary follows  by  combining  Lemma \ref{lem:k_Copositive} and Corollary~\ref{cor-conv-sl}. 
\begin{corollary}\label{cor-best}
Let $n\geq 3$, $A=A\tp\in\R^{n\times n}$ and   $\lambda_1<\lambda_{2}\leq \ldots \leq \lambda_n$  the eigenvalues of $A$.   Suppose that $\lambda_2\le
	(\lambda_1+\lambda_3)/2$ and $\lambda_2I_n-A$ is $\mc K$-copositive. If either $\mathcal K\cap -\mathcal L_{\lambda_2}=\{0\}$ or $\mathcal 
	K\cap \mathcal L_{\lambda_2}=\{0\}$, then $q_A$ is spherically quasi-convex.
\end{corollary}
\begin{proof}
Let $n\geq 3$, $A=A\tp\in\R^{n\times n}$ and   $\lambda_1<\lambda_{2}\leq \ldots \leq \lambda_n$  the eigenvalues of $A$. 
If we have that $\lambda_2\le (\lambda_1+\lambda_3)/2$ and suppose either $\mathcal K\cap -\mathcal L_{\lambda_2}=\{0\}$ or $\mathcal  K\cap \mathcal L_{\lambda_2}=\{0\}$ by Corollary \ref{cor-conv-sl} we obtain that $[\varphi_A\leq c]$ is convex for any $c\notin(\lambda_2, \lambda_n)$. 

Suppose that $\lambda_2I_n-A$ is $\mc K$-copositive, by Lemma \ref{lem:k_Copositive} we have $[\varphi_A\leq c] = \inte({\cal K})$ is convex for $c\ge \lambda_2$. Therefore, by Corollary \ref{cor:cor} we conclude that $q_A$ is spherically quasi-convex.
\hfill$\square$\par\end{proof}

In the following two theorems we present classes of quadratic quasi-convex functions defined  in
spherically subdual convex sets, which
include as particular instances in Example~\ref{Ex:countergg} and \ref{Ex:CountExMulEngfm}.
\begin{theorem} \label{th:countergg}
Let $n\geq  3$,     $A=A\tp\in\R^{n\times n}$ and $\{v^1, v^2,\dots,v^n\}$ be an orthonormal system of eigenvectors of $A$  corresponding to the eigenvalues $\lambda_1 , \lambda_2 ,  \ldots , \lambda_n$ , respectively.    Assume that   $\lambda:=\lambda_1$, $\mu:=\lambda_2=\ldots= \lambda_{n-1}$,  $\eta:=\lambda_n$ and
\begin{equation} \label{eq:ciqc-2}
 v^1-\sqrt{\frac{\eta-\mu}{\mu-\lambda}}|v^{n}|^{\cal K}\in {\cal K}^*,   \qquad  \lambda<\mu<\eta, 
\end{equation}
where $ |\cdot |^{\cal K}$ is defined in  \eqref{eq:npav}.	 Then,  the quadratic function    $q_A$ is spherically quasi-convex.
\end{theorem}
\begin{proof}
By using the spectral decomposition of $A$, we have
\begin{equation} \label{eq:equivA}
A=\sum_{i=1}^n\lambda_i v^i(v^i)\tp=\lambda v^1(v^1)\tp+ \mu \sum_{j=2}^{n-1} v^j(v^j)\tp+\eta v^n(v^n)\tp.
\end{equation}
Hence,  for any   $x\in{\cal K}$,  by using $\|x\|^2=\sum_{i=1}^n\lng v^i,x\rng^2$ and \eqref{eq:equivA}, we obtain
\begin{align}
	\langle A x,x \rangle-\mu\|x\|^2 &= \lf(\lambda \lng v^1,x\rng^2 + \mu\sum_{i=2}^{n-1}\lng v^i,x\rng^2 + \eta \lng v^n,x\rng^2 \rg) - \mu\sum_{i=1}^n\lng v^i,x\rng^2 \nonumber\\
	& =(\mu-\lambda) \left[\frac{\eta-\mu}{\mu-\lambda}\langle v^n, x\rangle^2-\langle v^1, x\rangle^2 \right]. \label{eq:ceq-2}
\end{align}
From \eqref{eq:ciqc-2} we conclude that 
\begin{equation}\label{eq:proj_com}
	0 \leq \lng v^1-\sqrt{\frac{\eta-\mu}{\mu-\lambda}}|v^{n}|^{\cal K},x\rng,\qquad \forall ~ x\in{\cal K}. 
\end{equation}

To procced with the proof we  note that  \eqref{eq:npav} implies that  $|v^n|^{\cal K}\in{\cal K}+{\cal K}^*$ and,  owing to  ${\cal K} \subseteq{\cal K}^*$, we conclude that   $|v^n|^{\cal K}\in {\cal K}^*$. Combing \eqref{eq:proj_com} we have
\[
	0 \leq  \sqrt{\frac{\eta-\mu}{\mu-\lambda}} \lng |v^n|^{\cal K},x\rng\leq \langle v^1, x\rangle, \qquad \forall ~ x\in{\cal K}.
\]

Hence,  for any   $x\in{\cal K}$, the last inequality  yields 
\begin{align}
	\frac{\eta-\mu}{\mu-\lambda}\langle v^n, x\rangle^2 -\langle v^1, x\rangle^2 &\leq \frac{\eta-\mu}{\mu-\lambda}\left[\lng v^n,x\rng^2-\lng|v^n|^{\cal K},x\rng^2\right]\nonumber \\
	&=\frac{\eta-\mu}{\mu-\lambda}\lng v^n+|v^n|^{\cal K},x\rng\lng v^n-|v^n|^{\cal K},x\rng.\label{eq:ie1}
\end{align}

On the othet hand, by using $|{v^{n}}|^{\cal K}={\rm P}_{\cal K}(v^{n})+{\rm P}_{{\cal K}^*}(-v^{n})$,  ${v^{n}}={\rm P}_{\cal K}({v^{n}})-{\rm P}_{{\cal
K}^*}(-{v^{n}})$, $P_{\cal K}(v^n)\in K\subseteq K^*$, we obtain $\lng v^n+|v^n|^{\cal K},x\rng\lng v^n-|v^n|^{\cal K},x\rng=-4\lng P_{\cal K}(v^n),x\rng\lng P_{{\cal K}^*}(-v^n),x\rng\leq 0$, for any   $x\in{\cal K}$. Thus, due to $\lambda<\mu<\eta$,  the previous inequality together  \eqref{eq:ie1}  implies 
\begin{equation}\label{eq:ie2}
	\frac{\eta-\mu}{\mu-\lambda}\langle v^n, x\rangle^2 -\langle v^1, x\rangle^2 \leq 0, \qquad \forall ~ x\in{\cal K}.
\end{equation}
	Thus,  considering  that $\lambda<\mu$,  the combination of   \eqref{eq:ceq-2} with   \eqref{eq:ie2}, implies that  
\[
	\langle A x,x \rangle-\mu\|x\|^2 \leq 0, \quad \forall x\in {\cal K},
\]
which means	$\mu {\rm I_n}-A$ is ${\cal K}$-copositive.  Taking into account that  $|v^n|^{\cal K}\in {\cal K}^*$, \eqref{eq:ciqc-2} implies $v^1\in{\cal K}^*$.
	Therefore, we can apply the item (iii) of Theorem~\ref{th:best}   to conclude that $q_A$ is  spherically quasi-convex. 
 \hfill$\square$\par\end{proof}
The following example satisfies the assumptions of Theorem~\ref{th:countergg}.
\begin{example} \label{ex:countergg}
Letting  $ {\cal K}=\R^{n}_{+}$  and   $\lambda<(\lambda+\eta)/2<\mu<\eta$, the unit vectors    $v^1=(e^1+e^n)/\sqrt{2}, v^2=e^2, ~ \ldots,
~ v^{n-1}=e^{n-1}, v^n=(e^1-e^n)/\sqrt2$ are pairwise orthogonal  and  satisfy the condition  \eqref{eq:ciqc-2}. Now,  taking $\cal{K} = \cal{L}$
and denoting  $v^n = \lf((v^n)_1, (v^n)^2\rg)$, by using Lemma~\ref{l:projude}, condition \eqref{eq:ciqc-2} can be written as
\footnotesize
\[
		v^1-\sqrt{\frac{\eta-\mu}{\mu-\lambda}}\f1{\|(v^n)^2\|}\Big{(} \max\left(|(v^n)_1|,\|{(v^n)^2}\|\right)\|(v^n)^2\|,~\min\left(|(v^n)_1|,\|{(v^n)^2}\|\right)\sgn((v^n)_1)(v^n)^2\Big{)}\in {\cal K},
\]
\normalsize
and $\lambda<\mu<\eta$.
 The vectors  $v^1=(e^1+e^n)/\sqrt{2}, v^2=e^2, ~ \ldots, ~ v^{n-1}=e^{n-1}, v^n=(-e^1+e^n)/\sqrt2$  are pairwise orthogonal  and   satisfy the last inclusion.
\end{example}
\begin{theorem} \label{th:CountExMulEngfm}
Let $n\geq  3$,     $A=A\tp\in\R^{n\times n}$ and $\{v^1, v^2,\dots,v^n\}$ be an orthonormal system of eigenvectors of $A$  corresponding to the
eigenvalues $\lambda_1 , \lambda_2 ,  \ldots , \lambda_n$, respectively, such that $v^1\in\inte({\cal K}^*)$.  Let
\[
	\alpha:=\min\{ \langle v^1, y\rangle^2:~y\in \SP^n\cap{\cal K}\}>0,
\]
\begin{equation}\label{ineq:alpha}
	\eta:=\max\lf\{\frac{\sum_{i=3}^n\langle v^i,y\rangle^2}{\langle v^1,y\rangle^2}:~y\in S^n\cap{\cal K}\rg\}>0.
\end{equation}
Assume that
\begin{equation} \label{eq:CondEngfm-2}
 \qquad \lambda_1<\lambda_2\leq \cdots \leq \lambda_n \leq\lambda_2+\delta(\lambda_2-\lambda_1), \qquad \delta\in\{\alpha,1/\eta\}.
\end{equation}
Then,    $\lambda_{2} {\rm I_n}-A$ is ${\cal K}$-copositive. Consequently,  the quadratic function $q_A$ is spherically quasi-convex.
\end{theorem}
\begin{proof}
Note that the spectral decomposition of $A$ implies $A=\sum_{i=1}^n\lambda_i v^i(v^i)\tp$. Thus, considering that  $\|x\|^2=\sum_{i=1}^n\lng v^i,x\rng^2$, for any   $x\in{\cal K}$,  we conclude that
\begin{equation} \label{eq:equivnormA}
  \langle A x,x \rangle-\lambda_{2}\|x\|^2= \sum_{i=1}^n(\lambda_i-\lambda_{2})\langle v^{i}, x\rangle^2.
\end{equation}
	Since  \eqref{eq:CondEngfm-2}   implies   $\lambda_{2}-\lambda_1>0$ and   $0\leq \lambda_{j}-\lambda_{2}\leq \lambda_n-\lambda_{2}$,  for any $j\in\{3, \ldots, n\}$, it follows   from  \eqref{eq:equivnormA}  that 
\begin{equation}\label{ineq:coprf}
	  \langle A x,x \rangle-\lambda_{2}\|x\|^2 \leq (\lambda_{2}-\lambda_1)\left[ \frac{\lambda_n-\lambda_{2}}{\lambda_{2}-\lambda_1} \sum_{i=3}^n\langle  v^{i}, x\rangle^2 -\langle v^1, x\rangle^2 \right].
\end{equation}

Since \eqref{ineq:alpha} implies  $ \sum_{i=3}^n\langle  v^{i}, x\rangle^2\leq \eta \langle v^1, x\rangle^2$, the inequality \eqref{ineq:coprf} becomes
\begin{equation} \label{eq:coprfm-2}
 \langle A x,x \rangle-\lambda_{2}\|x\|^2 \leq (\lambda_2-\lambda_1)  \left[\lf(\eta\frac{\lambda_n-\lambda_{2}}{\lambda_{2}-\lambda_1}-1\rg)\langle v^1,x\rangle^2\right].
\end{equation}
First we suppose that $\delta=1/\eta$. Thus,  the last inequality in  \eqref{eq:CondEngfm-2} implies $\eta (\lambda_n-\lambda_{2})/(\lambda_{2}-\lambda_{1})\leq 1$, which combined with     \eqref{eq:coprfm-2} yields 
\begin{equation} \label{eq:cpla}
  \langle A x,x \rangle-\lambda_{2}\|x\|^2 \leq 0, \quad \forall x\in{\cal K}.
\end{equation} 
Next, suppose that $\delta=\alpha$. First of all,  noting that for any $y\in S^n$ we have $\sum_{i=3}^n\lng v^i,y\rng^2\leq \sum_{i=1}^n\lng v^i,y\rng^2=\|y\|^2=1$. Thus, using \eqref{ineq:alpha},  we conclude that 
\[
	\eta=\max\lf\{\frac{\sum_{i=3}^n\langle v^i,y\rangle^2}{\langle v^1,y\rangle^2}:~ y\in S^n\cap K\rg\}
	\le \max\lf\{\frac1{\langle v^1,y\rangle^2}:y\in S^n\cap K\rg\}=\frac1\alpha.
\]
 Hence,   it follows from  \eqref{eq:coprfm-2}  that 
\begin{equation} \label{eq:sccp}
  \langle A x,x \rangle-\lambda_{2}\|x\|^2 \leq  (\lambda_2-\lambda_1)\left[\lf(\frac1\alpha\frac{\lambda_n-\lambda_{2}}{\lambda_{2}-\lambda_1}-1\rg)\langle v^1,x\rangle^2\right].
  \end{equation} 
Due  to  $\delta=\alpha$, the last inequality in  \eqref{eq:CondEngfm-2}  implies $ (\lambda_n-\lambda_{2})/[\alpha (\lambda_{2}-\lambda_{1})]\leq 1$, which together with \eqref{eq:sccp} also implies  \eqref{eq:cpla}.  Hence, we conclude  that
$\lambda_{2} {\rm I_n}-A$ is ${\cal K}$-copositive. Therefore, since $v^1\in {\cal K}^*$ and it  is an eigenvector  of
$A$ corresponding to the eigenvalue $\lambda_1$,  by applying item (iii) of Theorem~\ref{th:best}, we can conclude that the function $q_A$ is
spherically quasi-convex.  
 \hfill$\square$\par\end{proof}
 
In the following we present an example   satisfying the assumptions of Theorem~\ref{th:CountExMulEngfm}.
\begin{example}\label{Ex:CountExMulEngfm-2}
 Let ${\cal L}$ be the Lorentz cone,  $v^i=e^i$, for any $i\in\{1, \ldots, n\}$, and  $\lambda_1<\lambda_2\leq \ldots \leq
\lambda_n<\lambda_2 +(1/2)(\lambda_2-\lambda_1)$ satisfy condition \eqref{eq:CondEngfm-2}. Note that in this case $\alpha=1/2$.
\end{example}

\begin{theorem} \label{th:TwoEng}
Let $n\ge 3$ and   $A=A\tp\in\R^{n\times n}$. Suppose that $A$ has only two distinct eigenvalues, and the smaller one has multiplicity one.
If there exists an eigenvector of $A$ corresponding  to the smaller eigenvalue belonging to ${\cal K}^*$, then  $q_{A}$ is spherically
quasi-convex.
\end{theorem}
\begin{proof}
Let  $\{v^1, v^2,\dots,v^n\}$ be an orthonormal system of eigenvectors of $A$  corresponding to the eigenvalues $\lambda_1 , \lambda_2 ,  \ldots ,
\lambda_n$ , respectively.  Without loss of generality, we assume that $\lambda_1=:\lambda <\mu:=\lambda_2=\dots=\lambda_n$ and
$v^1 \in {\cal K}^*$. Thus,  using the spectral decomposition of $A$, we have
\begin{equation} \label{eq:TwoEng}
	A= \lambda v^1(v^1)\tp+\sum_{j=2}^n\mu v^j(v^j)\tp.
\end{equation}
Since  $\|x\|^2=\sum_{i=1}^n\lng v^i,x\rng^2$,    for any   $x\in \R^{n}$,  by using \eqref{eq:TwoEng}  and $\lambda <\mu$, we  obtain that
\begin{equation} \label{eq:fetw}
\mu\|x\|^2-\langle A x,x \rangle= (\mu-\lambda)\langle v^{1}, x\rangle^2\geq 0, \qquad \forall x\in \R^{n}.
\end{equation} In particular,   \eqref{eq:fetw} implies that $\mu I_n-A$ is ${\cal K}$-copositive. Thus, since $ v^1 \in
{\cal K}^*$, by applying item (iii) of Theorem~\ref{th:best} with $\lambda_2=\mu$ we can conclude that the function $q_A$ is spherically
quasi-convex. 
 \hfill$\square$\par\end{proof}
 
In the next example we show how to generate  matrices  satisfying the assumptions of Theorem~\ref{th:TwoEng}
and consequently generate  spherically quasi-convex functions on  spherically subdual convex sets.

\bigskip \begin{example}
	The Householder matrix associated to  $v\in\inte({\cal K}^*)$ is defined by $H:={\rm I_n}-2vv^{T}/\|v\|^2$. We know that   $H$ is a symmetric and nonsingular matrix. Furthermore,  $Hv=-v$ and  $Hu=u$ for any $u\in {\cal S}$, where  $ {\cal S}:=\{ u\in \R^n~:~ \langle v, u \rangle=0\}$. It is easy to verify that the
dimension of ${\cal S}$ is $n-1$, then  we have  that $1$ and $-1$ are  eigenvalues of $H$  with multiplicities $n-1$ and $1$, respectively.
Moreover,  considering that $v\in\inte({\cal K}^*)$, Theorem~\ref{th:TwoEng} implies that   $q_{H}(x)=\langle Hx,x\rangle$ is spherically
quasi-convex.
\end{example}\index{Householder matrix}

\section{Spherically Quasi-Convex Quadratic Functions on the Spherical Lorentz Convex Set} \label{sec:qcqflc}

In this section we present a  condition partially characterising the  spherical   quasi-convexity of  quadratic functions on  spherically convex sets associated to the second order cone (Lorentz  cone). We remark that for the second order cone ${\cal L}$, since by Lemma \ref{lorcop}, we have a characterisation of ${\cal L}$-copositive matrices. By using the idea of item (iii) of Theorem~\ref{th:best},we can provide a more general result than Theorem \ref{th:CountExMulEngfm}:
\begin{theorem}\label{th:lor-sqc}
	Let ${\cal L}$ be the second order cone, $n\geq  2$,  $A=A\tp\in\R^{n\times n}$,  $\lambda_1\le\lambda_2\leq\dots\leq\lambda_n$ be  the
	eigenvalues  of $A$, $v^1$ be   an eigenvector  of $A$ corresponding to   $\lambda_1$ and $J=\diag(1,-1,\dots,-1)\in\R^{n\times n}$.
	If $v^1\in {\cal L}$ and there exists an $\rho\ge0$ such that $\lambda_2I_n-A-\rho J$ is positive semidefinite, then $q_A$ is
	spherically quasi-convex.
\end{theorem}
\begin{proof}
If  there exists an $\rho\ge0$ such that $\lambda_2I_n-A-\rho J$ is positive semidefinite, then it follows from Lemma~\ref{lorcop}  that $\lambda_2I_n-A$ is a ${\cal L}$-copositive matrix. Therefore,  considering that  $v^1\in {\cal L}={\cal L}^*$ and it  is an eigenvector  of
$A$ corresponding to the eigenvalue $\lambda_1$,  by applying item (iii) of Theorem~\ref{th:best}, we  conclude that $q_A$ is spherically quasi-convex.  
\hfill$\square$\par\end{proof}

 The next result  is a version  of Theorem \ref{thm:2neg-v} for the second order cone.
\begin{theorem}\label{thm:2neg-v2}
Let $n\ge 3$ and  $A=A\tp\in\R^{n\times n}$ . Assume that  $A$ has only two distinct eigenvalues and the smaller one has multiplicity one.
	Then, $q_{A}$ is a spherically quasi-convex function if and only if the second order cone ${\cal L}$ contains an eigenvector of $A$ corresponding to the smaller eigenvalue.
\end{theorem}
\begin{proof}
If  there exists an eigenvector  of $A$ corresponding  to the smaller eigenvalue belonging to   ${\cal L}$, then Theorem~\ref{th:TwoEng} implies
that $q_{A}$ is spherically quasi-convex.  Conversely, suppose that $q_{A}$ is spherically quasi-convex.  Let    $\lambda_1, \lambda_2,\dots,
\lambda_n$  be the  eigenvalues of $A$  corresponding to an orthonormal system
of eigenvectors  $\{v^1, v^2,\dots,v^n\}$, respectively.  Then, without loss of generality, we assume that $\lambda_1=:\lambda <\mu:=\lambda_2=\dots=\lambda_n. $
	Thus, by using the spectral decomposition of $A$, we have
\begin{equation} \label{eq:VLVT-2}
	A= \lambda v^1(v^1)\tp+\sum_{j=2}^n\mu v^j(v^j)\tp.
\end{equation}
We can also assume,  without loss of generality,  that  $v^1_ 1\geq 0$.
Let $x\in \p{\cal L}\setminus\{0\}$ and  note that    $y=2x_1e^1-x\in \p{\cal L}\setminus\{0\}$.  Since $\sum_{i=1}^nv^i(v^i)\tp=I_n$ (i.e., the spectral
decomposition of $I_n$) and $\langle x, y\rangle=0$, \eqref{eq:VLVT-2} implies that
\begin{equation}\label{eq:l}
	\lng Ax,y\rng=\lf\lng\lf[\mu\sum_{i=1}^nv^i(v^i)\tp+(\lambda-\mu)(v^1)(v^1)\tp\rg] x,y\rg\rng=(\lambda-\mu) \lng v^1,x\rng\lng v^1,y\rng.
 \end{equation}
 Since  $x, y\in {\cal L}$,  $\langle x, y\rangle=0$ and ${\cal L}$ is a self-dual cone, it follows  from  Corollary~\ref{pr:qpcf} that  $ \lng Ax,y\rng\leq 0$.
 Thus, considering that $\lambda <\mu$ and  $y=2x_1e^1-x$, equation \eqref{eq:l} yields
\begin{equation}\label{ineq:l}
	0\leq  \lng v^1,x\rng\lng v^1,y\rng=\lng v^1,x\rng [2v^1_1x_1-\lng v^1,x\rng].
\end{equation}
On the other hand, due to  $x\in{\cal L}$, we have $x^1\ge0$. Thus,   since  $v_1^1\ge0$, if $\lng v^1,x\rng<0$, then $\lng
v^1,x\rng[(2v^1_1x_1-\lng v^1,x\rng]<0$, which
contradicts \eqref{ineq:l}. Hence $\lng v^1,x\rng\ge0$, where $x$ can be chosen arbitrarily in $\p{\cal L}\setminus\{0\}$. Therefore, we proved that $v^1\in {\cal L}$. 
 \hfill$\square$\par\end{proof}

\section{Conclusions and comments}
In this chapter, we present several conditions characterising the spherical quasi-convexity of quadratic functions. In Section \ref{sec:qcqfcs}, we started from finding conditions characterising quadratic spherically quasi-convex function on a general spherically convex set. Then we move forward to the researches about the properties and  the conditions implying spherical quasi-convexity of quadratic function defined on the spherical positive orthant. Minimising a quadratic function on the spherical nonnegative orthant is of particular interest because the nonnegativity of the minimum value is equivalent to the copositivity of the corresponding matrix \cite[Proposition~1.3]{UrrutySeeger2010} and to the nonnegativity of its Pareto eigenvalues \cite[Theorem~4.3]{UrrutySeeger2010}.

In Section \ref{sec:qcqfsdcs}, we extended our results obtained in Section \ref{sec:qcqfcs} to the general subdual convex sets. We studied the sufficient conditions for spherical quasi-convex functions on spherically subdual convex sets. Also, we proved a condition for the partial characterisation of spherical quasi-convexity on spherical second order sets (Lorentz sets) in Section \ref{sec:qcqflc}. 

There are still many interesting questions to be answered in this topic: 

\begin{enumerate}
\item First of all, we note that  for all our classes of  spherically quasi-convex quadratic  functions $q_A$ on
the spherically subdual convex set  ${\cal C}=\SP^{n-1}\cap\inte({\cal K})$, the matrix $A$ has  the smallest eigenvalue  with  multiplicity one and the associated eigenvector belongs  to the dual ${\cal K}^*$ of the subdual cone ${\cal K}$. We conjecture that this condition is necessary and sufficient to characterise spherically quasi-convex quadratic  functions. 

\item We also remark  that, in Theorem~\ref{thm:2neg-v2} we present a partial characterisations of spherically
quasi-convex quadratic functions on the  spherical Lorentz convex set.  However, the general question remains open even for this specific set.

\item An even more challenging problem is to develop  efficient algorithms for constrained quadratic optimisation problems on spherically convex sets. Minimising a quadratic function on the intersection of the second order cone with the sphere is a particularly relevant, related topic, since the nonnegativity of the minimum value is equivalent to the Lorentz-copositivity of the corresponding matrix, see \cite{LoewySchneider1975, GajardoSeeger2013}.  

\item In general, replacing the second order cone with an arbitrary closed convex cone $K$ leads to the more general concept of ${\cal K}$-copositivity. By considering the intrinsic geometrical properties of the sphere, interesting perspectives for detecting the  general copositivity of matrices emerge.
\end{enumerate}

{
\chapter{Final remarks} \label{cht:FinalRemarks}

In this thesis, we studied the complementarity and it related problems. We now flash through the results of this thesis, summarise our contributions, and present possible future works. 

\section{Summary of results}
We start this study from solving a linear complementarity problem on an extended second order cone. In Chapter \ref{cht:lesoc}, we convert an ESOCLCP to a MixCP on nonnegative orthant, therefore we can exploit complementarity function methods to solve the original problem. Then we provide two approaches: reformulate to a system of nonlinear equations; and, reformulate to an unconstrained minimisation problems.  Based on these two approaches, algorithms are provided for solving the problem. It must be recognised that such conversion of a linear complementarity problem to a mixed complementarity problem may increase the complexity of the original problem. Finally, we furnish this chapter by providing a numerical example. 

In Chapter \ref{cht:stesoc}, we studied the stochastic linear complementarity problems on extended second order cones. Similarly, we convert the stochastic ESOCLCP to a corresponding stochastic mixed complementarity problem on nonnegative orthant. We introduced the Conditional Value-at-Risk method to measure the loss of complementarity in the stochastic case. Unlike the ESOCLCP discussed in Chapter \ref{cht:lesoc}, we can only use the unconstrained minimisation approach to solve this problem. We also introduce an algorithm and provide a numerical example for this problem. 

In Chapter \ref{cht:psp}, we introduced the mean-Euclidean  norm (MEN) model for portfolio optimisation problem. This new model is based on the mean-absolute deviation (MAD) model. The KKT condition of the MEN model is a nonlinear complementarity problem on extended second order cone. Based on the results obtained in Chapter \ref{cht:basic_esoc} and Chapter \ref{cht:lesoc}, we find an analytical solution to the MEN model.

Chapter \ref{cht:sphere} presents several conditions characterising the spherical quasi-convexity of quadratic functions. We started from finding conditions characterising quadratic spherically quasi-convex function on a general spherically convex set. Then we move forward to the researches about the properties and  the conditions implying spherical quasi-convexity of quadratic function defined on the spherical positive orthant. Further, we extended our results on to the general subdual convex sets. The sufficient conditions for spherical quasi-convex functions on spherically subdual convex sets are presented. In addition, a condition for the partial characterisation of spherical quasi-convexity on spherical  Lorentz sets is provided.

\section{Contributions}

We now outline our main contributions: 

\begin{enumerate}
	\item The study about the linear complementarity problem on extended order cone. We successfully find an solution to this problem. Algorithms and numerical examples are provided. These results are published on our paper \cite{nemeth2018linear}.

	\item The study about the stochastic linear complementarity problem on extended order cone. The CVaR-based minimisation problem is used on the   merit function to measure the loss of complementarity in the stochastic case. Since the objective function of the CVaR-based minimisation problem is convex in some neighbourhoods (it is not globally convex), given an initial point that close enough to the optimal point, we can still solve this problem by an algorithm. Hence, we introduced an algorithm for solving this problem and provide a numerical example. 
	
	\item The introduction of the mean-Euclidean  norm (MEN) model. Since the KKT condition of this model is a nonlinear complementarity problem on extended second order cone, we innovatively used the results in Chapter \ref{cht:basic_esoc} and Chapter \ref{cht:lesoc} to find the analytical solution to this model. 
	
	\item Chapter \ref{cht:sphere} is a continuation of \cite{FerreiraIusemNemeth2014, FerreiraIusemNemeth2013, FerreiraNemeth2017}, where intrinsic properties of the spherically quasi-convex sets and functions were studied. As far as we know this is the pioneering study of spherically quasi-convex quadratic functions on spherically subdual convex sets. The results about the conditions characterising spherical quasi-convexity of quadratic function defined on the spherical positive orthant  are published on \cite{FerreiraNemethXiao2018}. 

\end{enumerate}

}

\backmatter

\printglossary[type=\acronymtype,title=Acronyms,style=long3colheader]


\nomenclature[S]{$\R^n$}{the $n$ dimensional Euclidian space}
\nomenclature[V]{$\langle \cdot, \cdot \rangle$}{the inner product of two vectors}
\nomenclature[V]{$\|\cdot\|$}{the norm in Euclidian space}
\nomenclature[S]{$\R^n_+$}{the nonnegative orthant}
\nomenclature[S]{$\R^n_{++}$}{the positive orthant}
\nomenclature[V]{$\bs x\\u\es$}{the pairs of vectors of $x$ and $u$}
\nomenclature[Se]{${\cal K} $}{a cone, a convex cone, or a convex set}
\nomenclature[M]{$S_0$}{the $S_0$ matrix}
\nomenclature[P]{$\SCP(F,{\cal K})$}{the solution set of $\CP(F,{\cal K})$}
\nomenclature[f]{$F(x)$}{the linear function $F(x)=Tx+r: \R^n\rightarrow \R^n$ }
\nomenclature[P]{$\ICP(G,F,{\cal K})$}{the implicit complementarity problem defined by $G$, $F$, and ${\cal K}$}
\nomenclature[P]{$\MixCP(F_1,F_2,{\cal S})$}{the mixed complementarity problem defined by function $F_1$, $F_2$, and set ${\cal S}$}
\nomenclature[P]{$MixICP(F_1,F_2,G_1,{\cal S})$}{the mixed implicit complementarity problem defined by function $F_1$, $F_2$, and $G_1$, as well as the set ${\cal S}$}
\nomenclature[P]{$\SMixICP(F_1,F_2,G_1,{\cal S})$}{the solution set of the mixed complementarity problem $\MixICP$ $(F_1,F_2,G_1,{\cal S})$}
\nomenclature[P]{$\SMixCP(F_1,F_2,{\cal S})$}{the solution set of $\MixCP(F_1,F_2,{\cal S})$}
\nomenclature[P]{$\SICP(G,F,{\cal K})$}{the solution set of $\ICP(G,F,{\cal K})$}
\nomenclature[Se]{$\C({\cal K})$}{the complementarity set of ${\cal K}$}
\nomenclature[Se]{$L(k,\ell)$, $L$}{the extended second order cone}
\nomenclature[Se]{$M(k,\ell)$, $M$}{the dual cone of extended second order cone}
\nomenclature[M]{$T=\bs A & B\\C & D \es$}{the T matrix with component matrices $A$, $B$, $C$, and $D$}
\nomenclature[Se]{${\cal N}_{\overline{x}}$}{the neighbourhood of ${\overline{x}}$}
\nomenclature[Se]{${\cal K} ^*$}{the dual cone of ${\cal K}$}
\nomenclature[M]{$I_n$}{the $n\times n$ identity matrix}
\nomenclature[P]{$\CP(F,{\cal K})$}{the complementarity problem defined by function $F$ and cone ${\cal K}$}
\nomenclature[P]{$\SLCP(T,r,{\cal K})$}{the solution set of $\LCP(T,r,{\cal K})$ }
\nomenclature[F]{$\partial\mathbb{F}^{\MixCP}_{FB}(x,u,t)$}{the generalised Jacobian set of $\mathbb{F}^{\MixCP}_{FB}(x,u,t)$ }
\nomenclature[M]{$ \mathcal{D}_a$}{denoted by diag$\big(a_1(x,u,t)$, $\dots$, $a_k(x,u,t)\big)$ the $k\times k$ diagonal matrices }
\nomenclature[M]{$ \mathcal{D}_b$}{denoted by diag$\big(b_1(x,u,t)$, $\dots$, $b_k(x,u,t)\big)$ the $k\times k$ diagonal matrices }
\nomenclature[F]{$\theta^{\MixCP}_{FB}(x,u,t)$}{the associated merit function based on FB C-function}
\nomenclature[Se]{$\mathcal{C}$}{the complementarity index set}
\nomenclature[Se]{$\mathcal{R}$}{the residual index set}
\nomenclature[Se]{$\mathcal{P}$}{the positive index set}
\nomenclature[Se]{$\mathcal{N}$}{the negative index set}

\nomenclature[F]{$JF$}{$\equiv\frac{\partial F_i}{\partial x_j}$ the $m\times n$ Jacobian of a mapping $F:~\R^n\rightarrow \R^m$}
\nomenclature[F]{$J_xF(x,y)$}{$:=\lf(\frac{\partial F_i}{\partial x_j}\rg)$ the $m\times n$ partial Jacobian of a mapping $F:~\R^k\times\R^{\ell}\rightarrow \R^m$}
\nomenclature[F]{${\cal r}\theta$}{$:=\lf(\frac{\partial \theta}{\partial x_j}\rg)$ the gradient of a function $\theta:~\R^n \rightarrow \R$}
\nomenclature[F]{$\mathbb{F}^{\MixCP}_{FB}(x,u,t)$}{FB-based C-function of MixCP problem}
\nomenclature[F]{$\psi_{FB}(a,b)$}{scalar form of FB C-function }
\nomenclature[M]{$\mathcal{A}$}{$\partial\mathbb{F}^{\MixCP}_{FB}(x,u,t)$ an element in the generalised Jacobian set of \\ $\mathbb{F}^{\MixCP}_{FB}(x,u,t)$}
\nomenclature[Se]{$\in,~\not\in$}{element membership, non-membership in a set}
\nomenclature[Se]{$\emptyset$}{the empty set}
\nomenclature[Se]{$\subseteq,~\subset$}{subset inclusion, subset}
\nomenclature[Se]{$\cup,~\cap$}{set union, intersection}
\nomenclature[Se]{$S_1\setminus S_2$}{set $S_1$ minus set $S_2$}
\nomenclature[S]{$\mathbb{V}^{\perp}$}{the orthogonal complement of a subspace $\mathbb{V}$ in $\R^n$}
\nomenclature[Se]{$\Omega$}{the sample set of possible outcomes}
\nomenclature[Se]{$\mathcal{F}$}{$\subseteq 2^{\Omega}$, a $\sigma$-algebra generated by $\Omega$ (or a collection of all subsets of $\Omega$)}
\nomenclature[F]{$\mathcal{P}$}{$: \mathcal{F}\rightarrow [0,1]$, a function maps from events to probabilities.}
\nomenclature[V]{$\omega$}{$\in\Omega$ an $n$-dimensional random vector}
\nomenclature[M]{$A$}{$:=(a_{ij})$ a matrix with entries $a_{ij}$}
\nomenclature[M]{$A\tp$}{the transpose of a matrix $A$}
\nomenclature[M]{$A^{-1}$}{the inverse of a matrix $A$}
\nomenclature[M]{$M\/A$}{the Schur complement of $A$ in $M$}
\nomenclature[F]{$\E[\cdot]$}{the expected value of the random vector}
\nomenclature[F]{$\Phi$}{$:\R^n \times\Omega\rightarrow\R^n$ the multi dimensional C-function}
\nomenclature[F]{$\phi$}{$: \R\times\R \rightarrow\R$ a scalar C-function such that $\phi(a,b) = 0 ~\Leftrightarrow~ a\ge 0,~b\ge 0, ~ab = 0$}
\nomenclature[V]{$x\perp y$}{$x$ and $y$ are perpendicular}
\nomenclature[V]{$x \circ y$}{$:=(x_i y_i)^n_{i=1}$ the Hadamard product of $x$ and $y$}
\nomenclature[F]{$\mathbf{1}_{S}(x)$}{$: S\subseteq \R^n\rightarrow\{0,1\}$ the indicator function with $\mathbf{1}_{S}(x)=1$ if $x\in S$ and $\mathbf{1}_{S}(x)=0$ otherwise}
\nomenclature[Se]{$\arg \min_{a}F(x)$}{the set of constrained minimisers of $F$ on $a$}
\nomenclature[F]{$[t]_+$}{$:= \max\{0,t\}$ the nonnegative operator}
\nomenclature[F]{$p(\cdot,\mu)$}{a smoothing function of the nonnegative operator}
\nomenclature[F]{$i(x)$}{$:\R^n\setminus\{0\} \rightarrow \R^n\setminus \{0\}$; $i(x)=\frac{x}{\|x\|^2}$ the inversion operator}
\nomenclature[F]{$\mathcal{I}(F)$}{$:{\cal K}\rightarrow \R^n$ the inversion (of pole 0) of the mapping $F: \R^n\rightarrow \R^n$}
\nomenclature[F]{$\underline{F}^{\#}(x_0, {\cal K})$}{the lower scalar derivative of the mapping $F: {\cal K} \rightarrow \R^n$ at $x_0$}
\nomenclature[M]{$\lambda_{max}(A)$}{the largest eigenvalue of a matrix $A$}
\nomenclature[S]{$T_x\SP^{n-1}$}{$:=\left\{v\in \R^n :~ \langle x, v \rangle=0, ~ x\in\SP^{n-1} \right\}$ the tangent hyperplane at point $x\in\SP^{n-1}$}
\nomenclature[F]{$\gamma$}{$:[x,y]\rightarrow\SP^{n-1}$ the geodesic}
\nomenclature[F]{$\gamma_{xy}(t)$}{$:[0,1]\rightarrow\SP^{n-1}$ the minimal geodesic from $x$ to $y$}
\nomenclature[F]{$\grad f(x)$}{the gradient on the sphere of a differentiable function \(f: \Omega \to  \R\) at a point \(x\in \Omega\)}
\nomenclature[F]{$Df(x)$}{the gradient a differentiable function \(f: \Omega \to  \R\) at a point $x$}
\nomenclature[Se]{${\cal K}_{\cal S}$}{a cone spanned by ${\cal S}$}
\nomenclature[Se]{$\textbf{B}(x,\delta)$}{a spherically open ball center in point $x$}
\nomenclature[Se]{$\bar{\textbf{B}}(x,\delta)$}{a spherically closed ball center in point $x$}
\nomenclature[Se]{$[f\leq c]$}{the sub-level set of a function \(f:~{\cal S} \to \R\)}
\nomenclature[F]{$f\circ g$}{the composition of function $f$ and $g$}
\nomenclature[F]{$q_A$}{$: {\cal S}\to\R$ the associated quadratic function defined by the symmetric matrix $A=A^T\in\R^{n\times n}$}
\nomenclature[F]{$\varphi_A$}{$ : \inte({\cal K}) \to \R$ the Rayleigh quotient function restricted on cone $\inte({\cal K}) $ defined by matrix $A$}
\nomenclature[Se]{$\inte {\cal S}$}{the interior of a set ${\cal S}$}
\nomenclature[Se]{${\cal L}_c$}{$:=\left\{x\in\R^n:~\lng v^1,x\rng\geq \sqrt{\sum_{i=2}^{n} \theta_i(c) \lng v^i,x\rng^2}\right\}$ the convex cone defined by $c$}

\printnomenclature[2in]

\printindex

\bibliographystyle{plain}
\bibliography{thesiseg}
\end{document}